# d-Complete Posets Generalize Young Diagrams for the Jeu de Taquin Property

Robert A. Proctor*
Department of Mathematics
University of North Carolina
Chapel Hill, NC 27599
rap[at]email:unc:edu

## Abstract

The jeu de taquin process produced a standard Young tableau from a skew standard Young tableau by shifting its entries to the northwest. We generalize this process to posets: certain partial numberings of any poset are shifted upward. A poset is said to have the jeu de taquin property if the numberings resulting from this process do not depend upon certain choices made during the process. Young diagrams are the posets which underlie standard Young tableaux. These posets have the jeu de taquin property. d-Complete posets are posets which satisfy certain local structual conditions. They are mutual generalizations of Young diagrams, shifted Young diagrams, and rooted trees. We prove that all d-complete posets have the jeu de taquin property. The proof shows that each d-complete poset actually has the stronger "simultaneous" property; this may lead to an algebraic understanding of the main result. A partial converse is stated: "Non-overlapping" simultaneous posets are d-complete.

* Supported in part by NSA grant MDA 904-95-H-1018.

1991 *Mathematics Subject Classifications:* Primary 06A07, 05C85, 05E10; Secondary 20F55, 17B67.

March 12, 2003
Page 11 revised 10/30/04, 8/7/06, 7/15/15

## 1. Introduction

Please consider the following problem in structural combinatorics:

**Upper Echelon Replacement Problem.** *The organizational chart for a corporation has $n$ positions which are initially filled by $n$ employees with distinct seniorities of $1, 2, \ldots, n$ years. Any position may be subordinate to more than one other position: the chart is the (Hasse) diagram of a poset. No employee may be supervised by another employee with less seniority: If we denote the employees by their seniorities, then the initial assignment of employees to positions is an "order extension" of the poset. If a position becomes vacant, then the most senior of the employees in the immediately subordinate positions is shifted up into the vacant position. This procedure is repeated until only a minimal position in the chart is vacant. If multiple positions are simultaneously vacated, then this process is first applied to fill the position held by the most junior of the departed employees. It is then successively applied to the other vacated positions in the increasing order of the seniorities of the departed employees. This is the most natural way to fill simultaneous vacancies.*

*Now suppose that some "upper echelon" (order filter) of employees from the chart are invited to an executive retreat at a ski resort, and that an avalanche completely buries the resort. For most charts, the final assignment pattern of the uninvited survivors will depend upon the assignment pattern of the victims. In Figure 1, the hollow numbers indicate the victims. The final assignment of the employee with 1 year of seniority depends upon the former assignments of the departed employees with 2 and 4 years of seniority. (A phenomenon of a similar flavor arises when two spouses are mortally injured in an accident: then the distribution of their estate can depend upon the relative timing of their deaths.) However, a result of Schützenberger and Thomas implies that if the chart arises from a Young (Ferrers) diagram, then the final assignment pattern of the survivors surprisingly depends upon only the initial assignment pattern of the survivors. Since the former seniority pattern of the departed supervisors should be irrelevant after the disaster, we will call a chart "fair" if the final survivor assignment pattern always (for any upper echelon of victims) depends only upon the initial assignment pattern of the survivors. From the extension of the results of [Sch] [Tho] by Sagan and Worley [Sag] [Wor], it is known that charts arising from shifted Young diagrams are also fair. Can other large families of fair charts be described?*

The notion of fair chart is easily equivalent to the "jeu de taquin" (jdt) property for posets defined in Section 2: A poset has the jdt property if the "jeu de taquin emptying procedure" (introduced by Schützenberger for shapes) always yields the same survivor pattern, independently of the order in which the vacancies are "slid out". This was the property actually confirmed by the authors cited above for shapes (Young diagrams) and shifted shapes (shifted Young diagrams): the skew standard tableaux in those papers correspond to the initial survivor assignment patterns. (The jdt emptying process is only distantly related to the "promotion" and "evacuation" procedures, which were also defined by Schützenberger in terms of sliding manipulations.)

d-Complete posets are posets which first arose [Pr3] in the representation theory of Kac-Moody algebras in a manner which naturally related and mutually generalized the obvious poset structures formed by the boxes in shapes and shifted shapes. However, the resulting statement of the definition [Pr2] of d-complete poset is purely combinatorial; it is couched (Section 2) in terms of local structural conditions.

Combinatorial methods are used in this paper to prove that all d-complete posets have the jeu de taquin property. This is one of the three major combinatorial results obtained so far in which d-complete posets generalize shapes, and the first to be fully written up. Dale Peterson and this author have jointly obtained [Pr5] a hook length product formula for the number order extensions of a d-complete poset. This formula is a mutual generalization of the formulas of Frame-Robinson-Thrall and of Knuth for the numbers of standard and standard shifted Young tableaux. Thirdly, Kawanaka has announced [Kaw] a complete Sprague-Grundy analysis of a version of Nim in which each move consists of the removal of a "rim hook" from the bottom of one of k given d-complete posets, rather than removing a subpile from one of k given piles of coins. This generalizes 1950's results of Sato and Welter for shapes.

The jeu de taquin (jdt) process for shapes is closely related to the Robinson-Schensted correspondence. Some of the nicest proofs of the Littlewood-Richardson rule for the product of two Schur functions refer to the jdt process or to Schensted's algorithm. This author has recently obtained one generalization of the Littlewood-Richardson rule for each jdt poset whose order dual is also jdt. Versions or relatives of the jdt process have been used in various papers concerned with combinatorial bijections, symmetric functions, representation theory, and algebraic geometry.

Only 236 of the 14,512 connected posets with 8 elements have the jdt property, and d-complete posets account for 181 of these 236. Visually comparing the diagrams for the jdt posets to those for the non-jdt posets suggests only fragmentary conjectures predicting which posets are jdt. To see that the shape consisting of two rows of three boxes apiece is a fair chart, one must consider 50 scenarios: there are 5 possible initial assignment patterns and 10 possible upper echelons. Several shifts are required to compute the outcome of a typical scenario. Although Lie theoretic analogies heuristically led to the conjecture that d-complete posets are jdt, no algebraic explanation of this phenomenon is known.

Kimmo Eriksson gave a new proof [Eri] that shapes are jdt. Our proof proceeds by generalizing each of the two stages of his proof. We say that a poset is "simultaneous" if our generalization of the second stage of his proof works. The stronger version of our main result is: all d-complete posets are simultaneous. It is hoped that the simultaneous property will lead to an algebraic explanation of the jdt property for d-complete posets. A poset is "non-overlapping" if it satisfies the last of the three d-complete structural conditions. We have proved that a non-overlapping poset is simultaneous if and only if it is d-complete.

The process of repeatedly filling one vacancy as it shifts downward within the chart can be thought of as "sliding out" [Sc1] the "ghost" of the departed employee. The words "promotion" (for the survivors) and "evacuation" (for the ghosts) cannot be used in the jeu de taquin context since they have unrelated pre-existing meanings. By "jeu de taquin emptying procedure" we mean the sliding out of all of the ghosts in the order of their seniorities. There have been only a few papers which have used the slideout operation in the general setting of posets, and each of those seemed to be mainly concerned with the evacuation procedure [Sc2]. An early occurrence of the terminology "jeu de taquin" is in Foata's review [Foa] of [Sc3]. There it is noted that although the jdt property for shapes is not explicitly noted in [Sc3], it is a consequence of Proposition 3.7 of [Sc3]. Thomas provided an early accessible proof [Tho].

One early step in the sequence of developments which led to Littelmann's path description of Kac-Moody characters was the paper "G/P-I" in the Laksmibai-Seshadri program for describing characters of simple Lie algebras; it dealt with certain "minuscule" cases. In 1993 we proved that the combinatorial reformulation of the G/P-I basis theorem existed for certain $\mathbb{Z}$-modules based upon a colored poset if and only if that poset satisfied certain local structural conditions. That was the first appearance of the conditions which now constitute the definition of d-complete poset. That basis result [Pr3] is a fourth major result concerning d-complete posets; it is still only in preliminary manuscript form. Dale Peterson had earlier formulated the notion of $\lambda$-minuscule Weyl group element: If $\lambda$ is an integral weight for a Kac-Moody Weyl group W, then $w \in W$ is *$\lambda$-minuscule* if it has a reduced decomposition which subtracts one positive simple root at each stage as its simple reflections are successively applied to $\lambda$ and its descendants. When W is simply laced, we proved [Pr2] that the integral weights between $\lambda$ and $w\lambda$ inclusive form a distributive lattice L. That paper also showed that when $\lambda$ is further assumed to be dominant, then the poset P of meet irreducibles of L is a d-complete poset, and that all d-complete posets arise from dominant integral simply laced $\lambda$-minuscule w's in this fashion. Coming from the direction of his independently instigated study of fully commutative elements of Coxeter groups, Stembridge extended this characterization of the posets of meet irreducibles to the non-simply laced and non-dominant cases. He also realized these posets in the dominant case with certain subsets of positive roots [St2]. Much earlier, the *minuscule posets* $a_n(j)$, $d_n(n)$, $d_n(1)$, $e_6(1)$, and $e_7(1)$ arose when this process was applied to the longest $\lambda$-minuscule w for each minuscule weight $\omega_j$ for the finite simply laced Weyl groups W of rank n [Pr4]. The minuscule posets $a_{g+h+1}(g+1)$, $d_{h+3}(h+3)$, $d_5(1)$, and $e_7(1)$ appear in Figures 8.1, 8.2, 6, and 8.15 of this paper. After reading [Pr4], Peterson used a computer to confirm the jdt property for $e_6(1)$ and $e_7(1)$. That empirical result led us to perform enough computations in 1994 to conjecture that every d-complete poset is jdt.

The "double tailed diamond" minuscule posets $d_n(1) =: \Delta_{n-2,n-2}$ (Section 2) arising from the Weyl group $D_n$ play a key role in the definition of d-complete poset; this is the origin of the "d". These posets are now known to have appeared in Stanley's thesis. After generalizing the notion of partition of an integer to the notion of P-partition of an integer with respect to a fixed poset P, Stanley found analogs of a partition generating function identity of Euler when P was a shape, double tailed diamond, or rooted tree. Gansner later confirmed Stanley's conjecture that such an identity also holds for any shifted shape, and this author found [Pr4] such identities for $e_6(1)$ and $e_7(1)$. These four families of posets (shapes, shifted shapes, rooted trees, minuscules) were the precursors to the set of all d-complete posets. The "$x = 1$" specializations of these identities yielded hook length product formulas for the number of order extensions of posets of these kinds. The joint work with Peterson mentioned above actually produces such a P-partition generating function identity in x for any d-complete poset [Pr5]. Hence all d-complete posets are "hook length" posets. Ishikawa and Tagawa have announced a combinatorial proof of this result [I-T]; the original proof [Pr5] employed algebraic geometry at one step.

The paper [Pr1] classified d-complete posets with Dynkin diagrams which are contained in their order diagrams. The non-trivial building blocks of connected d-complete posets, the "irreducible components", fall into 15 classes shown in Figures 8.1 - 8.15 below. Class 1 is comprised of shapes and Class 2 is comprised of shifted shapes. Any connected d-complete poset can be formed by joining irreducible components together with rooted trees.

We extend the first stage of Eriksson's re-proof for shapes to general posets in Proposition 6.2 and strengthen it in Theorem 6.5. This general sufficient condition for the jdt property could be the main result of its own short paper; it describes what could be the most efficient way to test a particular poset P for this property. Theorem 6.5 and an (two) application(s) of the second stage of Eriksson's proof together form what may be the "book" combinatorial proof of the jdt property for shapes (shifted shapes).

Even with this reduction, there does not seem to be much hope for deducing the jdt property (which is a global condition) directly from the d-complete axioms (which are local structural conditions). However, repeated application of the d-complete requirements in [Pr1] yielded the description of the possible global structures for connected d-complete posets via the list of the possible irreducible components. Rather than repeating this local-to-global iterative work, we will instead refer to this list of possible irreducible components. In fact, the conditions of Theorem 6.5 are simpler to check for a d-complete irreducible component than for most posets. Nonetheless, the second stage of Eriksson's proof must be repeatedly applied to see that an irreducible component is jdt. We say a poset has the simultaneous property if each challenge posed by Theorem 6.5 can be answered in this specific iterative manner. We convert the question of whether the simultaneous property holds for a connected d-complete poset to a set of slightly stronger conditions for its irreducible components. Checking these conditions for the various irreducible components finishes the proof of our main result.

Since the 15 classes of irreducible components may be grouped into two batches of closely related cases, checking the necessary conditions takes only a dozen pages once some preliminary results have been established. Support for this proof of Theorem 5.1 and for our emphasis on the simultaneous property is given by the existence of the partial converse, Theorem 5.2, mentioned above. This result suggests that any proof that a d-complete poset has the simultaneous property should be fairly "tight", thereby diminishing the prospects for easily(loosely)/non-algebraically deducing the jdt property directly from the d-complete axioms. In fact this tightness can be seen when the details of the various cases are read. In each case it is initially remarkable to see how the simultaneous conditions are barely satisfied. But this is actually forced by the existence of the partial converse, since most perturbations of the structure of the poset at hand will violate the first two thirds of the d-complete definition. Comments and examples pertaining to the interplay between the d-complete property and the jdt or simultaneous properties appear at the ends of Sections 2 and 7. If these examples are kept in mind, then it will be clear that the spirit of the first two thirds of the d-complete definition is the driving force for checking the simultaneous conditions in the various classes; reference to [Pr1] for the global structure in each class is essential since the simultaneous property is global in nature.

The seventh paragraph of Section 2 outlines a method of constructing "cheap" jdt posets: a global structural "cure" makes up for arbitrarily bad local structures in the lower "halves" of these posets. This construction indicates that it should not be possible to completely characterize jdt posets with only local conditions. Hence the fact that d-complete posets account for only 181 of the 236 posets mentioned above should not be disappointing. Two of the other 55 posets appear in Figure 4.

It would be nice to have an algebraic proof that d-complete implies jdt; this might help to explain some of the mysterious roles played by the jdt procedure in Lie combinatorics. (But it is also desirable to have a purely combinatorial proof of a purely combinatorial statement!) The present proof should be viewed as a first step toward an algebraic proof, rather than as a substitute for an algebraic proof. The existence of arbitrary local structures in jdt posets indicate that the prospects for algebraizing the jdt property are poor. Replacing the jdt property by the simultaneous property as the center of attention should improve the prospects for finding algebraic insights. Once the iterative verifications of the simultaneous property in the various cases have been read, many readers may find plausible our belief that some unknown underlying algebraic mechanism should account for the remarkable satisfaction of the stringent simultaneous conditions. Moreover, every statement in this paper may be recouched as a statement regarding reduced decompositions of $\lambda$-minuscule Weyl group elements via [Pr2] or [St2]. Although this viewpoint has not led us to an algebraic proof, it does imply that the main result can be restated as a theorem in the subject of Weyl groups concerning certain manipulations of reduced decompositions. Since the uniqueness of the jdt emptying procedure can be viewed as a kind of commutativity statement, the main result of this paper gives another instance of a commutativity phenomenon arising for $\lambda$-minuscule Weyl group elements. It was already known that each $\lambda$-minuscule element is fully commutative in the sense of Stembridge [St2]. Let $g$ be a simple Lie algebra, and let $b$ be a fixed Borel subalgebra of it. Peterson showed that the abelian ideals of $b$ correspond to the $\omega_0$-minuscule elements of $g'$, where $\omega_0$ is the fundamental weight for the additional Dynkin diagram node adjoined when the affine Lie algebra $g'$ associated to $g$ is formed. d-Complete posets can be defined algebraically [St2], and non-overlapping simultaneous posets are equivalent to d-complete posets by Theorem 5.2. Hence at least non-overlapping simultaneous posets can be thought of as being objects of an algebraic nature. They are classified by Dynkin diagrams.

The definitions needed to understand the statement "d-complete implies jdt" comprise Section 2. The main results and an outline of their proofs are presented in Section 5. Browsers should refer back to the additional "first half" definitions and facts in Sections 3 and 4 as much as is needed to understand this outline. After extending and strengthening Eriksson's reduction in Section 6, it is possible to present the recursive definition of the simultaneous property in Section 7. The first paragraph of Section 8 notes that the strengthened reduction becomes relatively easy to check for the irreducible components. After attaining some familiarity with the simultaneous poset definition, browsers should next read the first two paragraphs of Section 12. Here the remarkable intricacy of the repair mechanics is described, and a bit of mysterious Lie combinatorial numerology is noted. Section 8 shows that reassembling the "slant irreducible" components of a connected d-complete poset produces a simultaneous poset, provided it is known that these components are "strongly" simultaneous. So at this point the problem has been reduced to showing that every irreducible d-complete component is strongly simultaneous. The remaining definitions needed for the second half of this paper appear in Sections 7 and 9 and in the first paragraph of Section 10. Results regarding the paths of slid-out vacancies in irregularly shaped convex grid regions of the plane are developed in Sections 9 and 10. After Section 9, no references are made to tableau entries: Shorter proofs using qualitative arguments for vacancy paths are given. (It should be possible to define "semistandard" numberings on some d-complete posets (beyond shapes and shifted shapes) which possess the jdt property. Extending our verification of the jdt property to such semistandard numberings should be facilitated by this emphasis on vacancy paths.) Sections 11 - 14 show that every irreducible component is strongly simultaneous by studying the choreography of up to n vacancy paths. Given the preliminary results of Sections 9 - 11, the proof of strongly simultaneous in Theorem 12.1 takes only a couple of sentences for Class 1 and only a page for Classes 1 - 7 combined.

## 2. Definitions of Jeu de Taquin and d-Complete Posets

Let $P$ be a poset and let $x,y \in P$. We say $x$ *is covered by* $y$ and write $x \to y$ if $x < y$ and $x < z \leq y$ implies $z = y$. The *diagram* of $P$ is the directed graph on $P$ whose edges are these covering relations. The arrowheads on the upper ends of the edges in the diagrams of Figures 1 - 6 and on the western and northern ends of the edges in Figures 7 - 8 have been suppressed. A poset is *connected* if its diagram is connected. A subset $I \subseteq P$ is an *ideal* if $y \in I$ and $x \leq y$ implies $x \in I$. A subset $F \subseteq P$ is a *filter* if $x \in F$ and $y \geq x$ implies $y \in F$. A *(shifted) shape* is a poset which can be realized as a filter of one of the posets depicted in Figure 8.1 (Figure 8.2). If $x,y \in P$, define the *interval* $[x,y] := \{z: x \leq z \leq y\}$. Elements $t_1,\ldots,t_n \in P$ form a *chain* if $t_1 \to t_2 \to \ldots \to t_n$ in P. A *rooted tree* is a poset with a unique maximal element in which each non-maximal element is covered by one other element. Let $b,n \geq 0$. Define the poset $\Delta_{b,n}$ to consist of two incomparable elements $x_0$ and $y_0$ together with a chain of b elements $a_b \to \ldots \to a_2 \to a_1$ below $x_0$ and $y_0$ and a chain of n elements $t_1 \to t_2 \to \ldots \to t_n$ above $x_0$ and $y_0$. The poset $\Delta_{3,3}$ appears on the left in Figure 2.

To define the jdt property, we must introduce analogs for posets of some kinds of Young tableaux. Let P be a poset. Let $D \subseteq P$ and set $d := |D|$. A *numbering* $\varphi$ *of* D is a bijection from D to a *chain poset* $\{C_1 \to C_2 \to \ldots \to C_d\}$ such that $x \leq y$ for $x,y \in D$ implies $\varphi(x) \leq \varphi(y)$. Let $\gamma$ and $\rho$ be numberings of some $D_G \subseteq P$ and some $D_R \subseteq P$ to the chain posets $\{G_1 \to G_2 \to \ldots \to G_g\}$ and $\{R_1 \to R_2 \to \ldots \to R_r\}$ respectively. Then $\gamma$ and $\rho$ together form a *bi-numbering* $\gamma/\rho$ *of* P if $D_G \cup D_R = P$ and $D_G \cap D_R = \varnothing$. The element $G_i$ is called the *ith bubble* and the element $R_j$ is called the *jth label.* The bi-numbering $\gamma/\rho$ should be visualized by writing the number i in green ink next to $\gamma^{-1}(G_i)$ and the number j in red ink next to $\rho^{-1}(R_j)$. Each of the first four diagrams of Figure 3 depict bi-numberings of $\Delta_{3,2}$: here outlined (bold) numbers inside the elements' circles indicate the green (red) bubble (label) subscripts. We will blur the distinctions between the bubble $G_i$ and its location $\gamma^{-1}(G_i)$ and between the label $R_j$ and its location $\rho^{-1}(R_j)$. For example, we refer to the elements of the set $\{\rho(x): x \to \gamma^{-1}(G_i)\}$ as the *labels covered by* $G_i$. Let $R(.)$ be the operation which extracts the red numbering from a bi-numbering of P, namely $R[\gamma/\rho] := \rho$.

Let $\gamma/\rho$ be a bi-numbering of P. We want to define an operator $J$ which will successively "slide out" each of the bubbles $G_1, G_2, \ldots G_g$ as far as possible. The result of applying the *move operator* $M_1$ to $\gamma/\rho =: \gamma^{(0)}/\rho^{(0)}$ is the bi-numbering $M_1(\gamma/\rho)$ obtained by interchanging $G_1$ with the largest of the labels it covers if it covers one or more labels; otherwise define $M_1(\gamma/\rho)$ to be $\gamma/\rho$. Let $\gamma/\rho$ denote the bi-numbering shown in the first diagram of Figure 3. Then the second diagram shows $M_1(\gamma/\rho)$. This operator $M_1$ is further defined similarly upon the resulting bi-numberings $\gamma^{(k+1)}/\rho^{(k+1)} := M_1(\gamma^{(k)}/\rho^{(k)})$. The result of applying the *slideout operator* $S_1$ to $\gamma/\rho$ is the bi-numbering $S_1(\gamma/\rho)$ produced by repeatedly applying $M_1$ to $\gamma/\rho$ and its successors $\gamma^{(k)}/\rho^{(k)}$ until $M_1(\gamma^{(k)}/\rho^{(k)}) = \gamma^{(k+1)}/\rho^{(k+1)}$. The third diagram in Figure 3 shows $S_1(\gamma/\rho)$. For $i \geq 2$, the actions of the *move and slideout operators* $M_i$ and $S_i$ are analogously defined, starting with the bi-numbering $S_{i-1}\circ \ldots \circ S_2 \circ S_1(\gamma/\rho)$. Define the result of applying the *emptying operator* $J$ to $\gamma/\rho$ to be $J(\gamma/\rho) := S_g\circ \ldots \circ S_2\circ S_1(\gamma/\rho)$. The fourth and fifth diagrams in Figure 3 show $J(\gamma/\rho)$ and $R[J(\gamma/\rho)]$.

**Definition.** A poset P has the *jeu de taquin property* (or is *jdt*) if:
For any ideal $I \subseteq P$ and for all numberings $\rho$ of I,

$$R[J(\gamma_1/\rho)] = R[J(\gamma_2/\rho)] \text{ for any two numberings } \gamma_1 \text{ and } \gamma_2 \text{ of the filter } P - I.$$

Suppose that a poset P is jdt. To see that P is a fair chart, note that within each scenario we may as well renumber the victims and the survivors with green and red numbers respectively, starting with "1" each time. Suppose that a chart P is fair. To see that the poset defined by P is jdt, the g green bubbles in a scenario with r red labels should be renumbered as r+1, r+2, …, r+g. The fair chart formulation is less efficient than the jdt formulation.

The poset $\Delta_{3,2}$ does not have the jdt property: Create a second initial green numbering $\gamma'$ by interchanging $G_1$ and $G_2$, and note that $R_1$ now ends up to the right in $R[J(\gamma'/\rho)]$ instead of to the left as in Figure 3. (However, the poset $\Delta_{3,3}$ shown in Figure 6 is jdt.) Figure 2 shows the smallest "three-rowed doubly shifted shape". It also fails to be jdt; start to show this by taking I to consist of all elements strictly less than v or y. Apart from these two non-jdt posets, all of the posets appearing in the figures of this paper are jdt.

It is not hard to see that each connected component of a jdt poset must be jdt, and that a connected jdt poset must possess a unique maximal element. But it may be difficult to obtain a set of local structural conditions which are both necessary and sufficient for a poset to be jdt: Let $P_1$ be any poset with n elements. Let $P_2$ be a chain with n elements. Let P be the poset formed by dropping edges down to the maximal elements of $P_1$ from the minimal element of $P_2$. Such a poset P will clearly be jdt, and it can possess any kind of local structure within its $P_1$ portion. The two non-d-complete connected 8-element jdt posets shown in Figure 4 have "necks" with 2 and 1 elements respectively. Only 3 of the 55 non-d-complete connected 8-element jdt posets have no neck elements [Wil].

The d-complete property gets its name from the minuscule posets $d_k(1)$. For $k \geq 3$, $d_k(1) \cong \Delta_{k-2,k-2}$. Let P be a poset. Let $k \geq 3$. An interval in P is a $d_k$-*interval* if it is isomorphic to $\Delta_{k-2,k-2}$. A $d_3^-$-*interval* [w;x,y] consists of two elements x and y which both cover a third element w. It is *completed* if there exists an element z covering x and y such that [w,z] is a $d_3$-interval. See the second diagram in Figure 4. Two $d_3^-$-intervals [w;x,y] and [w´;x´,y´] *overlap* if x = x´, y = y´, and w ≠ w´. See the first diagram in Figure 5. Let $k \geq 4$. A $d_k^-$-*interval* is an interval which is isomorphic to $\Delta_{k-2,k-3}$. The interval [w,y] on the right in Figure 2 forms a $d_4^-$-interval. A $d_k^-$-interval [w,y] is *completed* if there exists an element z covering y such that [w,z] is a $d_k$-interval. The interval [w,y] in Figure 2 is not completed since $[w,z] \ncong \Delta_{2,2}$. Suppose [x,y] is a $d_{k-1}$-interval and that x covers distinct elements w and w´. If [w,y] and [w´,y] are both $d_k^-$-intervals, then they *overlap*. See the second diagram in Figure 5. A poset P is d-*complete* if for every $k \geq 3$ each of the following the three axioms is satisfied:

(D1) Every $d_k^-$-interval is completed.

(D2) If [w,z] is a $d_k$-interval, then z covers only elements in [w,z].

(D3) There are no overlapping $d_k^-$-intervals.

The poset P is $d_3$-*complete* if it satisfies each axiom for k = 3. In any poset a set {w,x,y,z} is a *diamond* if z covers x and y, and both of those cover w. In a $d_3$-complete poset: It can be seen that a diamond {w,x,y,z} must be a $d_3$-interval, namely [w,z], and that z covers only x and y. A violation of Axiom D2 appears in the second diagram of Figure 4 for k = 3. The posets $\Delta_{n,n}$ are d-complete. Rooted trees are d-complete. Each poset depicted in Figures 8.1 - 8.15 can be seen to be d-complete. Clearly, any filter of a d-complete poset is d-complete. Thus shapes and shifted shapes are d-complete. But the smallest three-rowed doubly shifted shape (Figure 2) is not d-complete.

The core aspect of the interplay between the d-complete and jdt properties can be illustrated using the posets $\Delta_{b,n}$; this will highlight the importance of Axiom D1. Consider the ideal $I := \{a_b,\ldots,a_2,a_1\}$; it has only one numbering $\rho$. There are only two numberings $\gamma$ and $\gamma'$ of $P - I$; these are analogous to the $\gamma$ and $\gamma'$ for $\Delta_{3,2}$ above. If $b > n$, then it can be seen that $R[J(\gamma/\rho)] \neq R[J(\gamma'/\rho)]$ as for $3 > 2$ above. But $b \leq n$ will imply that $R[J(\gamma/\rho)] = R[J(\gamma'/\rho)]$, since the red labels will end up totally ordered within the "top chain". (Slide $G_3$ out in the second and third diagrams of Figure 6.) Applying similar reasoning for ideals $I \subset \{a_b,\ldots,a_2,a_1\}$ leads to the conclusion that $\Delta_{b,n}$ has the jdt property if and only if $b \leq n$. Since $\Delta_{b,n}$ satisfies Axiom D1 for every $k \geq 3$ if and only if $b \leq n$, it can be seen that $\Delta_{b,n}$ is also d-complete if and only if $b \leq n$.

## 3. Test Emptyings, Crucial Challenges, Solutions, Snapshots, Paths

Let $P$ be a poset. For any $n \geq 0$, let $\{G_A, G_B, G_1 \to G_2 \to \ldots \to G_n\}$ be the quasi-ordered set such that $G_A < G_B$, $G_B < G_A$, $G_A < G_1$, $G_B < G_1$. A *test numbering* $\gamma$ of some $D \subseteq P$ is a bijection from $D$ to this quasi-ordered set such that $x \leq y$ for $x,y \in D$ implies $\gamma(x) \leq \gamma(y)$. Here no requirements are placed on the relationship between $\gamma^{-1}(G_A)$ and $\gamma^{-1}(G_B)$. We call $G_A$ and $G_B$ *test bubbles*. (The word "test" is short for "A/B-test".) A *test bi-numbering* $\gamma/\rho$ of $P$ consists of a test numbering $\gamma$ of some subset $D_G \subseteq P$ and a numbering $\rho$ of some subset $D_R \subseteq P$ to the chain poset $\{R_1 \to R_2 \to \ldots \to R_r\}$ such that $D_G \cap D_R = \varnothing$ and such that $D_G \cup D_R$ is an ideal of $P$. The first diagram in Figure 6 shows a test bi-numbering $\gamma/\rho$ of $\Delta_{3,3}$ in which $n = 3$ and $r = 3$. The definitions of $M_i$ and $S_i$ can be readily extended to test bi-numberings as $i$ runs through A, B, 1, 2, …, n or through B, A, 1, 2, …, n. The second and third diagrams in Figure 6 show $S_2 \circ S_1 \circ S_B \circ S_A(\gamma/\rho)$ and $S_2 \circ S_1 \circ S_A \circ S_B(\gamma/\rho)$ respectively. Given a test bi-numbering $\gamma/\rho$ of $P$, define the *test emptyings* $J_{BA}(\gamma/\rho) := S_n \circ \ldots \circ S_2 \circ S_1 \circ S_B \circ S_A(\gamma/\rho)$ and $J_{AB}(\gamma/\rho) := S_n \circ \ldots \circ S_2 \circ S_1 \circ S_A \circ S_B(\gamma/\rho)$.

A *challenge* $(I,\rho,(x,y))$ *for* $P$ consists of an ideal $I \subseteq P$, a numbering $\rho$ of $I$, and a pair $(x,y)$ of distinct minimal elements of the filter $P - I$. Here an *xy-test numbering* $\gamma$ of an ideal $J$ of $P - I$ containing $x$ and $y$ is a test numbering of $J$ such that $\gamma(x) = G_A$ and $\gamma(y) =$ $G_B$. We say $\gamma$ is a *solution* to the challenge $(I,\rho,(x,y))$ if $R[J_{BA}(\gamma/\rho)] = R[J_{AB}(\gamma/\rho)]$. If $x,y \in P$, define the filter $\langle x,y\rangle := \{z \in P\text{: } z \geq x \text{ or } z \geq y\}$. If x and y are incomparable, then $I_{xy} := P - \langle x,y\rangle$ is an ideal such that x and y are the only minimal elements of $P - I_{xy} = \langle x,y\rangle$. Two incomparable elements $x,y \in P$ form a *crucial pair* $(x,y)$ for P if every element which covers one of these elements also covers the other, and if there exists at least one $w \in P$ such that $w < x$ and $w < y$. If $(x,y)$ is a crucial pair for P, then the challenge $(I_{xy},\rho,(x,y))$ is *crucial* and will be denoted more simply $((x,y),\rho)$.

The interchanges which occur during the calculations of test emptyings may be successfully executed in some other orders. A *snapshot* of a test emptying is any test bi-numbering produced during any valid calculation of the test emptying. Suppose a bubble E starts at $v \in P$ in a test bi-numbering $\gamma/\rho$ and then is slid out. The *path* $T(E)$ of E in P for $\gamma/\rho$ consists of the sequence $v, \ldots$ of elements of P which are successively occupied by E.

## 4. Other d-Complete Definitions and Facts

Let P be a connected d-complete poset. By Part F1 of Proposition 3 of [Pr1], we know that P has a unique maximal element t. The *top tree* T of P consists of all $x \in P$ such that $[x,t]$ is a chain. Note that T is a filter of P which is a rooted tree under the order inherited from P. An element $y \in P$ is *acyclic* if $y \in T$ and it is not in the upper chain of any $d_k$-interval in P for any $k \geq 3$.

Let $P_1$ be a connected d-complete poset containing an acyclic element y, and let $P_2$ be a connected d-complete poset. Let x denote the unique maximal element of $P_2$. Then the *slant sum* of $P_1$ with $P_2$ at y, denoted $P_1 \,{}^{y}\backslash_{x}\, P_2$, is the poset formed by creating a covering relation $x \to y$. By Proposition 4B of [Pr1], this poset is also connected and d-complete. A connected d-complete poset P is *slant irreducible* if it cannot be expressed as a slant sum of two non-empty d-complete posets. Suppose P is a connected d-complete poset with top tree T. An edge $x \to y$ of P is a *slant edge* if $x,y \in T$ and y is acyclic. By Proposition 4B, the edge $x \to y$ in $P_1 \,{}^{y}\backslash_{x}\, P_2$ is a slant edge. An *irreducible component* is a slant irreducible d-complete poset which has at least two elements. Let Q be a connected d-complete poset.

Theorem 4 of [Pr1] states that removal of all slant edges from Q produces a collection of slant irreducible d-complete posets $P_i$. These are the *slant irreducible components of* Q. Each $P_i$ is either a one element poset or an irreducible component. The overall global form of a general connected d-complete poset Q is "tree-like": As is indicated in Figure 3 of [Pr1], such a poset is formed by joining together many one element posets and many irreducible components with slant edges, whereupon the irreducible components look like bunches of grapes hanging from the vine-like connective tissue formed from the one element slant irreducible components.

Let P be an irreducible component. By Theorem 5 of [Pr1], the top tree T of P can be formed by identifying the maximal elements of *left and right chains* which have $g+1 \geq 1$ and $h+1 \geq 1$ elements respectively with the minimal element of a *top chain* which has $f+1 \geq 1$ elements. It will be assumed without loss of generality that $g \leq h$. By Proposition 4C of [Pr1], the only potential acylic elements of P are the minimal elements $l$ and $r$ of the left and right chains of T. We say P is *maximal* if it is not a filter of any other irreducible component whose top tree is T. From Section 7 of [Pr1], it can be seen that the maximal irreducible components fall into 15 families. Some of the restrictions on f, g, and h specified in Table 1 ensure the disjointness of these families. The $+45^\circ$ rotated diagrams of these maximal irreducible components appear in Figures 8.1 - 8.15 below. There $n := f+1$, and the *'s indicate relatively unimportant unnamed elements. Table 1 specifies which of the elements $l$ and $r$ is actually acyclic for a given irreducible component.

## 5. Main Results and Overview of Proof

The original definition of the jdt property requires many checks. Eriksson's preliminary reduction reduced the number of checks required for shapes; the general poset version of that result is Proposition 6.2 below. The number of checks will be further reduced by our

**Theorem 6.5.** *A poset is jdt if and only if every crucial challenge has a solution.*

It is easy to see that a poset is jdt (d-complete) if and only if each of its connected components is jdt (d-complete). This will also be clear for the "simultaneous" property. Therefore we only need to consider connected posets in the proofs of Theorems 5.1 and 5.2. By removing all of its slant edges, a connected d-complete poset can be further decomposed into a disjoint union of one element posets and irreducible components coming from the 15 classes. Let P be an irreducible component. From Figures 8.1 - 8.15, it can be seen that P has only one crucial pair $(x_0,y_0)$. Hence the crucial challenges for P are simply the numberings $\rho$ of the ideal $P - \langle x_0,y_0 \rangle$. The filter $\langle x_0,y_0 \rangle$ consists only of $x_0$, $y_0$, and the chain $t_1 \rightarrow t_2 \rightarrow \ldots \rightarrow t_n$. There is one xy-test numbering of each of the n+1 ideals of $\langle x_0,y_0 \rangle$; these are denoted by $\gamma_m$ for $0 \leq m \leq n$. In this context Theorem 6.5 becomes

**Proposition 8.1a.** *Let* P *be an irreducible component. Then* P *is jdt if and only if for every numbering* $\rho$ *of* $P - \langle x_0,y_0 \rangle$ *there exists some* $0 \leq m \leq n$ *such that* $R[J_{BA}(\gamma_m/\rho)] = R[J_{AB}(\gamma_m/\rho))]$.

If the paths formed by sliding the test bubbles $G_A$ and $G_B$ out from $x_0$ and $y_0$ through a given $\rho$ do not intersect, then it will be easy to see that $\gamma_0$ is a solution to $\rho$. For shapes P, Eriksson's second key idea was a "repair" procedure which could be used to show that $\gamma_1$ is a solution when the paths of $G_A$ and $G_B$ do intersect. With shifted shapes P, we will show that two repair procedures can be used to see that $\gamma_2$ is a solution if neither $\gamma_0$ nor $\gamma_1$ are solutions. For the various classes of irreducible components, iteration of Eriksson's repair procedure up to n times will be used to show that for each $\rho$, some $\gamma_m$ is a solution. In Section 7, we define a poset to be "simultaneous" if each of its crucial challenges can be solved in such an iterative fashion. Any simultaneous poset is automatically jdt.

Knowing that each of two irreducible components is simultaneous is not enough to conclude that a slant sum formed from the two components is also simultaneous: the addition of the connecting slant edge could conceivably disrupt a repair procedure in the upper component. This problem cannot arise if the repair procedures in the upper component satisfy certain conditions. The notion of "strongly simultaneous" for irreducible components is defined in Section 8 by adding these conditions; this leads to:

**Proposition 8.3.** *Let* Q *be a connected d-complete poset. If each irreducible component appearing in* Q *is strongly simultaneous, then* Q *is simultaneous.*

At this point the main problem will have been reduced to confirming that the irreducible components in each of the 15 classes are strongly simultaneous. This is accomplished by showing that the equality of Proposition 8.1a can always be satisfied with the iterative repair technique, and by checking the added conditions denoted by "strongly":

**Theorem 12.1.** *If* P *is an irreducible component in Classes 1 - 7 or 11, then it is strongly simultaneous.*

**Theorem 14.2.** *If* P *is an irreducible component in Classes 8 - 10 or 12 - 15, then it is strongly simultaneous.*

Combining Theorems 12.1 and 14.2 with the classification of irreducible components, Proposition 8.3, and the reduction-to-connected-posets comment proves our main result

**Theorem 5.1.** *Every d-complete poset is simultaneous, and is therefore jdt.*

Figure 5 shows two simultaneous posets which are not d-complete because they do not satisfy Axiom D3. A partial converse can be obtained by assuming one-third of the desired conclusion: Define a poset to be *non-overlapping* if it satisfies Axiom D3 for $k \geq 3$.

**Theorem 5.2.** *A non-overlapping poset is simultaneous if and only if it is d-complete.*

A ten page proof of this partial converse will be given in a future paper.

If P is a poset, let P* denote the *order dual* poset obtained by flipping the diagram of P. A poset P is *doubly jdt (doubly d-complete)* if both P and P* are jdt (d-complete). We have recently obtained a generalization of the Littlewood-Richardson rule for numberings of convex subsets of doubly jdt posets. It can be seen by the inspection of the diagrams in [Wil] that none of the non-d-complete jdt posets with 8 or fewer elements have jdt order duals. Moreover, almost all of these small posets have "necks". The discussion of the jdt property for double tailed diamond posets at the end of Section 2 leads to the expectation that the addition of "neck" elements above P to produce a larger poset P´ with improved jdt prospects will worsen the prospects for P´* to be jdt. So it would be surprising to find a non-d-complete jdt poset with 9 or more elements which has a jdt order dual. So if P is doubly jdt, it seems likely that it is doubly d-complete. Moreover, it can be seen using Figures 8.1 - 8.15 that the only connected doubly d-complete posets are the minuscule posets.

**Conjecture 5.3.** *The minuscule posets are the only connected doubly jdt posets.*

## 6. Reformulation of the Jeu de Taquin Property

The second result of this section extends Eriksson's reduction from shapes to posets; the last result strengthens this reduction.

**Lemma 6.1.** *A poset is jdt if and only if every one of its filters is jdt.*

Proof. Let P be a poset. Ignoring the trivial direction, let $F \subseteq P$ be a filter of P. Fix some numbering $\rho_0$ from $P - F$ to $\{1, 2, \ldots, |P{-}F|\}$. Let I be an ideal of F, let $\rho$ be a numbering of I, and let $\gamma_1$ and $\gamma_2$ be numberings of $F - I$. Note that $I \cup (P - F)$ is an ideal of P. Create a numbering $\rho'$ of this ideal by adding $|P - F|$ to the images of $\rho$ and then extending by $\rho_0$. Since P is jdt, we have $R[J(\gamma_1/\rho')] = R[J(\gamma_2/\rho')]$ on P. Since the red labels of $\rho_0$ are smaller than the increased original labels of $\rho$, their presence has no effect upon the final positions of those labels. Therefore $R[J(\gamma_1/\rho)] = R[J(\gamma_2/\rho)]$ on F, and so F is jdt. ▯

**Proposition 6.2.** *A poset is jdt if and only if every challenge has a solution.*

Proof. Let P be a poset. If a challenge $(I,\rho,(x,y))$ has no solutions, then any xy-test numbering $\gamma$ of $P - I$ will be such that $R[J_{BA}(\gamma/\rho)] \neq R[J_{AB}(\gamma/\rho)]$. Hence $\gamma$ can produce two numberings of $P - I$ which violate the jdt condition. Conversely, assume P is not jdt. By Lemma 6.1, there exists a filter F of P with an ideal I of F and a numbering $\rho$ of I which yields a jdt violation on F. From amongst all such choices, choose F, I, and $\rho$ such that $|F - I|$ is minimal. Let x be minimal in $F - I$. Let $\gamma$ and $\gamma'$ be numberings of $F - I$ such that $\gamma(x) = \gamma'(x) = G_1$. Set $\delta/\sigma := S_1(\gamma/\rho)$ and $\delta'/\sigma := S_1(\gamma'/\rho)$. Let $z := \delta^{-1}(G_1) = \delta'^{-1}(G_1)$ be the final location of $G_1$ in F. We claim $R[J(\gamma/\rho)] = R[J(\gamma'/\rho)]$ on F. Otherwise creating $\varepsilon_1$ and $\varepsilon_2$ on $F - I - \{x\}$ from $\gamma$ and $\gamma'$ by dropping $G_1$ and subtracting 1 from the other green labels would imply $R[J(\varepsilon_1/\sigma)] \neq R[J(\varepsilon_2/\sigma)]$ on $F - \{z\}$, violating the minimality of $|F - I|$. So every emptying starting at x yields the same result, $\Omega_x := R[J(\gamma/\rho)]$. Let $\gamma_1 := \gamma$. Since the jdt condition fails for I and $\rho$, there must be a distinct minimal $y \in F - I$ and a numbering $\gamma_2$ of $F - I$ such that $\gamma_2(y) = G_1$ and $\Omega_y := R[J(\gamma_2/\rho)] \neq \Omega_x$. Let $\varphi$ be an xy-test numbering of an ideal J of $F - I$. Let $\varphi'$ be an extension of $\varphi$ to $F - I$. Note that $R[J_{BA}(\varphi'/\rho)] = \Omega_x \neq \Omega_y = R[J_{AB}(\varphi'/\rho)]$. Since $R[J_{BA}(\varphi/\rho)] = R[J_{AB}(\varphi/\rho)]$ would imply $R[J_{BA}(\varphi'/\rho)] = R[J_{AB}(\varphi'/\rho)]$), the challenge (I,ρ,(x,y)) for F has no solutions. Extending the numbering ρ of I to a numbering ρ´ of $I \cup (P - F)$ as in the proof of Lemma 6.1 creates a challenge $(I\cup(P–F),\rho',(x,y))$ for P with no solutions. ▯

**Proposition 6.3.** *Let* (I,ρ,(x,y)) *be a challenge for a poset* P. *Let* γ *be the unique xy-test numbering of* {x,y}. *Let* $T(G_A)$ *denote the path of* $G_A$ *if it is slid out from* γ/ρ *first, and let* $T(G_B)$ *denote the path of* $G_B$ *if it is slid out from* γ/ρ *first. If* $T(G_A) \cap T(G_B) = \emptyset$, *then* γ *is a solution to* (I,ρ,(x,y)). *If* (I,ρ,(x,y)) *does not have a solution*, *there exists some* $w \in P$ *such that* $w \le x$ *and* $w \le y$.

Proof. If $G_B$ is slid out first, at each step it finds the label at the next element in $T(B)$ to be larger than any label which it might also then cover in $T(A)$. Now suppose $G_A$ is slid out first, followed by $G_B$. After $G_A$ slides out, the label now at any element of $T(A)$ will be smaller than the label which was originally there. Since the label which was originally there was not large enough to interest $G_B$ before, and since the movement of $G_A$ does not alter any of the labels on $T(G_B)$, we see that $G_B$ will still move along $T(G_B)$. Repeat the argument with A and B interchanged to see that $J_{BA}(\gamma/\rho) = J_{AB}(\gamma/\rho)$. If (I,ρ,(x,y)) has no solutions, we must have $T(G_A) \cap T(G_B) \neq \emptyset$. Any $w \in T(G_A) \cap T(G_B)$ will satisfy $w \le x$ and $w \le y$. ▯

**Corollary 6.4.** *Rooted trees are jdt.*

Proof. Let P be a poset which is not jdt. Then it has a challenge (I,ρ,(x,y)) with no solutions, which implies there exists $w \in P$ such that $w \le x$ and $w \le y$. Since x and y are incomparable, the poset P is not a rooted tree. ▯

**Theorem 6.5.** *A poset is jdt if and only if every crucial challenge has a solution.*

Suppose that we want to prove that a shape λ which is more than a chain is jdt. Here Theorem 6.5 becomes the following reduction (which was known already to Sheats): Note that λ has only one crucial pair. Let $\gamma_1$ and $\gamma_2$ be the two numberings of the "three box" subshape $\mu$ of $\lambda$. We only need to check that $R[J(\gamma_1/\rho)] = R[J(\gamma_2/\rho)]$ for every numbering $\rho$ of the skew shape $\lambda/\mu$.

Proof. Let P be a poset. The trivial direction is the same as in Proposition 6.2. Assume P is not jdt. By Proposition 6.2, there exist challenges for P which do not have solutions. Amongst these, let $(I,\rho,(x,y))$ be such that I is maximal. Let $r := |I|$ and $n := |P–I–\{x,y\}|$. Suppose there exists a minimal element w of $P – I$ which is distinct from x and y. Define $\rho'$ by extending $\rho$ with the additional value $\rho'(w) := R_{r+1}$. By the maximality of I, we know that $(I \cup \{w\},\rho',(x,y))$ has some solution $\gamma$. Let $\gamma'$ be an extension of $\gamma$ to all of $P – I – \{w\}$. The codomain of $\gamma'$ is $\{G_A,G_B,G_1,\ldots,G_{n-1}\}$. We have $R[J_{BA}(\gamma'/\rho')] = R[J_{AB}(\gamma'/\rho')]$. Hence $R_{r+1}$ will end up at the same maximal element t of P under either test emptying. Note that $G_A$ and $G_B$ start out at locations incomparable to w. Hence there is a unique earliest bubble $G_i$ amongst those which cover $R_{r+1}$ at its starting location w. When $G_i$ begins to slide out, it will trade places with $R_{r+1}$. There is no harm in interchanging $G_i$ and $R_{r+1}$ at the beginning of either test emptying calculation. In fact, we may move $R_{r+1}$ all of the way to t at the outset of either calculation by successively interchanging it with the earliest of the bubbles which covers it after each move. The bubbles $G_A$ and $G_B$ will not move during this process. Let $\delta$ denote the resulting xy-test numbering of $P – I – \{t\}$ with $\{G_A,G_B,G_1,\ldots,G_{n-1}\}$. Since $R[J_{BA}(\gamma'/\rho')] = R[J_{AB}(\gamma'/\rho')]$, sliding $G_A$, $G_B$, $G_1$, …, $G_{n-1}$ out the rest of the way in the two calculations will move $R_1$, …, $R_r$ up to the same respective locations distinct from t. Now change the $R_{r+1}$ at t to $G_n$ and define $\delta'$ to be the extension of $\delta$ with $\delta'(t) := G_n$. It is an xy-test numbering of $P – I$. Finally, sliding $G_n$ out will produce $R[J_{BA}(\delta'/\rho)]$ and $R[J_{AB}(\delta'/\rho)]$. These two red numberings are the same. This contradicts the fact that $(I,\rho,(x,y))$ had no solutions. Hence such a w cannot exist, and the elements x and y are the only minimal elements for I. Thus $I = I_{xy}$.

We further claim that there cannot exist an element w which covers x and is incomparable to y. Suppose such a w exists. Now extend $\rho$ with $\rho'(x) := R_{r+1}$. Note that w is minimal in $I_{xy} \cup \{x\}$. By the maximality of I, we know that $(I_{xy} \cup \{x\},\rho',(w,y))$ has some solution $\gamma$ which is a wy-test numbering of some ideal J of $P – I – \{x\}$. Let $\gamma'$ be an extension of $\gamma$ to all of $P - I - \{x\}$. We know that $R[J_{BA}(\gamma'/\rho')] = R[J_{AB}(\gamma'/\rho')]$. For its first step, the label $R_{r+1}$ will trade places with the earliest of the bubbles covering it: This would be $G_A$ at w, since $G_B$ is at y, which is incomparable to x. As before, it is alright to move $R_{r+1}$ to some t at the outset of either calculation. Here $G_A$ moves to x under the first interchange. Beginning with the definition of $\delta$ above, repeat seven sentences from above to obtain a contradiction. Hence every element which covers x must be comparable to y. Similarly, every element which covers y must be comparable to x.

Let z cover x. Then $z > y$. If there exists $z'$ such that $y \to z' < z$, then $z'$ would be incomparable to x, which is impossible. Hence z must cover y. Analogously, if z covers y, then z covers x. Proposition 6.3 states that there exists some element in P which is less than both x and y. Thus (x,y) is a crucial pair, and so $(I_{xy},\rho,(x,y))$ is a crucial challenge for P which has no solutions. ▯

## 7. Definition of Simultaneous Posets

This section will define a poset P to be "simultaneous" if for every crucial challenge $((x,y),\rho)$, a solution $\gamma$ can be found such that the two test emptying processes $J_{BA}$ and $J_{AB}$ for $\gamma/\rho$ are so closely related that they may be calculated "simultaneously". It will be evident that any simultaneous poset satisfies the requirements of Theorem 6.5 and hence is jdt.

Let P be a poset and let $((x,y),\rho)$ be a crucial challenge for P. Let $\gamma$ be an xy-test numbering of an ideal J of $\langle x,y\rangle$. We begin a series of recursive definitions by setting $m := 0$, renaming the test bubbles $A_1 := G_A$, $B_1 := G_B$, and setting $x_0 := x$, $y_0 := y$. Let $\gamma_0$ denote the restriction of $\gamma$ to $\{x_0,y_0\}$. Suppose $T(A_1) \cap T(B_1) = \varnothing$. Proposition 6.3 then implies that $J_{BA}(\gamma_0/\rho) = J_{AB}(\gamma_0/\rho)$. Since its proof indicates that $A_1$ and $B_1$ can be slid out simultaneously, we say $\gamma$ is a *0-simultaneous solution* to the challenge $((x,y),\rho)$.

We now begin to develop the definition of "m-simultaneous solution" for $m \geq 1$. Readers should take $m = 1$ during the first reading of this definition, and follow it with the $\Delta_{3,3}$ example which begins with the xy-test numbering shown in the first diagram of Figure 6.

Suppose $\mathcal{T}(A_m) \cap \mathcal{T}(B_m) \neq \emptyset$. Let $w_m$ denote the earliest common element in $\mathcal{T}(A_m)$ and $\mathcal{T}(B_m)$. We say the m*th collision* $\mathcal{C}_m$ occurs at the m*th collision site* $w_m$. Let $\sigma_m$ be the label currently at $w_m$. Let $x_m$ and $y_m$ denote the respective locations of $A_m$ and $B_m$ when they are each one move short of $w_m$. Use the reasoning of Proposition 6.3 to simultaneously move the test bubbles $A_m$ and $B_m$ down to $x_m$ and $y_m$ for the calculations of both $\mathcal{J}_{BA}$ and $\mathcal{J}_{AB}$, thereby producing just one snapshot. We say $\gamma$ is *capable of repairing the* m*th collision* $\mathcal{C}_m$ if all of the following are satisfied: (i) $|J| - 2 \geq m$, (ii) there exists some element $z_m$ which covers both $x_m$ and $y_m$, (iii) the *repair bubble* $G_m$ will now slide to the *repair site* $z_m$, and (iv) $\sigma_m$ is greater than the labels at any elements other than $x_m$ and $y_m$ which are covered by $z_m$ at the time $\sigma_m$ is defined. (If P is $d_3$-complete, then (ii) is satisfied uniquely via Axiom D1 and (iv) is satisfied vacuously via Axiom D2.) Whenever these four requirements are met, set $v_m := \gamma^{-1}(G_m)$ and let $\gamma_m$ denote the restriction of $\gamma$ to $\{x_0,y_0,v_1,\ldots,v_m\}$.

Assume $\gamma$ is capable of repairing the mth collision $\mathcal{C}_m$. Continue the calculations of $\mathcal{J}_{BA}(\gamma_m/\rho)$ and $\mathcal{J}_{AB}(\gamma_m/\rho)$: Slide $G_m$ down from $v_m$ to $z_m$. At the diamond $\{w_m,x_m,y_m,z_m\}$ we now have the assignments $\{\sigma_m,A_m,B_m,G_m\}$ shown in the leftmost array:

$$\begin{matrix} G_m & B_m \\ A_m & \sigma_m \end{matrix} \quad \to \quad \begin{matrix} G_m & B_m \\ \sigma_m & A_m \end{matrix} \quad \to \quad \begin{matrix} \sigma_m & B_m \\ G_m & A_m \end{matrix} \quad \to \quad \begin{matrix} \sigma_m & B_{m+1} \\ A_{m+1} & L_m \end{matrix}$$

Both $A_m$ and $B_m$ want to move to $w_m$ next. To further calculate the effect of $\mathcal{J}_{BA}$ when m is odd or of $\mathcal{J}_{AB}$ when m is even, first swap $A_m$ and $\sigma_m$ as shown. Let $\tau_m$ denote the label which will move to $y_m$ when $B_m$ departs $y_m$ under $\mathcal{J}_{BA}$. If $\tau_m$ earlier moved to $w_m$ as $A_m$ departed $w_m$, then $\tau_m < \sigma_m$ since the labels at any fixed location decrease over time. Otherwise, if $\tau_m$ was at an element other than $w_m$ covered by $y_m$, then we know that $\tau_m < \sigma_m$ because $B_m$ wanted to move to $w_m$. Now (iv) implies that $\sigma_m$ will be the largest element covered by $G_m$ after $B_m$ is slid out. Hence the next move of $G_m$ for $\mathcal{J}_{BA}$ will be to $x_m$. We may as well move it there now and bring $\sigma_m$ up to $z_m$. We have produced the snapshot depicted locally by the third array above; denote it $\Omega_{BA}$. On the other hand, performing the analogous steps for $\mathcal{J}_{AB}$ when m is odd or for $\mathcal{J}_{BA}$ when m is even will produce the snapshot depicted locally by the third array below; denote it $\Omega_{AB}$. In $\Omega_{BA}$, rename $A_m$ as the *leader bubble* $L_m$, $B_m$ as the *test bubble* $B_{m+1}$, and $G_m$ as the *test bubble* $A_{m+1}$: see the fourth array above. In $\Omega_{AB}$, rename $B_m$ as $L_m$, $A_m$ as $A_{m+1}$, and $G_m$ as $B_{m+1}$:

$$\begin{matrix} G_m & B_m \\ A_m & \sigma_m \end{matrix} \;\rightarrow\; \begin{matrix} G_m & \sigma_m \\ A_m & B_m \end{matrix} \;\rightarrow\; \begin{matrix} \sigma_m & G_m \\ A_m & B_m \end{matrix} \;\rightarrow\; \begin{matrix} \sigma_m & B_{m+1} \\ A_{m+1} & L_m \end{matrix}$$

From the rightmost arrays, $\Omega_{BA} = \Omega_{AB} =: \Omega_m$. We say the collision $C_m$ has been *repaired.* When we need to look more closely at the repair of $C_m$, in the *left repair view* we will refer to $A_m$ by the name $A_m/L_m$ as it is moved from $x_m$ to $w_m$ and to $G_m$ by the name $G_m/A_{m+1}$ as it is prematurely moved to from $z_m$ to $x_m$. Similar notation will be used within the *right repair view.* Combining the move arguments above with a tracking of the renamings yields

**Proposition 7.1.** *Let* $m \geq 1$. *Suppose* $\gamma$ *is capable of repairing the collision* $C_m$. *Slide* $G_m$ *down from* $v_m$ *to* $z_m$, *and produce* $\Omega_m$ *by replacing* $G_m$ *by* $\sigma_m$, $A_m$ *by* $A_{m+1}$, $B_m$ *by* $B_{m+1}$, *and* $\sigma_m$ *by* $L_m$. *Then sliding out* $L_m, B_{m+1}, A_{m+1}$ *from* $\Omega_m$ *and applying* $R[.]$ *will finish the calculation of* $R[J_{BA}(\gamma_m/\rho)]$ *when* $m$ *is odd and of* $R[J_{AB}(\gamma_m/\rho)]$ *when* $m$ *is even, and sliding out* $L_m, A_{m+1}, B_{m+1}$ *from* $\Omega_m$ *and applying* $R[.]$ *will finish the calculation of* $R[J_{AB}(\gamma_m/\rho)]$ *when* $m$ *is odd and of* $R[J_{BA}(\gamma_m/\rho)]$ *when* $m$ *is even.*

**Corollary 7.2.** *Let* $m \geq 1$. Suppose $C_m$ *has been repaired and* $L_m$ *has been slid out. If* $T(A_{m+1}) \cap T(B_{m+1}) = \varnothing$, *then* $R[J_{BA}(\gamma_m/\rho)] = R[J_{AB}(\gamma_m/\rho)]$.

The proof of the corollary is the same as for Proposition 6.3. If the corollary holds, we say $\gamma$ is a m-*simultaneous solution* to the challenge $((x,y),\rho)$. (After undoing the renamings, we have $J_{BA}(\gamma_m/\rho) \neq J_{AB}(\gamma_m/\rho)$ in contrast to the case $m = 0$: The application of $R[.]$ produces equality by ignoring the different ending locations of $A_1, B_1, G_1, G_2, \ldots$ in $J_{BA}(\gamma_m/\rho)$ when compared to $J_{AB}(\gamma_m/\rho)$.) If $T(A_{m+1}) \cap T(B_{m+1}) \neq \varnothing$ after $L_m$ has been slid out, increment $m$ and continue the recursive process with the definition of the next $w_m$. We continue to attempt to perform such repairs for a given $\rho$ until either we find $\gamma$ to be m-simultaneous for some $m \geq 0$ or until there is some collision $C_m$ which $\gamma$ is not capable of repairing.

**Definition.** A poset is *simultaneous* if every crucial challenge has an m-simultaneous solution for some $m \geq 0$.

The proof of Lemma 6.1 also shows that a poset is simultaneous if and only if each of its filters is simultaneous.

**Class 1 and Class 2 Examples.** Here we recall Eriksson's proof that shapes are jdt. Temporarily work in the context of Proposition 6.2 instead of Theorem 6.5. Let P be a shape; here $\rho$ is a skew tableau on a skew shape I which is such that x and y are minimal elements of $P - I$. Shapes are $d_3$-complete, and it was not difficult for Eriksson to choose $v_1 \in P - I - \{x,y\}$ and define $\gamma \equiv \gamma_1$ on $\{x,y,v_1\}$ such that $G_1$ would slide from $v_1$ to $z_1$. After performing the steps illustrated by the arrays above, the proof for shapes ends with an argument that the path of $L_1$ prevents $A_2$ from colliding with $B_2$. But if P is a shifted shape, it is possible for $A_2$ and $B_2$ to collide for some $\rho$. Figure 7 shows the paths of bubbles which arise for one such $\rho$. Here the 1's, 2's, 3's, 4's, 5's, 6's, 7's, 8's, 9's, and 0's respectively indicate the paths of $A_1$, $B_1$, $G_1$, $L_1$, $A_2$, $B_2$, $G_2$, $L_2$, $A_3$, and $B_3$. The second collision $\mathcal{C}_2$ is repaired by $G_2$, and the path of $L_2$ then prevents $A_3$ and $B_3$ from colliding. Thus this $\gamma$ is a 2-simultaneous solution to that $\rho$. ▯

The comments made at the end of Section 2 concerning the posets $\Delta_{b,n}$ and the jdt property can be revamped to exhibit some interplay between Axiom D1 and the simultaneous property. To see some interplay between Axiom D2 and the simultaneous property, start by forming a poset P from a diamond $\{w,x,y,z\}$ by adjoining an element t above z. This P is d-complete and hence simultaneous by Theorem 5.1. Next form $P'$ from P by adjoining an element q below z. This $P'$ does not satisfy Axiom D2, and there exists a numbering $\rho$ such that the the label at q will divert a repair bubble. This causes a failure of the simultaneous property, since requirement (iv) for the capability of $\gamma$ is not satisfied.

## 8. Reduction and Restriction to Irreducible Components

The hypothesis of Theorem 6.5 and the definition of simultaneous poset become much simpler to satisfy for irreducible components. Let P be an irreducible component. Then P must be a filter of one of the posets depicted in Figures 8.1 - 8.15 for some values of f, g, and h. Note that the elements $x_0$ and $y_0$ in each poset form a crucial pair. By inspection, there are no other crucial pairs in these posets. Hence naming these specific elements with the symbols $x_0$ and $y_0$ does not conflict with the usage of the symbols $x_0$ and $y_0$ in Section 7. Set $I_0 := I_{x_0y_0} = P - \langle x_0,y_0 \rangle$. The crucial challenges $((x,y),\rho)$ are now just the numberings $\rho$ of $I_0$. The filter $\langle x_0,y_0 \rangle - \{x_0,y_0\} =: \{t_1 \to t_2 \to \ldots \to t_n\}$ is the top chain of the top tree T of P, where $n := f+1$. Each of the xy-test numberings of ideals $J \subseteq P - I_0 = \{x_0,y_0,t_1,\ldots,t_n\}$ are initial portions of the xy-test numbering $\gamma_n$ defined by $\gamma_n(x_0) = G_A$, $\gamma_n(y_0) = G_B$, $\gamma_n(t_i) = G_i$ for $1 \le i \le n$. Here the definition of simultaneous poset becomes:

**Proposition 8.1b.** *Let* P *be an irreducible component. Then* P *is simultaneous if and only if for every crucial challenge* $\rho$ *there exists some* $0 \le m \le n$ *such that* $\gamma_n$ *is an* m-*simultaneous solution to* $\rho$.

Let P be an irreducible component. If P appears as a slant irreducible component of a connected d-complete poset Q, then bubbles can "escape" from P only along the slant sum edges which descend from acyclic elements of P to other slant irreducible components of Q. Let $\rho$ be a numbering of $I_0 \subseteq P$ and let $\gamma$ be an m-simultaneous solution to $\rho$ in P for some $m \ge 0$. We say $\gamma$ is a *strong* m-*simultaneous solution to* $\rho$ if for $1 \le i \le m$ no repair bubble $G_i$ nor leader bubble $L_i$ reaches an acyclic element of P during the simultaneous emptying calculations. (The bubble $G_i$ is not considered a repair bubble after it departs $z_i$ as $G_i/A_{i+1}$ or $G_i/B_{i+1}$.) And P is *strongly simultaneous* if for every numbering $\rho$ of $I_0$ there exists a solution $\gamma$ which is a strong m-simultaneous solution to $\rho$ for some $m \ge 0$. A filter of a strongly simultaneous irreducible component is not necessarily strongly simultaneous, since new acylic elements may arise for the filter. (See the footnotes to Table 1).

Let P be an irreducible component. Note that P is $d_3$-complete. Let $i \ge 1$. If a collision $\mathcal{C}_i$ occurs at $w_i$ during the test emptyings for $\gamma_n/\rho$, then $w_i$ must be the bottom element of a diamond whose unique top element $z_i$ will be the repair site for $\mathcal{C}_i$. Requirement (iv) for the ith capability of $\gamma_n$ is vacuous. The repair bubble $G_i$ starts at $v_i = t_i$ for $1 \le i \le n$.

**Proposition 8.2.** *Let* P *be an irreducible component. Then* P *is simultaneous if and only if for every numbering* $\rho$ *of* $I_0$ *there exists an* $0 \le m \le n$ *such that whenever a collision* $C_i$, $i \ge 1$, *arises during the calculation of the test emptyings for* $\gamma_n/\rho$, *then* $i \le m$ *and* $G_i$ *slides from* $t_i$ *to* $z_i$. *Further, the irreducible component* P *is strongly simultaneous if and only if in addition each such* $G_i$ *does not reach an acyclic element of* P *when sliding to* $z_i$.

Proof. Let $i \ge 1$. A repair bubble $G_i$ will be available to repair $C_i$ if and only if $i \le n$, and $C_i$ can be repaired if and only if $G_i$ reaches $z_i$. For the second statement, note that the leader bubble $L_i$ moves down from $w_i$. Each $w_i$ must be covered by at least two elements. Acyclic elements of P are located in the top tree of P, and so no element above an acyclic element in P can be covered by two or more elements. Thus $L_i$ cannot reach an acyclic element. Hence we only need to rule out the $G_i$ reaching acyclic elements to obtain *strongly* simultaneous. ▯

**Proposition 8.3.** *Let* Q *be a connected d-complete poset. If each irreducible component appearing in* Q *is strongly simultaneous, then* Q *is simultaneous.*

Proof. Let (x,y) be a crucial pair for Q. The poset Q is an iterated slant sum of its slant irreducible components. If two incomparable elements in a slant sum $P_1 {}^{z}\backslash_{w} P_2$ are each greater than a common third element, then both of the elements must be in $P_1$ or in $P_2$. Hence both x and y must be elements of the same slant irreducible component P of Q. Since P has at least two elements, it is an irreducible component. Let $\langle x,y\rangle_P$ and $\langle x,y\rangle_Q$ denote filters respectively in P and in Q. Set $I_P := P - \langle x,y\rangle_P$ and $I_Q := Q - \langle x,y\rangle_Q$. Let $\rho$ be any numbering of $I_Q$, and let $\rho'$ denote the restriction of $\rho$ to $I_P$. The image of $\rho'$ will no longer consist of consecutive labels, but this is harmless. The irreducible component P is strongly simultaneous. Let $\gamma$ be an xy-test numbering of some ideal J of $\langle x,y\rangle_P$ which is a strong m-simultaneous solution to the challenge $\rho'$ within P for some $m \ge 0$. Note that J is an ideal of $\langle x,y\rangle_Q$. Thus $\gamma$ could be an $m'$-simultaneous solution to the challenge $\rho$ within Q for some $m' \ge 0$. The original test bubbles $A_1$ and $B_1$ and the repair bubbles $G_i$, $1 \leq i \leq m$, begin within P. The calculations of $\mathcal{J}_{BA}(\gamma/\rho)$ and $\mathcal{J}_{AB}(\gamma/\rho)$ within Q will be the same as for $\mathcal{J}_{BA}(\gamma/\rho')$ and $\mathcal{J}_{AB}(\gamma/\rho')$ within P until some bubble $A_i$, $B_i$ for $1 \leq i \leq m+1$ or $G_i$, $L_i$ for $1 \leq i \leq m$ reaches an acylic element of P and leaves P. But since $\gamma$ is a strong m-simultaneous solution, only test bubbles may reach acyclic elements of P. If for some $0 \leq i \leq m$ a test bubble $A_{i+1}$ (or $B_{i+1}$) leaves P after reaching an acyclic element u, then no collision $\mathcal{C}_{i+1}$ can arise: for $\mathcal{C}_{i+1}$ to occur outside of P, then $B_{i+1}$ (respectively $A_{i+1}$) would have to also pass through u, which would contradict the definition of $\mathcal{C}_{i+1}$ as the earliest intersection between $\mathcal{T}(A_{i+1})$ and $\mathcal{T}(B_{i+1})$. Hence $\gamma$ is an m´-simultaneous solution to $\rho$ within Q for some $m' \leq m$. Therefore Q is simultaneous. ▯

## 9. Grid Regions and Paths of Bubbles

Sections 9 - 14 prove every irreducible component is strongly simultaneous using Proposition 8.2. Let P be an irreducible component whose top tree is parameterized by f, g, and h. Then P is a filter of some maximal irreducible component P´ appearing in Figures 8.1 - 8.15 whose top tree is also parameterized by f, g, and h. Except for the *strongly* simultaneous aspects of Theorems 12.1 and 14.2, any property that we consider in Sections 9 - 14 will immediately be true for each filter P of a maximal irreducible component P´ once we know that it is true for P´. So while our statements will refer to irreducible components P, readers may as well think in terms of the maximal irreducible components P´ shown in the figures.

In this section we study how the path taken by one bubble can constrain the path taken by the next bubble. The *grid region* H of an irreducible component P consists of all elements which are members of diamonds, together with all edges which lie between two such elements. After the orientation of the topmost diamond $\{a_1,x_0,y_0,t_1\}$ has been fixed in the plane by the assumption $g \leq h$, the orientation of every other diamond in H is determined. Coordinatize only the elements of H as in a matrix: Then $x_0$ is at (2,1) and $y_0$ is at (1,2).

The following statement can be confirmed by inspection:

**Lemma 9.1.** *Let* $v_1,v_2 \in H$ *and suppose* $v_1$ *and* $v_2$ *are locations in the path* $\mathcal{T}(E)$ *of a bubble* E. *Then every* $v \in \mathcal{T}(E)$ *between* $v_1$ *and* $v_2$ *is in* H.

We often consider two bubbles E and F in some snapshot such that F is to be slid out immediately after E is slid out. When the jdt process is studied in the Class 1 context of skew tableaux, it is often noted that the path of F cannot cross the path of E. The 5's and 7's in Figure 7 illustrate this phenomenon for $A_2$ and $G_2$ after $A_2$ is created and before $A_2$ reaches an element just short of the southwestern boundary. (After $A_2$ is created during the calculation of $J_{BA}(\gamma/\rho)$, the bubble $G_2$ can be regarded as the next bubble slid out after $A_2$ in one region of P, if $G_2$ is slid out as early as possible.) The next definitions will be used to describe how the path of E limits the path of F in the irregular grid regions of Classes 2 - 15

Suppose E and F are located at elements of H with coordinates $(i_1,j_1)$ and $(i_2,j_2)$. We say F is *to the left (right)* of E if $i_2 \geq i_1$ and $j_2 < j_1$ (respectively $i_2 < i_1$ and $j_2 \geq j_1$). (This choice of terminology is motivated by temporarily unrotating the diagrams 45° so that the lines of constant $i + j$ are horizontal.) We say $T(F)$ *intersects* $T(E)$ *transversely from the left (right) at* $(i,j) \in H$ if $(i,j) \in T(F) \cap T(E)$ and $(i,j+1) \notin T(E)$ (respectively $(i+1,j) \notin T(E)$). Let $v_1$ and $v_k$ be two locations in $T(F)$. We say F *follows to the left (right) of* E *from* $v_1$ *to* $v_k$ if there exists some u in $T(E)$ such that $v_1$ is to the left (right) of u and if the portion of $T(F)$ from $v_1$ to $v_k$ does not intersect $T(E)$ transversely from the left (right).

Let $\alpha$ denote the vertical (horizontal) edge $u \rightarrow v$ in H. We say $\alpha$ is a *vertical (horizontal) swingout* if it is not the right (bottom) edge of a diamond and if v is the bottom element of a diamond. Let $\beta$ denote an edge $u \rightarrow v$ in P which is not in H. We say $\beta$ is a *tail swingout* if v is the bottom element of a diamond. In Figure 8.15, the edges $\kappa$ and $\nu$ are vertical swingouts, $\lambda$ and $\mu$ are horizontal swingouts, and $\pi$ is a tail swingout.

The 7's and 5's in Figure 7 illustrate the next result as follows: The path of $G_2$ intersects the path of $A_2$ transversely from the left only after $A_2$ uses a vertical swingout.

**Proposition 9.2.** *Let* E *and* F *be two bubbles in* H *such that* F *is slid out immediately after* E *is slid out. Suppose* F *starts to the left (right) of* E. *Further suppose* $T(F)$ *intersects* $T(E)$ *transversely from the left (right) for the first time at* y. *If* y *is not the last element of* $T(E)$ *in* H, *then it is the top element of a vertical (horizontal) swingout used by* E.

Proof. Denote the elements in $\mathcal{T}(E)$ and $\mathcal{T}(F)$ by $u_1, u_2, \ldots$ and $v_1, v_2, \ldots$ respectively. Note that $v_1, u_1, y \in H$. Let $(i_r,j_r)$ and $(p_s,q_s)$ denote the coordinates of $u_r$ and $v_s$ respectively if these elements lie in H. Let $k \geq 1$ and $t \geq 1$ be such that $y = u_t = v_k$. Suppose y is not the last element of $\mathcal{T}(E)$ in H. Then $u_{t+1} \in H$. The transverse intersection assumption implies that the coordinates of $u_{t+1}$ must be $(i_t+1,j_t)$. Since F and E began at distinct locations in H with F to the left of E, we have $k > 1$. We claim that $v_{k-1}$ was at $(i_t,j_t-1)$. If $t = 1$, then $v_1$ being to the left of $u_1 = u_t$ implies that $p_1 = \ldots = p_k = i_t$. Thus $p_{k-1} = i_t$, which forces $q_{k-1} = j_t-1$. Suppose $t > 1$. By the earliest left transverse intersection assumption we cannot have both $v_{k-1}$ and $u_{t-1}$ at $(i_t-1,j_t)$. Having $v_{k-1}$ at $(i_t-1,j_t)$ and $u_{t-1}$ at $(i_t,j_t-1)$ can be ruled out by the same assumption by tracing back along each path and using an inequality argument. This will verify the claim. Both $u_{t+1}$ and $u_t$ are in H, and so $u_{t+1} \rightarrow u_t$ is a vertical edge in H. Both $v_k$ and $v_{k-1}$ are in H, and so $v_k \rightarrow v_{k-1}$ is a horizontal edge in H. Suppose this edge is the top edge of a diamond. Then there exists $x \in H$ with coordinates $(i_t+1,j_t-1)$. Before E arrived, the label $\rho(u_{t+1})$ at $u_{t+1}$ was less than the label $\rho(x)$ at x. When E moved from $u_t$ to $u_{t+1}$, it lifted $\rho(u_{t+1})$ to $u_t$. When F is descending from $v_{k-1}$, it must move to the location of the largest label which it covers, which would be x, not y. This contradiction implies that $v_k \rightarrow v_{k-1}$ is not the top edge of a diamond. Since this edge is in H, it must be bottom edge of a diamond. Hence y is the bottom element of a diamond. If $u_{t+1} \rightarrow u_t$ is the right edge of a diamond, there would exist $x \in H$ with coordinates $(i_t+1,j_t-1)$. But this was just ruled out. So $y = u_t$ is the top element of a vertical swingout. ▯

The 7's and 5's in Figure 7 illustrate the next result as follows: The bubble $G_2$ follows to the left of $A_2$ before $A_2$ uses a vertical swingout.

**Corollary 9.3.** *Let* E *and* F *be two bubbles in* H *such that* F *is slid out immediately after* E *is slid out. Suppose* F *starts at* $v_1$ *to the left (right) of* E. *If* E *uses a vertical (horizontal) swingout as it slides out, let* $(i_h,j_h)$ *be the coordinates of the top element* $u_h$ *of the first such edge which it uses. Otherwise let* $(i_h,j_h)$ *be the coordinates of the last element* $u_h$ *of* H *which* E *reaches. Let* $v_k$ *be the last element in* $\mathcal{T}(F) \cap H$ *whose coordinates* $(p_k,q_k)$ *are such that* $q_k < j_h$ $(p_k < i_h)$. *Then* F *follows to the left (right) of E from* $v_1$ *to* $v_k$.

Proof. We just proved that $\mathcal{T}$(F) cannot transversely intersect $\mathcal{T}$(E) from the left any earlier than $(i_h,j_h)$. The element $v_k$ occurs in $\mathcal{T}$(F) before $(i_h,j_h)$ can be reached since $q_k < j_h$. Hence F follows to the left of E from $v_1$ to $v_k$. ▯

## 10. Collisions and Repairs

In Sections 10 - 14 we work in the following context: Let P be an irreducible component with top tree T parameterized by (f,g,h), grid region H, and crucial pair $(x_0,y_0)$. Let $n := f+1$, and set $I_0 := P - \langle x_0,y_0\rangle$. The xy-test numbering $\gamma_n$ defined before Proposition 8.1b assigned the test bubbles $G_A$ and $G_B$ and the repair bubbles $G_1, \ldots, G_n$ to the elements of $\langle x_0,y_0\rangle$. Let ρ be a numbering of $I_0$. We will want to show that $\gamma_n$ is a strong m-simultaneous solution to ρ for some $m \geq 0$ by satisfying the requirements of Proposition 8.2. Each result in these sections is to be understood in the following context: The test emptyings $\mathcal{J}_{BA}(\gamma_n/\rho)$ and $\mathcal{J}_{AB}(\gamma_n/\rho)$ have been simultaneously calculated up until the situation described in the hypothesis of the result at hand in terms of the notation developed in Section 7. We often refer to test bubbles as *testers*, leader bubbles as *leaders*, and repair bubbles as *fixers*.

Next we begin to analyze how collisions between the test bubbles can arise. The results obtained later in this section will be used to guarantee that repair bubbles will arrive at the required repair sites. Individual labels were mentioned when the repair procedure was described in the definition of "m-simultaneous solution", and again when the non-crossing results were proved in the preceding section. No individual labels will be mentioned in Sections 10 - 14; instead we will consider the various possible bubble paths.

**Lemma 10.1.** *Suppose the test bubbles* $A_i$ *and* $B_i$ *collide. Then the entire paths of* $A_i$ *and* $B_i$ *and the convex hull in the plane of these paths are contained in* H. *If this collision is repaired, then the leader bubble* $L_i$ *is created at an element of* H.

Proof. Since P is $d_3$-complete, the creation sites $x_{i-1}$ and $y_{i-1}$ for $A_i$ and $B_i$, the ending sites $x_i$ and $y_i$ for $A_i$ and $B_i$, and the collision site $w_i$ for $A_i$ and $B_i$ are all elements of diamonds. Hence $x_{i-1}, y_{i-1}, x_i, y_i, w_i \in H$. By Lemma 9.1, the paths of $A_i$ and $B_i$ are contained in H. By inspection, there are no holes in the grid regions of irreducible components. The leader $L_i$ is created at $w_i$. □

Suppose the test emptyings $J_{BA}(\gamma_n/\rho)$ and $J_{AB}(\gamma_n/\rho)$ have been simultaneously calculated through k stages. Then we may regard the bubbles as having been slid out in $2^{2k+1}$ different possible orders so far: $A_1$ or $B_1$, $G_1$, $(A_1/L_1,G_1/A_2)$ or $(B_1/L_1,G_1/B_2)$, $L_1$, $A_2$ or $B_2$, ..., $L_k$, $A_{k+1}$ or $B_{k+1}$. Since each of these orders produce the same patterns of labels, we may choose any one of them when we constrain the motion of the "next" bubble with Corollary 9.3. For example, if we want to use the path of $B_1$ to limit the northeastern progress of $G_1$, we may regard $B_1$ as having been slid out after $A_1$: then the sliding of $B_1$ immediately precedes the sliding of $G_1$. And we can then reverse that viewpoint to use the path of $A_1$ to limit the southwestern progress of $G_1$.

When we outlined Eriksson's proof that shapes were jdt in Section 7, we noted that the last step should argue that the path of $L_1$ prevents $A_2$ and $B_2$ from colliding. A simpler version of Corollary 9.3 can be used for that argument. More generally, the path of $L_i$ will prevent $A_{i+1}$ and $B_{i+1}$ from colliding in any irreducible component P as long as $L_i$ remains within the "interior" of the grid region H of P. In Figure 7, the 8's, 9's, and 0's illustrate the fact that the path of $L_2$ prevents $A_3$ and $B_3$ from colliding. The next result gives our formal description of this phenomenon. The 4's, 5's and 6's in Figure 7 illustrate the fact that after the path of $L_1$ uses a vertical swingout as it reaches the southwestern boundary, it no longer prevents $A_2$ and $B_2$ from colliding.

**Proposition 10.2.** *Let* $i \geq 1$. *Suppose* $L_i$ *does not use any swingouts prior to reaching* $v \in H$. *If* $A_{i+1}$ *and* $B_{i+1}$ *collide at* $w_{i+1}$, *then* $v \geq w_{i+1}$.

Proof. Suppose $w_{i+1}$ is the last element of $T(L_i)$ in H. Then $w_{i+1}$ is the top element of a tail swingout used by $L_i$, and hence $v \geq w_{i+1}$. Otherwise $w_{i+1}$ is not the last element of $T(L_i)$ in H. If $w_{i+1} = w_i$, then either $T(A_{i+1})$ or $T(B_{i+1})$ intersects $T(L_i)$ transversely at $w_{i+1}$. Here let $u := w_{i+1}$. Otherwise $w_{i+1} < w_i$. The slideout path of $L_i$ begins at $w_i$ within what will become the convex hull of the paths of $A_{i+1}$ and $B_{i+1}$. By Lemma 10.1, this convex hull is contained in H. Planar topology implies that one of $\mathcal{T}(A_{i+1})$ or $\mathcal{T}(B_{i+1})$ intersects $\mathcal{T}(L_i)$ transversely from the left or from the right respectively at a location no lower than it intersects the other of $\mathcal{T}(A_{i+1})$ or $\mathcal{T}(B_{i+1})$. Let $u \in H$ be the location of the earliest such transverse intersection. So $u \geq w_{i+1}$. In either case, Proposition 9.2 implies that u is the top element of a vertical or horizontal swingout used by $L_i$. Hence $v \geq u \geq w_{i+1}$. ∎

There are n fixers $G_1, G_2, \ldots, G_n$ available in P. The 1's, 2's, and 3's in Figure 7 illustrate the fact that the paths of $A_1$ and $B_1$ force $G_1$ to arrive at the repair site $z_1$ for $C_1$.

**Proposition 10.3.** *Let* $1 \leq i \leq n$. *Suppose* $A_i$ *and* $B_i$ *collide. Further suppose there is some element* $u \in \mathcal{T}(A_i)$ *and some element* $u' \in \mathcal{T}(G_i)$ *such that* $u'$ *is to the right of* u. *Also suppose there is some element* $v \in \mathcal{T}(B_i)$ *and some element* $v' \in \mathcal{T}(G_i)$ *such that* $v'$ *is to the left of* v. *Then the repair bubble* $G_i$ *will arrive at the repair site* $z_i$ *for the collision* $C_i$.

Proof. Any horizontal swingout which $A_i$ could conceivably use would be along the northeastern boundary of the convex hull of the paths of $A_i$ and $B_i$ in H. But to use such an edge, the tester $A_i$ would first have to enter the path of $B_i$. Hence $A_i$ cannot use any horizontal swingouts before it stops short of the collision with $B_i$ at $w_i \in H$. Similarly, the tester $B_i$ cannot use any vertical swingouts. By Corollary 9.3, we see that $G_i$ will follow to the right of $A_i$ starting from $u'$, and to the left of $B_i$ starting from $v'$. By Lemma 10.1, we see that $G_i$ is guided by the paths of $A_i$ and $B_i$ to $z_i$. ∎

The fixer $G_1$ starts at $t_1$. When $2 \leq i \leq n$, the fixer $G_i$ automatically slides down to $t_1$ after $C_{i-1}$ is repaired.

**Corollary 10.4.** *Let* $1 \leq i \leq n$. *Any first collision* $C_1$ *will be repaired, any collision* $C_i$ *occuring at* $a_1$ *will be repaired, and the earliest collision* $C_i$ *to occur below* $a_1$ *will be repaired. In these cases, the repair bubble* $G_i$ *will not reach a minimal element of* T.

Proof. The testers $A_i$ and $B_i$ in all three situations are created at $x_0$ and $y_0$. Note that $t_1$ is to the right of the starting site $x_0$ for $A_i$ and to the left of the starting site $y_0$ for $B_i$. Apply

Proposition 10.3. The path of $A_i$ will prevent $G_i$ from reaching *l* and the path of $B_i$ will prevent $G_i$ from reaching *r*. ▯

Suppose a bubble F is slid out immediately after a bubble E is slid out. Further suppose F was following to the left of E when E uses the vertical swingout $u \to v$. The southwestern boundary of P then often forces F to reach v, where it *crosses over* and starts to follow to the right of E. Analogous mirror statements can be formed for each of the next two results. The 7's and 5's in Figure 7 illustrate the next result as follows: After the fixer $G_2$ reaches a location which is to the left of $z_1$, the tester $A_2$ uses a vertical swingout and $G_2$ crosses over to follow to the right of $A_2$, whence it is forced to reach $z_2$.

**Proposition 10.5.** *Let* $1 \leq i \leq n-1$. *Suppose* $G_{i+1}$ *reaches some* $v \in H$ *which is to the left of* $z_i$. *If* $C_{i+1}$ *arises, it will be repaired if and only if* $G_i/A_{i+1}$ *or* $A_{i+1}$ *uses a vertical swingout and* $G_{i+1}$ *crosses over to follow to the right of* $A_{i+1}$.

Proof. In the right view of the repair of $C_i$, the fixer $G_i$ becomes $B_{i+1}$. Hence $G_{i+1}$ follows to the left of $B_{i+1}$ after reaching v. In the left view, the fixer $G_i$ becomes $A_{i+1}$ and so $G_{i+1}$ is initially following to the left of $G_i/A_{i+1}$ and $A_{i+1}$. The repair site $z_{i+1}$ for $G_{i+1}$ is to the right of the ending site $y_{i+1}$ for $A_{i+1}$. Hence $G_{i+1}$ cannot reach $z_{i+1}$ unless $G_i/A_{i+1}$ or $A_{i+1}$ uses a vertical swingout. If this happens and $G_{i+1}$ crosses over to follow to the right of $A_{i+1}$, then we can apply Proposition 10.3. ▯

The 4's, 7's, 3's, and 5's in Figure 7 illustrate the next result when $i = 1$ for $L_1$, $G_2$, $G_1$, and $A_2$ as follows: Take $v = (1,2)$ and $u = (1,3)$. Then $t = (2,2)$. After $G_1$ becomes $A_2$, it follows to the left of $L_1$. Then this tester $A_2$ is forced to use the next higher vertical swingout above the vertical swingout used by $L_1$. This in turn forces $G_2$, which has been following to the left, to cross over and follow to the right of $A_2$ until it reaches $z_2$.

**Proposition 10.6.** *Let* $1 \leq i \leq n-1$. *Suppose* $L_i$ *uses a vertical swingout* $s_1 \to s_2$ *as its first swingout. If this swingout is not* $\kappa$, *also assume all of the following suppositions are true. Suppose there exists* $v \in \mathcal{T}(G_{i+1})$ *and* $u \in \mathcal{T}(G_i) \cap H$ *such that* v *is to the left of* u *and such* *that* $G_i$ *uses no vertical swingouts between* u *and its repair site* $z_i$. *If* v *is at* $(i_v,j_v)$ *in* H, *let* $t \in H$ *at* $(i_t,j_v)$ *be such that* $i_t$ *is maximal with respect to the requirement that every* $(i,j_v) \in$ H if $i_v \le i \le i_t$. *Suppose* $t \ge s_2$. *Suppose there exists* $q_1 \in H$ *such that the interval* $[s_1,q_1]$ *in* P *is as shown below. Then any collision* $C_{i+1}$ *which arises will be repaired.*

$$\begin{array}{ccccc} q_1 & \leftarrow & r_2 & \leftarrow & s_3 \\ & & \uparrow & & \uparrow \\ & & r_1 & \leftarrow & s_2 \\ & & & & \uparrow \\ & & & & s_1 \end{array}$$

Proof. If the vertical swingout $s_1 \to s_2$ is $\kappa$, then $w_i = s_2 = a_1$, and we may apply Corollary 10.4. Otherwise, by inspection, any such interval $[s_1,q_1]$ has the property that $r_2$ is the bottom of a diamond. Hence $r_1 \to r_2$ will be a vertical swingout. The collision $C_i$ occurs at $w_i \ge s_2$, and so $z_i \ge r_2$. After $G_i$ arrives at $z_i$, in the left view of the repair it becomes $G_i/A_{i+1}$ as it passes to $x_i$. It then becomes $A_{i+1}$, which follows to the left of $L_i$ until $L_i$ uses its first swingout. By hypothesis, the leader $L_i$ cannot use the vertical swingout $r_1 \to r_2$. So it must be the case that $L_i$ uses $s_2 \to s_3$ if it is not created at $s_2$. Thus $A_{i+1}$ is forced to use the vertical swingout $r_1 \to r_2$ (or $G_i/A_{i+1}$ uses this swingout if $w_i = s_2$). Since $G_i$ does not use any vertical swingouts between u and $z_i$, we see that $G_{i+1}$ follows to the left of $G_i$ until $G_i$ arrives at $z_i$. Then $G_{i+1}$ follows to the left of $G_i/A_{i+1}$ and $A_{i+1}$ until one of these uses a vertical swingout. Since $t \ge s_2$, we can form a path of southerly and easterly steps from t down to $r_1$ which describes the local southwestern boundary of H. The fixer $G_{i+1}$ must stay between this path and the paths of $G_i$, $G_i/A_{i+1}$, and $A_{i+1}$ until one of the latter two uses a vertical swingout. Then $G_{i+1}$ must cross over and begin to follow to the right of $A_{i+1}$. By Proposition 10.5, we see that any $C_{i+1}$ will be repaired. □

## 11. Preliminary Results for Classes 1 - 7 and 11

We continue to work in the context of P, T, (f,g,h), n, H, $(x_0,y_0)$, $I_0$, $\gamma_n$, $\rho$, and $\gamma_n/\rho$ established at the beginning of Section 10. Now we will refer more frequently to specific structural aspects of the 15 types of irreducible components. If an element name occurs in more than one of Figures 8.1 - 8.15, then the intervals $[v,t_1]$ generated by the elements v with that name in the various irreducible components are the same. If an edge name occurs in more than one of Figures 8.1 - 8.15, then the intervals $[v,t_1]$ generated by the bottom elements v of the edges with that name in the various irreducible components are the same. If an element or an edge of a maximal irreducible component does not exist in a particular irreducible component P, then any statement refering to that element or edge may be ignored when considering P.

To prove that an irreducible component P is simultaneous using Proposition 8.2, we must confirm two things: No more than n collisions may arise in P for any $\rho$, and every collision which arises will be repaired. Tracking how the various swingout edges are used as first swingout edges by the leader bubbles $L_i$ will play a key role in these confirmations. Suppose a collision $C_{i+1}$ arises for some $i \geq 1$ while we are attempting to prove that $\gamma_n$ is an m-simultaneous solution to $\rho$ for some $0 \leq m \leq n$. By Proposition 10.2, we know that the site $w_{i+1}$ for $C_{i+1}$ must be such that $w_{i+1} \leq v_i$, where $v_i$ is the top element of the first swingout used by $L_i$. This implies that all of the elements which $L_{i+1}$ passes through must be weakly below $v_i$. Hence the top elements $v_i$ of the first swingouts used by the leaders $L_i$ form a weakly decreasing sequence in H, and $v_i = v_{i+1}$ can happen only if $C_{i+1}$ occurs at $v_i$.

**Proposition 11.1.** *Let* P *be an irreducible component. Let* $m_1$ *and* $m_2$ *respectively be the number of times that the swingouts* $\kappa$ *and* $\lambda$ *are used as first swingouts by the leaders* $L_1, L_2,$ … . *Then* $m_1 \leq n-1$ *and* $m_2 \leq 1$.

Proof. To use $\kappa$, an $L_i$ must arise at $a_1$. Suppose $L_1$ uses $\kappa$. Let $q \leq n$ be maximal such that $L_1$ reaches $a_q$. Suppose $L_2$ uses $\kappa$. Since $L_1$ either stops at $a_q$ or moves to $b_q$, the furthest $a_j$ which $L_2$ can reach is $a_{q-1}$. By induction we can see that $L_{q-1}$ would be the last $L_i$ to possibly use $\kappa$. (Once an $L_i$ moves from $a_1$ to $b_1$, the tester $B_{i+1}$ must move to $b_0$, implying that no further collisions occur at $a_1$.) Thus $m_1 \leq q-1$, and so $m_1 \leq n-1$. Next, let $0 \leq i \leq n-1$. Suppose $L_{i+1}$ uses $\lambda$ as its first swingout. Then $w_{i+2} \leq b_2$. Either $L_{i+1}$ stops at $c_2$ or it is forced to move to the only element covered by $c_2$. If $w_{i+2} = b_2$, either of the possibilities just described for $L_{i+1}$ will prevent $L_{i+2}$ from using $\lambda$. Hence $m_2 \leq 1$. ▯

**Proposition 11.2.** *Let* P *be an irreducible component in Classes* 5 - 15. *Let* $s \geq 0$ *be the number of collisions which occur at* $a_1$. *If* $s \leq n-2$ *and a collision* $C_{s+2}$ *arises, the repair bubble* $G_{s+2}$ *starts out* $t_1, x_0, a_1$ *and does not reach the right minimal element* *r* *of* T. *If it is also known that* $C_{s+1}$ *occurs at* $b_1$ *or* $b_2$, *then* $C_{s+2}$ *will be repaired. If* $s \leq n-3$ *and a collision* $C_{s+3}$ *arises, the repair bubble* $G_{s+3}$ *starts out* $t_1, y_0$. *If* $s \leq n-1$ *and* $C_{s+1}$ *occurs at* $b_3$ *or below, then the repair bubble* $G_{s+1}$ *starts out* $t_1, y_0, a_1, a_2$ *and does not reach a minimal element of* T.

Proof. Suppose $s \leq n-2$ and $C_{s+2}$ arises. Then $w_{s+1} = b_1$, $w_{s+1} = b_2$, or $w_{s+1} \leq b_3$. The movement of $G_{s+1}/A_{s+2}$ or $G_{s+1}$ from $y_0$ to $a_1$ when $C_{s+1}$ respectively occurs at $b_1$ or below $b_1$ forces $G_{s+2}$ to move to $x_0$. It must then move to $a_1$. If $C_{s+1}$ occurs at $b_1$, the absence of an immediate swingout for $L_{s+1}$ implies that $A_{s+2}$ and $B_{s+2}$ will shift down to $a_2$ and $b_1$, which is where they are respectively created if $C_{s+1}$ occurs at $b_2$. In both cases Proposition 10.3 may be applied when $G_{s+2}$ arrives at $a_1$. Suppose $s \leq n-3$ and $C_{s+3}$ arises. The movement of $G_{s+2}$ from $x_0$ to $a_1$ forces $G_{s+3}$ to move to $y_0$. Suppose $s \leq n-1$ and $C_{s+1}$ occurs at $b_3$ or below. Here the only possible beginnings of $T(A_{s+1})$ and $T(B_{s+1})$ force $G_{s+1}$ to begin as claimed. □

## 12. Strong Simultaneity for Classes 1 - 7 and 11

Guaranteeing that every collision will be repaired can become complicated. Suppose the collision $C_i$ has been repaired. For $C_{i+1}$ to arise, Proposition 10.2 implies that $L_i$ must use a swingout. Suppose it does and that a $C_{i+1}$ does arise. For $G_{i+1}$ to eventually reach it repair site $z_{i+1}$, it must first reach the depth which is one rank less deep than $z_i$. Hence $G_{i+1}$ must be to the left or to the right of $G_i$ at this point, and $G_i$ becomes both $A_{i+1}$ and $B_{i+1}$. Proposition 10.5 implies that there must be at least one usage of a swingout by $A_{i+1}$ or $B_{i+1}$ if $G_{i+1}$ is to reach its repair site $z_{i+1}$ and repair $C_{i+1}$. Hence a sequence of m collisions and repairs can arise only if there is an associated sequence of m–1 uses of swingouts by leaders which is interwoven with a sequence of m–1 uses of swingouts by testers. Successful arguments will depend upon noting the interplay between these aspects. For example, the 4's, 5's, and 7's in Figure 7 illustrate how the usage of a swingout by $L_1$ forces $A_2$ to execute the swingout usage which is required to permit $G_2$ to reach $z_2$. We also must argue that a $\mathcal{C}_{i+1}$ cannot arise if $i+1 > n$. In Class 2 we have $n = 2$. For those P we must argue that $L_2$ prevents $A_3$ and $B_3$ from colliding. In every class there exist numberings $\rho$ for which sequences of collisions of length n do arise.

If one ignores "trivial" collisions which occur at $a_1$, then sequences of as many as 4 collisions may arise. The most difficult case is the maximal irreducible component P shown in Figure 8.15 (which is the minuscule poset $e_7(1)$), for which 4 collisions always arise. There are 78 numberings $\rho$ for which all 4 collisions take place at $a_1$. Here our qualitative arguments must guarantee that the paths of $L_1$, $L_2$, $L_3$, and $L_4$ somehow "pre-process" the arrangement of labels so that $A_5$ and $B_5$ will not collide after they leave $x_0$ and $y_0$. There are also 78 numberings $\rho$ for which all 4 collisions take place at $g_5$; these arise when $\mathcal{T}(A_1)$ and $\mathcal{T}(B_1)$ do not intersect above $g_5$. Here we must guarantee that $G_4$ arrives at $z_4 = d_4$ after various earlier bubble motions have taken place between $t_1$ and $d_4$. (The number of extensions of the closely related minuscule poset $e_6(1)$ to a total order is also 78.)

A filter P of a strongly simultaneous poset P´ will be strongly simultaneous if every acylic element of P is acyclic in P´. Let P´ be the maximal irreducible component in one of Classes 1 - 7 or 11 for a given triple of (f,g,h)-values. Any irreducible component P in the same class with the same triple of (f,g,h)-values will be a filter of P´. In Classes 1, 3 - 7, and 11, any such P and P´ will possess the same acyclic elements. In these classes we can continue to think in terms of only the maximal irreducible components. In Class 2, the left minimal element $l = x_0$ of T is not acyclic in the maximal irreducible components P when $h \geq 2$, but it is acyclic for filters $P \subseteq P'$ which have no more than 3 rows when $h \geq 1$. Hence we must consider some non-maximal irreducible components P in Class 2 separately when confirming the "strong" part of strong simultaneity.

**Theorem 12.1.** *Let* P *be an irreducible component in Classes 1 - 7 or 11. If* $\rho$ *is a numbering of* $I_0$*, then there exists* $0 \leq m \leq n$ *such that* $\gamma_n$ *is a strong* m*-simultaneous solution to* $\rho$.

Proof. Let $\rho$ be a numbering of $I_0$. When attempting to calculate the test emptyings for $\gamma_n/\rho$ simultaneously, we will automatically apply Corollary 10.4 to repair each first collision $\mathcal{C}_1$. This will rule out $G_1$ reaching $l$ or $r$.

Suppose P is in Class 1. Since there are no swingouts in P, at most one collision can arise. Under the conventions just stated, there is nothing else to do.

Suppose P is in Class 2. A $\mathcal{C}_2$ can arise only if $L_1$ uses one of the vertical swingouts and $A_2$ crosses over, or wants to cross over, the path of $L_1$. Suppose $L_1$ does not use $\kappa$. When the fixer $G_2$ reaches $t_1$, it is to the left of $G_1$ at $y_0$. The path of $A_1$ prevents $G_1$ from using any vertical swingouts. The requirements of Proposition 10.6 are satisfied here, as they are when $L_1$ does use $\kappa$. Hence any $\mathcal{C}_2$ will be repaired. Since no horizontal swingouts are available to $L_1$, the tester $B_2$ and then its renamed version $L_2$ will stay to the right of $L_1$. But the only swingouts are along the southwestern boundary of P. Hence $L_2$ can never use a swingout, and so no more than $n = 2$ collisions can arise. Since $G_2$ follows to the left of $G_1$, it is impossible for $G_2$ to reach $r$. Now suppose P has at most three rows, which implies that $l = x_0$ is acylic. Here the leader $L_1$ must use $\kappa$ for a $\mathcal{C}_2$ to arise. This implies $w_1 = a_1$, and then Corollary 10.4 rules out $G_2$ reaching $l$.

For P in Classes 3 - 7 and 11, suppose $\kappa$ is used $k \geq 0$ times. From Proposition 11.1 we know that $k \leq n-1$. The collisions $\mathcal{C}_1, \ldots, \mathcal{C}_k$ occur at $a_1$. Since $k+1 \leq n$, any $\mathcal{C}_{k+1}$ which arises will be repaired by Corollary 10.4 without $G_{k+1}$ reaching $l$ or $r$. We are done with Classes 3 and 4 since there are no other swingouts. In the remaining classes, the only acyclic element is $r$.

Continuing with Classes 5 - 7 and 11, first suppose $k = n-1$. From the first paragraph of the proof of Proposition 11.1, we see that there is a snapshot $\Omega$ in which $L_1, \ldots, L_k$ are respectively at $a_n, \ldots, a_2$. By considering the forced paths of $L_1, \ldots, L_{k-2}$ in succession, it can be seen that $L_{k-1}$ is forced to move to $a_3$, $b_3$, $c_3$, and $d_3$ respectively in Classes 6, 5, 7, and 11. This forces $L_k$ to reach $b_2$. Any further movements of $L_k$ then prevent usage by $L_n$ of $\lambda$ or, in Class 11, of any other swingout. Thus there can be no more than n collisions when $k = n-1$. Next suppose $k < n-1$. Further suppose $L_{k+1}$ uses $\lambda$ as its first swingout. Then $\mathcal{C}_{k+1}$ must have occured at $a_1, b_1$, or $b_2$. If $w_{k+1} = a_1$, then any $\mathcal{C}_{k+2}$ will be repaired by Corollary 10.4 without $G_{k+1}$ reaching $r$. If $w_{k+1} = b_1$ or $w_{k+1} = b_2$, then any $\mathcal{C}_{k+2}$ will be repaired with a $G_{k+2}$ which does not reach $r$ by Proposition 11.2. In Classes 5 - 7 there are no other possible swingout usages for $L_{k+1}$ or $L_{k+2}$ since $\lambda$ is the only swingout besides $\kappa$ and it can be used at most once by Proposition 11.1. Thus there are at most $k+2 \leq n$ collisions possible, and we are done with Classes 5 - 7.

For Class 11, continue to consider $k < n-1$. A $\mathcal{C}_{k+2}$ can arise only if $L_{k+1}$ uses some horizontal swingout. By forming the mirror analog of reasoning used for $L_2$ and $L_1$ in Class 2, it can be seen that any $L_{k+1}$ will prevent any $L_{k+2}$ from using a swingout. Thus any $\mathcal{C}_{k+2}$ will be the last collision to arise, and $k+2 \leq n$. Now suppose $L_{k+1}$ uses one of the horizontal swingouts below $\lambda$. Will $G_{k+2}$ arrive at $z_{k+2}$ without passing through $r = b_0$? The treatment above for $\mathcal{C}_{k+2}$ when $\mathcal{C}_{k+1}$ occurs at $a_1, b_1$, or $b_2$ is still good here. Suppose $w_{k+1} \leq b_3$. By Proposition 11.2, the fixer $G_{k+2}$ reaches $a_1$ and cannot reach $r$. And when $G_{k+2}$ is at $a_1$, it is to the right of $G_{k+1}$ at $a_2$. The path of $B_{k+1}$ prevents $G_{k+1}$ from using horizontal swingouts. The requirements of the mirror analog of Proposition 10.6 are satisfied. □

## 13. Preliminary Results for Classes 8 - 10 and 12 - 15

Let P be an irreducible component in Classes 8 - 10 and 12 - 15, and continue the conventions established at the beginning of Section 10.

**Lemma 13.1.** *Let* $i \geq 0$. *Suppose leaders* $L_{i+1}$ *and* $L_{i+2}$ *are located at* $d_3$ *and* $c_3$ *respectively. Any leader* $L_{i+3}$ *cannot use* $\mu$ *or* $\nu$.

Proof. If $L_{i+3}$ is to use $\mu$, it must be created at $c_3$ or reach $c_3$ from above. First assume $L_{i+2}$ moves from $c_3$ to $c_4$. For $L_{i+3}$ to arise here, the leader $L_{i+2}$ must have used $\kappa$ or $\lambda$. If $L_{i+2}$ used $\kappa$, then $L_{i+3}$ must arise at $a_1, b_1$, or $b_2$. If $L_{i+2}$ used $\lambda$, then $L_{i+3}$ must arise at $b_2$ or $b_3$. In all cases, the movement of $L_{i+2}$ to $c_4$ blocks $L_{i+3}$ from reaching $c_3$. So $L_{i+2}$ must move to $d_3$. This is impossible if $L_{i+1}$ moves from $d_3$ to $d_4$. Hence $L_{i+1}$ must move to $e_3$. Then $L_{i+1}$ stops at $e_3$ or moves to $e_4$. Either way, the leader $L_{i+2}$ must move to $d_4$. Now after $L_{i+3}$ is created at $c_3$ or arrives at $c_3$ to attempt to use $\mu$, it must move to $c_4$ instead.

If $L_{i+3}$ is to use $\nu$, it must reach $d_4$ from above or be created at $d_4$ from $A_{i+3}/L_{i+3}$ in the left view of the repair of $C_{i+3}$. Note that $A_{i+3}$ starts out following to the left of $L_{i+2}$. It will stay to the left of $L_{i+2}$ unless $L_{i+2}$ uses a vertical swingout. The only vertical swingout available prior to $\nu$ is $\kappa$. If $L_{i+2}$ used $\kappa$, then $L_{i+1}$ must have used $\kappa$ and moved to $a_3$. Then the only path for $L_{i+1}$ to $d_3$ would force $L_{i+2}$ to use $\lambda$. Then $C_{i+3}$ would have to occur at $a_1$ or $b_1$, after which $L_{i+3}$ would cross over and follow to the left of $L_{i+2}$. So we know that $A_{i+3}$ or $L_{i+3}$ begins to follow to the left of $L_{i+2}$ sometime after $L_{i+2}$ reaches $c_3$. In each class it can be seen that there is no motion of $L_{i+1}$ from $d_3$ which will permit $L_{i+2}$ to reach $e_3$ or $e_4$. Then it can be seen in each class that every motion of $L_{i+2}$ from $c_3$ which does not use $e_3$ or $e_4$ will either block $A_{i+3}$ from reaching or desiring to move to $d_4$, or else block $L_{i+3}$ from reaching $d_4$ or using $\nu$ if $L_{i+3}$ reaches $d_4$ or is created at $d_4$. ▯

**Lemma 13.2.** *For some* $q \geq 0$, *suppose each of the leaders* $L_1, \ldots, L_q$ *has reached one of* $c_2, c_3,$ *or* $c_4$ *or has been created somewhere below* $c_3$. *Then the swingout* $\pi$ *cannot be used as a first swingout by more than* $\max(3 - q, 0)$ *leaders* $L_i$ *with* $i > q$.

Proof. If $\pi$ exists, then P is in Classes 13 - 15 and $c_5$ does not exist. Hence each of $L_1, \ldots, L_q$ must reach $g_5$, or stop just short of $g_5$. In the most spacious situation of Class 15, the accumulation of leader bubbles below $g_5$ limits any usages of $\pi$ to $3 - q$ more leaders. ▯

**Proposition 13.3.** *Let* P *an irreducible component in Classes* 8 - 10 *or* 12 - 15 *and let* $\rho$ *be a numbering of* $I_0$. *Let* $m_1, \ldots, m_5$ *respectively be the number of times that the swingouts* $\kappa, \lambda, \mu, \nu, \pi$ *are used as first swingouts by the leaders* $L_1, L_2, \ldots$ . *Then* $m_3 \leq 2$, $m_4 \leq 1$, $m_5 \leq 3$, $m_2 + m_3 + m_4 \leq 2$, *and* $m_2 + m_3 + m_4 + m_5 \leq 3$.

Proof. First we show that $m_2 + m_3 \leq 2$. Suppose $L_{i+1}$ is the earliest leader to use either $\lambda$ or $\mu$ as its first swingout. For now suppose $L_{i+1}$ uses $\lambda$, after which it must reach $c_3$.

Here first suppose $L_{i+1}$ continues from $c_3$ to $c_4$. This prevents an $L_{i+2}$ from arising at $c_3$, since $A_{i+2}$ will want to move to $b_4$. This motion of $L_{i+1}$ will also block any $L_{i+2}$ arising at $b_2$ or $b_3$ from reaching $c_3$. So to avoid $m_3 = 0$ here, we have $L_{i+1}$ continuing from $c_3$ through $\mu$ to $d_3$. We can then combine this case with the case where $L_{i+1}$ uses $\mu$ as its first swingout: Now we are supposing that $\lambda$ and $\mu$ have been used a total of once as first swingouts so far, and we have $L_{i+1}$ located at $d_3$. If $m_2 + m_3 > 1$, then $L_{i+2}$ must arrive at $c_3$ to use $\mu$ as its first swingout. Apply Lemma 13.1 to see that no $L_{i+3}$ can use $\mu$. Hence $m_2 + m_3 \leq 2$, which implies $m_3 \leq 2$.

The argument for $m_4 \leq 1$ is similar to the argument for $m_2 \leq 1$ given in the proof of Proposition 11.1. Now the inequality $m_2 + m_3 + m_4 \leq 2$ can fail only if $m_2 + m_3 \geq 2$, which implies $m_2 + m_3 = 2$. Again let $L_{i+1}$ be the earliest leader to use either $\lambda$ or $\mu$ as its first swingout. The preceding paragraph indicates that $m_2 + m_3 = 2$ can only occur after a snapshot arises in which $L_{i+1}$ is at $d_3$ and $L_{i+2}$ is at $c_3$. Apply Lemma 13.1 to see that no $L_{i+3}$ can use $\nu$. Hence $m_4 = 0$ here and thus $m_2 + m_3 + m_4 \leq 2$ always. The last inequality stated holds if $m_5 = 0$. If $m_5 \geq 1$, then $\pi$ exists and $c_5$ does not exist. Then each of the leaders counted by $m_2 + m_3 + m_4$ reach one of $c_2, c_3$, or $c_4$ or was created somewhere below $c_3$. Since $m_2 + m_3 + m_4 \leq 2$, we have $3 - m_2 - m_3 - m_4 > 0$. Apply Lemma 13.2 to see that $m_5 \leq 3 - m_2 - m_3 - m_4$, which also implies $m_5 \leq 3$. ▯

**Proposition 13.4.** *Let* P *an irreducible component in Classes 8 - 10 or 12 - 15. Let* $s \geq 1$ *be the number of collisions which occur at* $a_1$. *Let* $s+1 \leq i \leq s+3$. *If* $i \leq n-1$, *any collision* $C_{i+1}$ *which arises will be repaired.*

Proof. Assume throughout that $s+1 \leq i \leq s+3$ and $i \leq n-1$. For a $C_{i+1}$ to arise, the leader $L_i$ must use $\lambda, \mu, \nu$, or $\pi$ as its first swingout. The swingout $\lambda$ cannot be used more than once as a first swingout, and the swingouts $\lambda, \mu$, and $\nu$ cannot be collectively used more than a total of two times as first swingouts. Successive first swingouts occur further and further down in the poset. So if the first swingout used by $L_i$ is $\lambda$, then $i = s+1$. If the first swingout used by $L_i$ is $\mu$ or $\nu$, then $s+1 \leq i \leq s+2$.

Suppose $L_i$ uses $\lambda$ as its first swingout. Then $i = s+1$ and $C_i$ occurs at $b_1$ or $b_2$. By Proposition 11.2, the collision $C_{i+1}$ will be repaired.

Since the part of Proposition 11.2 which was just quoted does not require $L_i$ to use $\lambda$, we can ignore the $w_i = b_1$ and $w_i = b_2$ cases below when $i = s+1$ and $L_i$ uses $\mu$, $\nu$, or $\pi$ as its first swingout. So when $i = s+1$, we will consider only $w_i \le b_3$. In this context another part of Proposition 11.2 implies that $G_{i+1}$ is following to the right of $G_i$ at $a_2$ when $G_{i+1}$ is at $a_1$. When $i = s+2$, Proposition 11.2 implies that when $G_{i+1}$ is at $y_0$ it is following to the right of $G_i$ at $a_1$.

Suppose $L_i$ uses $\mu$ as its first swingout. Then $s+1 \le i \le s+2$. Also, $w_i \ge c_3$ and $z_i \ge b_2$. Hence $G_i$ could not have used $\lambda$ or $\mu$. The requirements of the mirror analog of Proposition 10.6 are satisfied, and so $C_{i+1}$ will be repaired.

Suppose $L_i$ uses $\nu$ as its first swingout. Then $s+1 \le i \le s+2$. Also, $w_i \ge d_4$ and $z_i \ge c_3$. So $G_i$ cannot use $\mu$. We claim that $G_i$ cannot use $\lambda$ either. For $G_i$ to use $\lambda$, we must have $z_i \le c_3$, which implies $w_i \le d_4$. So we know for now that $w_i = d_4$ and $z_i = c_3$. If $i = s+1$, then $G_i$ follows to the left of $B_i$, which starts at $y_0$. Here the path of $B_i$ prevents $G_i$ from using $\lambda$. We are left with $i = s+2$. Since $\nu$ can only be used once as a first swingout, and since we are assuming that $L_i$ uses it this way, it must be the case that $L_{i-1}$ used $\lambda$ or $\mu$. Hence $w_{i-1} \ge c_3$. If $w_{i-1} = b_3$, then the usage of $\mu$ by $L_{i-1}$ forces $B_i$ to move from $c_2$ to $c_3$. Obviously this must also happen if $w_{i-1} = c_3$. The other possibilities for $w_{i-1}$ are $b_2$ and $b_1$. In these two subcases Proposition 11.2 implies that $G_i$ will reach $a_1$, whence it will begin to follow to the left of $B_i$. In all cases the path of $B_i$ prevents $G_i$ from using $\lambda$ when $i = s+2$. The claim has been confirmed. Since $L_i$ does not use $\mu$, it must pass from $c_4$ to $d_4$. If $w_i > d_4$, note that $B_{i+1}$ is following to the right of $L_i$. Thus it must use $\mu$. If $w_i = d_4$, then $G_i/B_{i+1}$ uses $\mu$. It may also be true that $B_{i+1}$ or $G_i/B_{i+1}$ used $\lambda$. Since $G_i$ does not use $\lambda$ or $\mu$, we see that $G_{i+1}$ follows to the right of $G_i$ until $G_i$ reaches $z_i$. Then $G_{i+1}$ follows to the right of $G_i/B_{i+1}$ and $B_{i+1}$ until one of these uses $\lambda$ or $\mu$. We just argued that such a usage must occur. Hence $G_{i+1}$ must cross over to follow to the left of $B_{i+1}$ and we can apply the mirror version of Proposition 10.5.

Suppose $L_i$ uses $\pi$ as its first swingout. Then we are in Classes 13 - 15 and $c_5$ does not exist. Also, $w_i \geq g_5$ and $z_i \geq d_4$. In all subcases we will show that $G_{i+1}$ is following to the left of $G_i$ after the vertical swingout $\kappa$ is no longer available to $G_i$. Then it can be seen that all of the requirements of Proposition 10.6 are satisfied, implying that $C_{i+1}$ will be repaired.

First assume $i \leq s+2$. Suppose $G_i$ uses $\lambda$. Recycling portions of the argument above which ruled this out when $L_i$ used $\nu$ as its first swingout, we see that $z_i \leq c_3$ and $w_i \leq d_4$. We also see that $i = s+2$. Dropping the assumption that $L_i$ uses $\nu$ as its first swingout, we can still follow through that earlier argument and see that usage of $\lambda$ by $G_i$ is impossible when $w_{i-1} \geq c_3$. Hence $w_{i-1} \leq b_4$ here. Hence $L_{i-1}$ must use $\nu$ or $\pi$, and so it must reach $g_5$. For $L_i$ to also use $\pi$, we must be in Classes 14 or 15. Thus $b_5$ does not exist. The eastern boundary of $P$ forces $G_{i+1}$ to cross over and begin to follow to the left of $G_i$ at $b_2$. Next suppose $G_i$ does not use $\lambda$. If $w_i \geq d_4$, then $G_i$ does not use $\mu$ either and the arguments given above when $L_i$ used $\nu$ as its first swingout can be applied here to see that $G_{i+1}$ must cross over to follow to the left of $B_{i+1}$. We are left with $w_i < d_4$. If $i = s+1$, then the path of $A_i$ forces $G_i$ to use $\mu$. If $i = s+2$, then Proposition 11.2 implies that $G_i$ reaches $a_1$, where it begins to follow to the right of $G_{i-1}$ if $w_{i-1} \leq b_3$. Even if $w_{i-1} > b_3$, after $C_{i-1}$ is repaired the fixer $G_i$ will follow to the right of $A_i$ as $A_i$ moves down to at least $d_4$. So for both values of i, the fixer $G_i$ must use $\mu$. This implies that $G_{i+1}$ must cross over and begin following to the left of $G_i$ at $c_3$.

We are left to consider the $i = s+3$ subcase when $L_i$ uses $\pi$ as its first swingout. Since they used $\lambda, \mu, \nu$, or $\pi$, the leaders $L_{i-2}$ and $L_{i-1}$ cannot reach $a_5$ or $b_5$ and must therefore reach $g_5$. Thus for $L_i$ to use $\pi$, we must be in Class 15 and $a_5$ does not exist. By Proposition 11.2, when the fixer $G_{s+4}$ starts out at $t_1$ it is following to the left of $G_{s+3}$ at $y_0$. The collision $C_{s+3}$ can occur no higher than $c_3$ since $\lambda$ can only be used once as a first swingout. Hence $z_{s+3} \leq b_2$ and so the fixer $G_{s+3}$ must use either the edge $\kappa$ or the edge $b_2 \rightarrow b_1$. If $G_{s+3}$ uses $b_2 \rightarrow b_1$, we are done. Suppose $G_{s+3}$ uses $\kappa$. Then $G_{s+4}$ crosses over and begins to follow to the right of $G_{s+3}$. We know from the $i = s+2$ portions of this proof that $G_{s+3}$ starts out following to the right of $G_{s+2}$, but then crosses over to follow to the left of $G_{s+2}$, $G_{s+2}/B_{s+3}$, or $B_{s+3}$. Therefore $G_{s+3}$ must use one of the edges $b_2 \to b_1$ or $c_3 \to c_2$ before arriving at its repair site. But since $G_{s+3}$ uses $\kappa$ here, it cannot use $b_2 \to b_1$. Hence $G_{s+3}$ must use $c_3 \to c_2$. This implies that $G_{s+3}$ uses $\lambda$. Here the eastern boundary of P forces $G_{s+4}$ to cross back over and resume following to the left of $G_{s+3}$. ▯

## 14. Strong Simultaneity for Classes 8 - 10 and 12 - 15

We continue the discussion of acyclic elements and the "strong" aspect of strong simultaneity which began just before the statement of Theorem 12.1. The irreducible components P in Classes 10 and 12 - 15 do not have acyclic elements. Let P be an irreducible component in Classes 8 or 9. In these classes, the element $r = b_0$ is acyclic when the element $e_3$ is not present in P. But any such P which does not contain $e_3$ can be regarded as a filter of the maximal irreducible component P´´ of Class 11 with the same top tree parameters f, g, and h. The proof for Class 11 in Theorem 12.1, wherein $r$ was acyclic, can therefore be applied to P. Hence we only need to prove simultaneity for Classes 8 - 10 and 12 - 15 to claim strong simultaneity. Before we complete the proof of our main result with Theorem 14.2 below, we must finish analyzing the usage of swingout edges.

**Proposition 14.1.** *The statement of Proposition 13.3 can be extended to include the bound* $m_1 + m_2 + m_3 + m_4 + m_5 \le n-1$.

Proof. We first use several facts from Propositions 11.1, Proposition 13.3, and Table 1. If $m_1 \le n-4$, then we can apply the inequality $m_2 + m_3 + m_4 + m_5 \le 3$. So assume $m_1 \ge n-3$. In Classes 9 - 10 and 12 - 15, we have $n-3 = f-2 \ge 1$. In Class 8, we have $n-3 = f-2 \ge 0$ and $m_4 = m_5 = 0$. Here the $m_1 = 0$ case can be handled by applying the inequality $m_2 + m_3 + m_4 \le 2$. So from now on we can assume $\min(1,n-3) \le m_1 \le n-1$.

From the first paragraph of the proof of Proposition 11.1, we see that there is a snapshot $\Omega$ in which $L_1, \ldots, L_{m_1}$ are respectively at $a_{m_1+1}, \ldots, a_2$. Let p be maximal such that $L_1, \ldots, L_p$ all use $\kappa$ and all reach at least one of $a_5$, $b_5$, or $c_5$. Clearly $p \le m_1$. For $1 \le i \le p-1$, each $L_{i+1}$ will begin at $a_1$ to follow to the right of $L_i$ at $a_2$. No crossing over to follow to the left of the preceding leader can occur from a usage of $\lambda$ amongst these leaders: If $L_i$ moves from $a_2$ to $b_2$ to access $\lambda$, this will imply that $i = m_1 = p$. Hence for $1 \le i \le p$, only $L_i$ can enter the ith row from the bottom of P. Since $a_5$, $b_5$, and $c_5$ are in the (n–4)th row from the bottom, we see that $0 \le p \le n-4$. If $p = m_1$, then $m_1 \le n-4$ contradicts our assumption that $m_1 \ge n-3$. So we have $p < m_1$. In $\Omega$, the leaders $L_{p+1}, \ldots, L_{m_1}$ are located within the elements from $a_4$ to $a_2$ inclusive. Hence $m_1 - (p+1) + 1 \le 3$. Thus $m_1-2 \le p+1 \le m_1$.

The leader $L_{p+1}$ starts at one of $a_4$, $a_3$, or $a_2$. Since it cannot reach $a_5$, $b_5$, or $c_5$, it must reach one of $c_2$, $c_3$, or $c_4$. Temporarily assume that the remaining zero to two leaders $L_j$ such that $p+1 < j \le m_1$ also can be shown to reach some $c_i$ for $2 \le i \le 4$. Any $L_k$ with $k \ge m_1+1$ which uses $\lambda$, $\mu$, or $\nu$ as its first swingout must reach one of $c_2$, $c_3$, or $c_4$ or be created somewhere below $c_3$. The number of such leaders is $m_2 + m_3 + m_4$. Then by Lemma 13.2, we see that $m_5 \le \max(3 - m_2 - m_3 - m_4 - (m_1 - p),0)$. If we can show that $m_1 - p \le 3 - m_2 - m_3 - m_4$, then $m_1 + m_2 + m_3 + m_4 + m_5 \le p + 3 \le n - 1$ and we will be done. So our goal will be to show that $m_2 + m_3 + m_4 \le 3 - (m_1 - p)$ for each of the possible values of $m_1-2$, $m_1-1$, or $m_1$ for $p+1$. We must also confirm that each $L_j$ with $p+1 < j \le m_1$ reaches some $c_i$ for $2 \le i \le 4$.

First assume $p+1 = m_1-2$. Here $3 - (m_1 - p) = 0$. We have $L_{p+1}$ at $a_4$, $L_{p+2}$ at $a_3$, and $L_{p+3}$ at $a_2$. Then $L_{p+1}$ must reach $d_4$ via $c_4$, which forces $L_{p+2}$ to reach $d_3$ via $c_3$, which forces $L_{p+3}$ to reach $c_2$ and then $c_3$. Then any $\mathcal{C}_{p+4}$ which arises must do so at $a_1$ or $b_1$. The movement of $L_{p+3}$ from $c_2$ to $c_3$ prevents $L_{p+4}$ from using $\lambda$. Lemma 13.1 implies that $L_{p+4}$ cannot use $\mu$ or $\nu$ either. A $\mathcal{C}_{p+5}$ can arise only if $L_{p+4}$ uses $\pi$ as its first swingout. Hence $m_2 = m_3 = m_4 = 0$.

Next assume $p+1 = m_1-1$. Here $3 - (m_1 - p) = 1$. We have $L_{p+1}$ at $a_3$ and $L_{p+2}$ at $a_2$. First assume $L_{p+1}$ reaches $d_3$ via $c_3$. Then $L_{p+2}$ reaches $c_3$ via $c_2$. Beginning with the consideration of $\mathcal{C}_{p+4}$, the second half of the preceding paragraph can be inserted here after the subscripts of this collision and all of the leaders are decreased by 1. Again $m_2 = m_3 = m_4 = 0$. Otherwise $L_{p+1}$ reaches $c_4$. Then $L_{p+2}$ must reach $c_3$. Suppose $L_{p+2}$ passes through $b_3$. Then any $\mathcal{C}_{p+3}$ can only occur at or above $b_2$, and $L_{p+3}$ will be forced to use $\lambda$ as its first swingout. The movets of $L_{p+1}$ to $d_4$, $L_{p+2}$ to $d_3$, and $L_{p+3}$ to $c_3$ are forced. Apply

Lemma 13.1 to see that any $L_{p+4}$ cannot use $\mu$ or $\nu$. A $C_{p+5}$ can arise only if $L_{p+4}$ uses $\pi$ as its first swingout. Hence $m_2 = 1$ and $m_3 = m_4 = 0$, and so $m_2 + m_3 + m_4 = 1$. The only remaining $p+1 = m_1-1$ subcase has $L_{p+1}$ reaching $c_4$ and $L_{p+2}$ reaching $c_3$ via $c_2$. Then movements of $L_{p+1}$ to $d_4$ and $L_{p+2}$ to $d_3$ are forced. The movement of $L_{p+2}$ from $a_2$ to $c_2$ implies that any $C_{p+3}$ can occur only at $b_1$ or $a_1$. Then $L_{p+3}$ must reach $b_3$ without using any swingouts. First suppose $L_{p+3}$ further reaches $b_5$ or $c_5$. Then no swingouts will be available to $L_{p+3}$, and so $m_2 = m_3 = m_4 = m_5 = 0$. Next suppose $L_{p+3}$ moves to $d_4$ via $c_4$. Now $L_{p+3}$ may use $\nu$ or $\pi$ as its first swingout. Since $\nu$ can be used as a first swingout at most once, a collision $C_{p+5}$ can arise only if $L_{p+4}$ uses $\pi$ as its first swingout. Hence $m_2 = m_3 = 0$ and $m_4 \leq 1$. Finally suppose $L_{p+3}$ moves to $c_3$. Apply Lemma 13.1 to see that no $L_{p+4}$ can use $\mu$ or $\nu$. Here $L_{p+3}$ may use $\mu$, $\nu$, or $\pi$ as its first swingout. A collision $C_{p+5}$ can arise only if $L_{p+4}$ uses $\pi$ as its first swingout. Hence $m_2 = 0$ and $m_3 + m_4 \leq 1$. In all of the sub-subcases in this paragraph, we have obtained $m_2 + m_3 + m_4 \leq 1$.

Finally assume $p+1 = m_1$. Here $3 - (m_1 - p) = 2$, and we can simply quote $m_2 + m_3 + m_4 \leq 2$ from Proposition 13.3. □

**Theorem 14.2.** *Let* P *be an irreducible component in Classes 8 - 10 or 12 - 15. If* $\rho$ *is a numbering of* $I_0$, *then there exists* $0 \leq m \leq n$ *such that* $\gamma_n$ *is a strong* m-*simultaneous solution to* $\rho$.

Proof. Let $\rho$ be a numbering of $I_0$. The length of the sequence of collisions which arises when we attempt to calculate the test emptyings for $\gamma_n/\rho$ simultaneously is $m_1 + m_2 + m_3 + m_4 + m_5$ or $m_1 + m_2 + m_3 + m_4 + m_5 + 1$. Proposition 14.1 states that $m_1 + m_2 + m_3 + m_4 + m_5 + 1 \leq n$. Thus if a collision $C_i$ arises, we know that $i \leq n$. Corollary 10.4 states that any collision occuring at $a_1$ will be repaired, and that the first collision to occur below $a_1$ will be repaired. Proposition 13.4 states that the second, third, and fourth collisions to occur below $a_1$ will be repaired. Proposition 13.3 states that the total number of usages of $\lambda$, $\mu$, $\nu$, or $\pi$ as first swingouts cannot exceed 3. These are the only swingouts available to leaders arising from collisions occuring below $a_1$. Hence at most a fourth collision below $a_1$ can arise. Therefore all collisions will be repaired. □

The author thanks Dale Peterson for telling us that every minuscule poset is jdt, John Stembridge for supplying the lists of all posets [St1] for our computations, Sarah Wilmesmeier for help in extending these computations [Wil], and Sergey Fomin for telling us about [Eri].

This packet contains the figures and tables for
"d-Complete posets generalize Young diagrams for the jeu de taquin property"
by R.A. Proctor, March 2003 version.

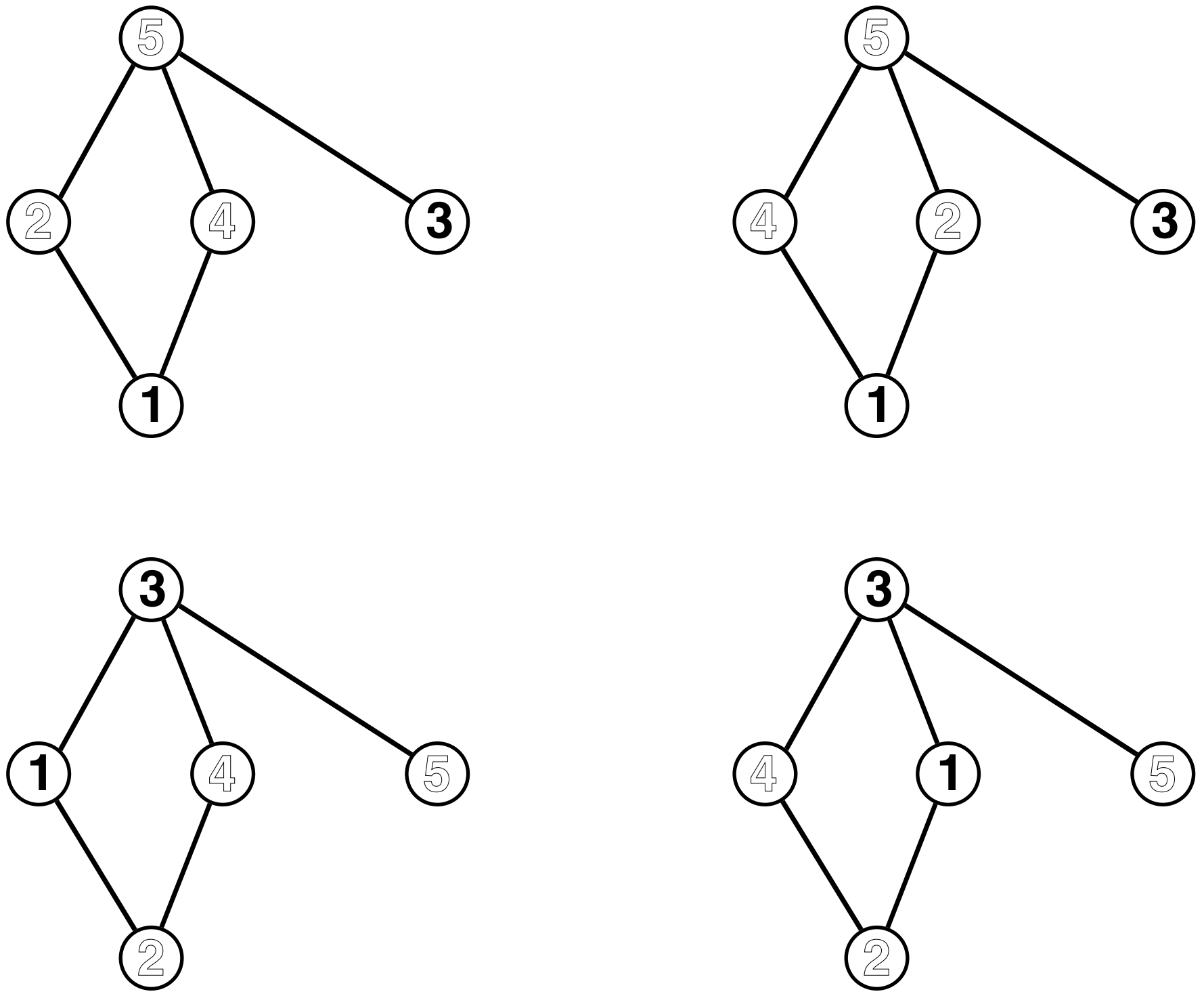

Figure 1

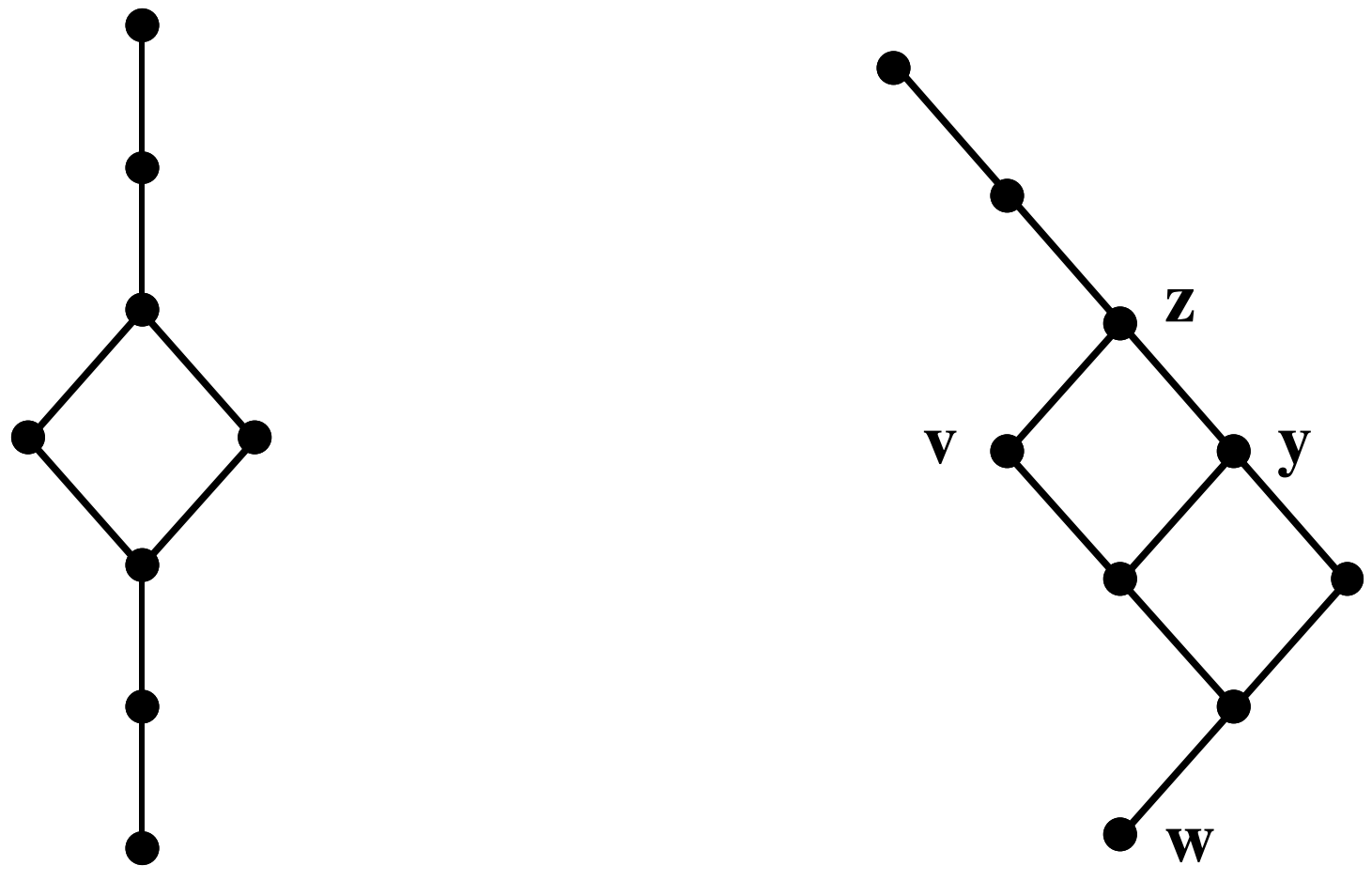

Figure 2

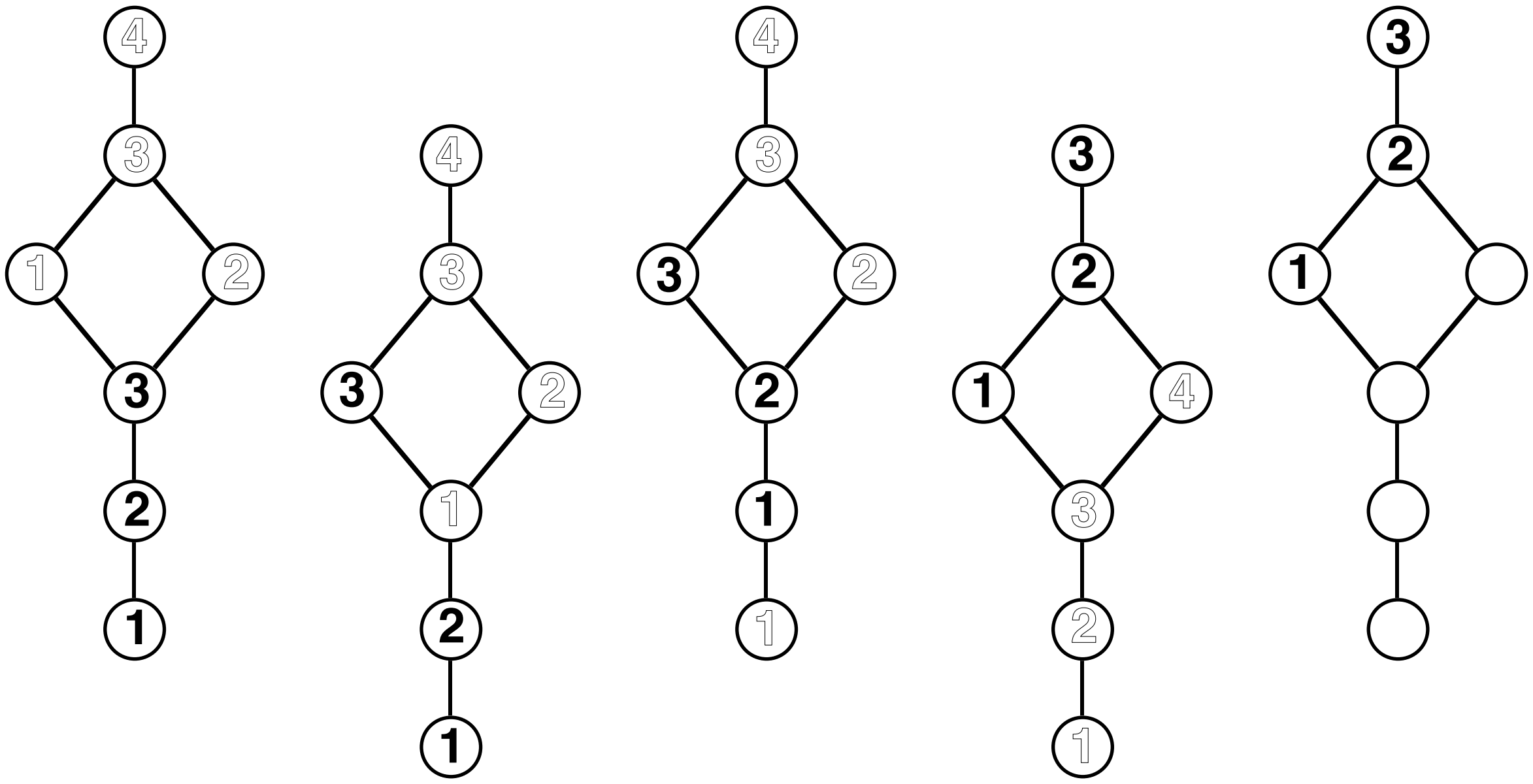

Figure 3

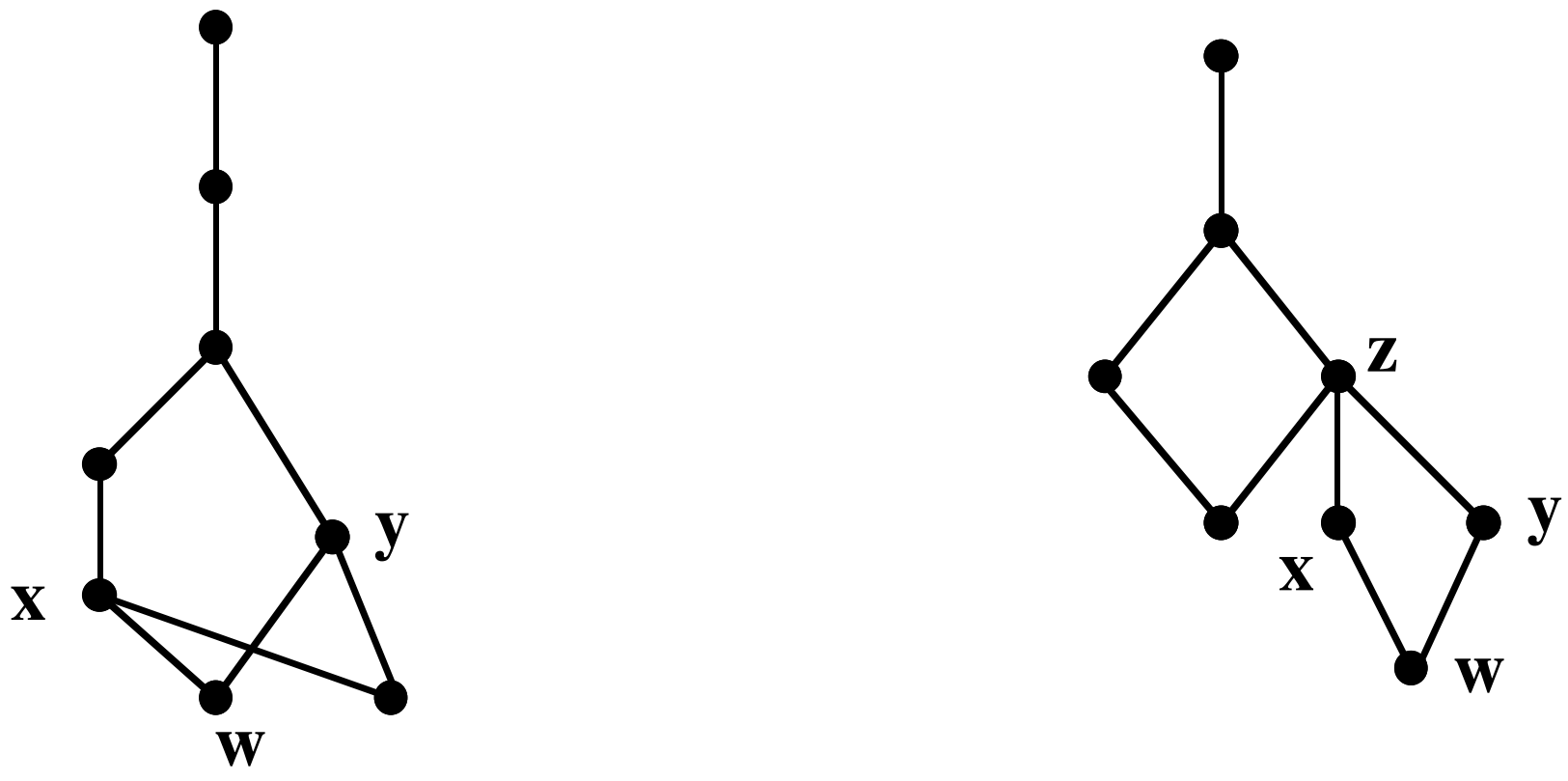

Figure 4

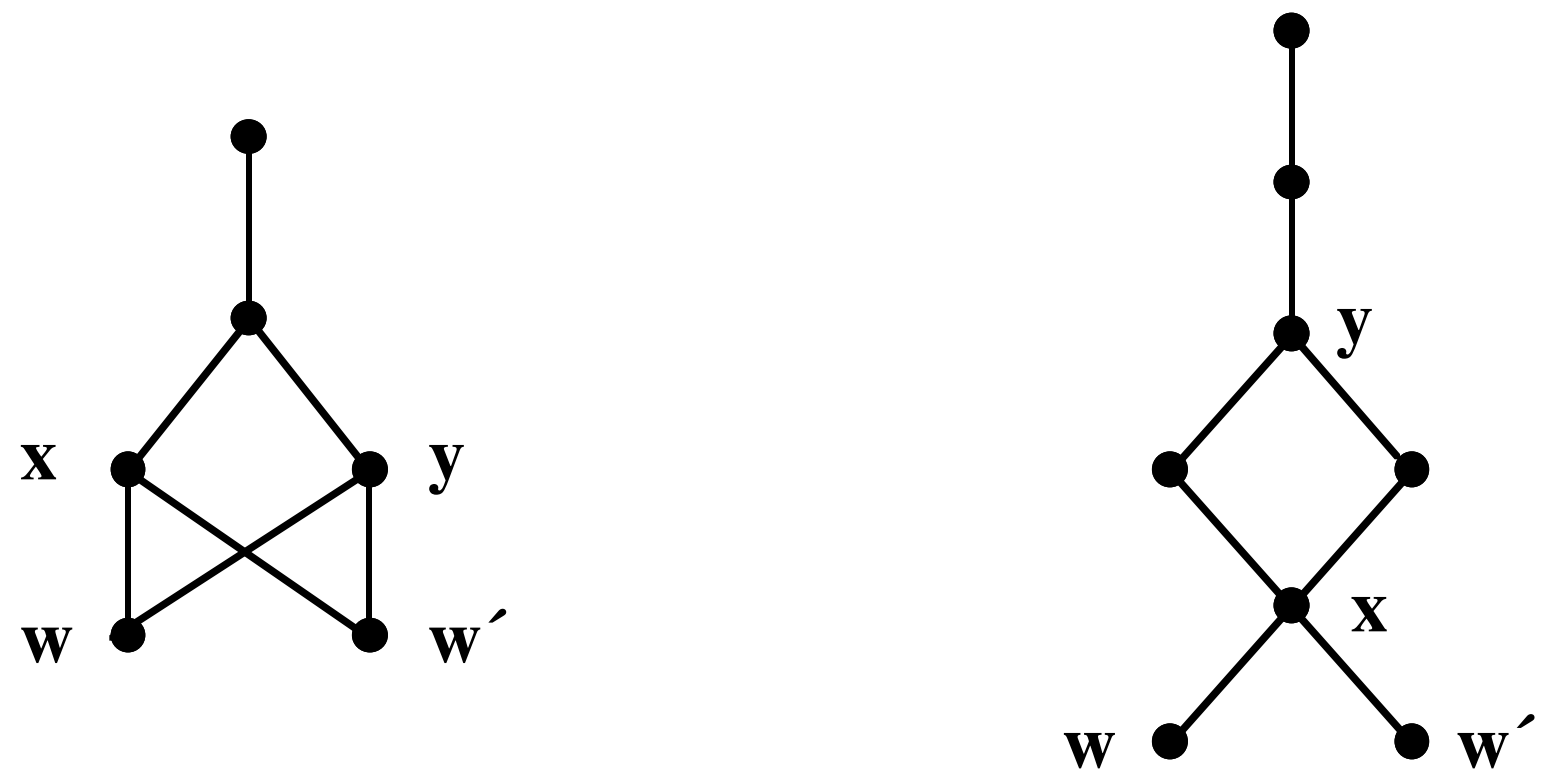

Figure 5

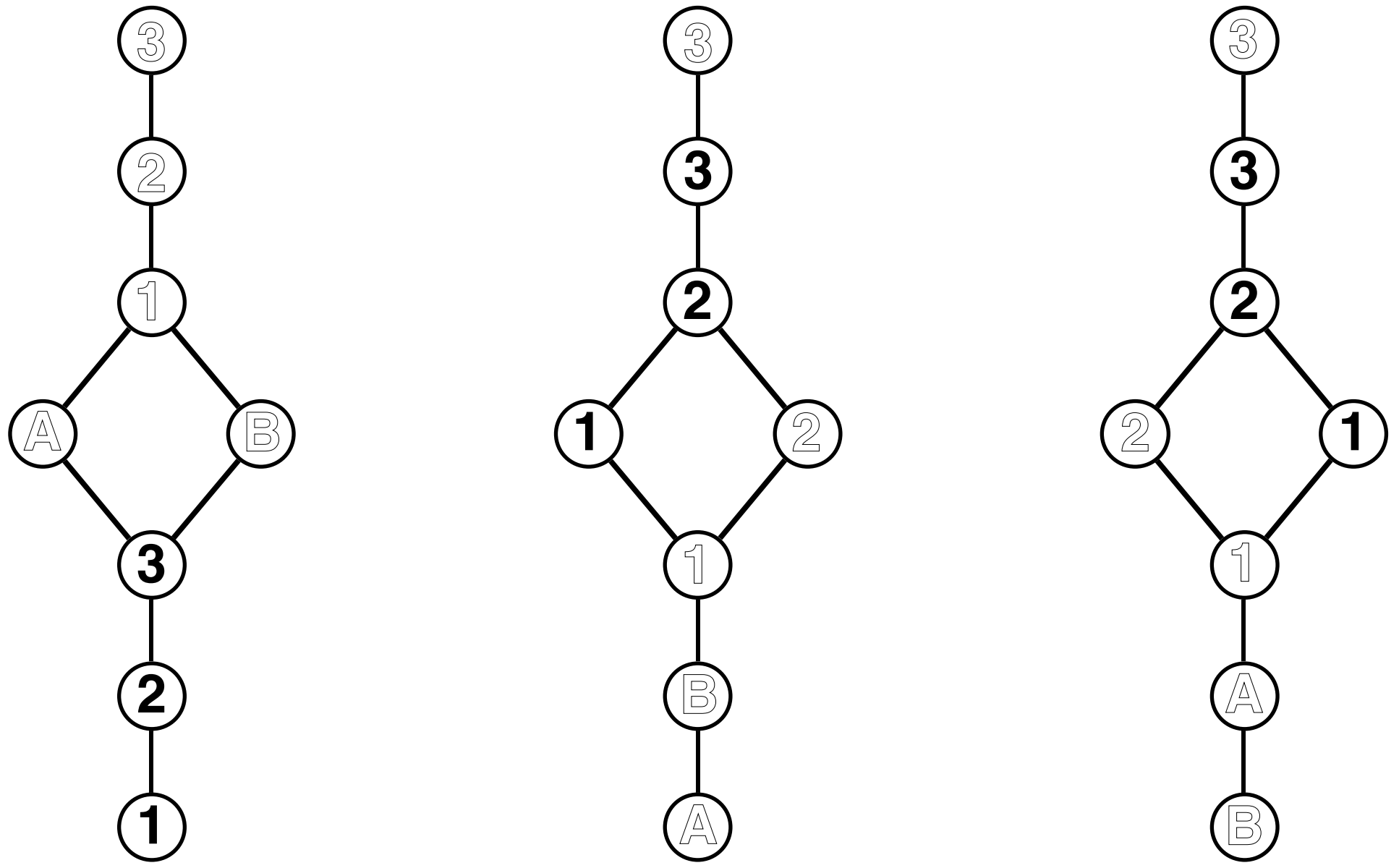

Figure 6

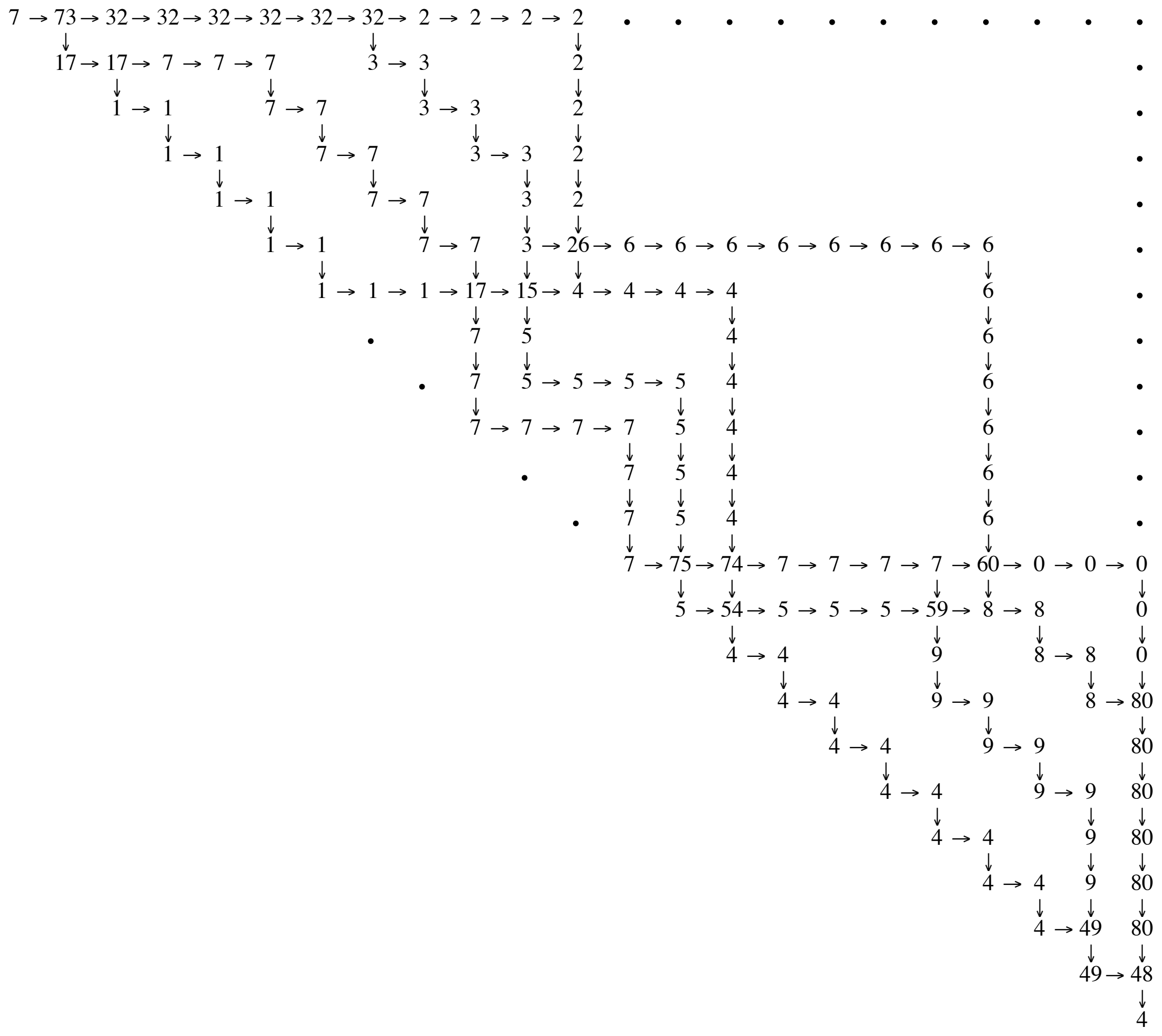
7 73 32 32 32 32 32 32 2 2 2 2
17 17 7 7 7 3 3 2
1 1 7 7 3 3 2
1 1 7 7 3 3 2
1 1 7 7 3 2
1 1 7 7 3 26 6 6 6 6 6 6 6 6
1 1 1 17 15 4 4 4 4 6
7 5 4 6
7 5 5 5 5 4 6
7 7 7 7 5 4 6
7 5 4 6
7 5 4 6
7 75 74 7 7 7 7 60 0 0 0
5 54 5 5 5 59 8 8 0
4 4 9 8 8 0
4 4 9 9 8 80
4 4 9 9 80
4 4 9 9 80
4 4 9 80
4 4 9 80
4 49 80
49 48
4

Figure 7

| Class | f | g | h | Acylics |
|---|---|---|---|---|
| 1 | = 0 | ≥ 1 | ≥ g | *l*, *r* |
| 2 | = 1 | = 1 | ≥ 1 | *l*[1], *r* |
| 3 | ≥ 1 | ≥ 2 | ≥ g | *l*, *r* |
| 4 | ≥ 2 | = 1 | ≥ 1 | *l*, *r* |
| 5 | ≥ 2 | = 1 | ≥ 2 | *r* |
| 6 | ≥ 2 | = 1 | ≥ 3 | *r* |
| 7 | ≥ 2 | = 1 | ≥ 3 | *r* |
| 8 | ≥ 2 | = 1 | = 2 | *r*[2] |
| 9 | ≥ 3 | = 1 | = 2 | *r*[2] |
| 10 | ≥ 4 | = 1 | = 2 | — |
| 11 | ≥ 4 | = 1 | = 2 | *r* |
| 12 | ≥ 3 | = 1 | = 2 | — |
| 13 | ≥ 3 | = 1 | = 2 | — |
| 14 | ≥ 3 | = 1 | = 2 | — |
| 15 | = 3 | = 1 | = 2 | — |

[1] The element $l = x_0$ is acyclic only when the element $c_2$ is not present.
[2] The element $r = b_0$ is acyclic only when the element $e_3$ is not present.

Table 1

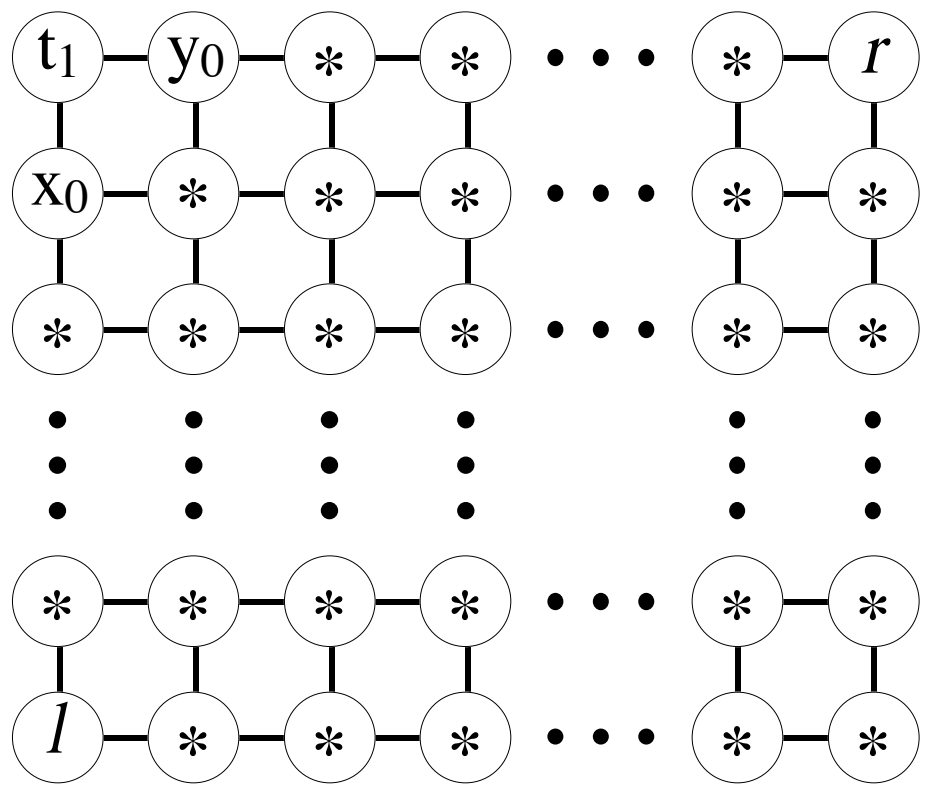

Figure 8.1

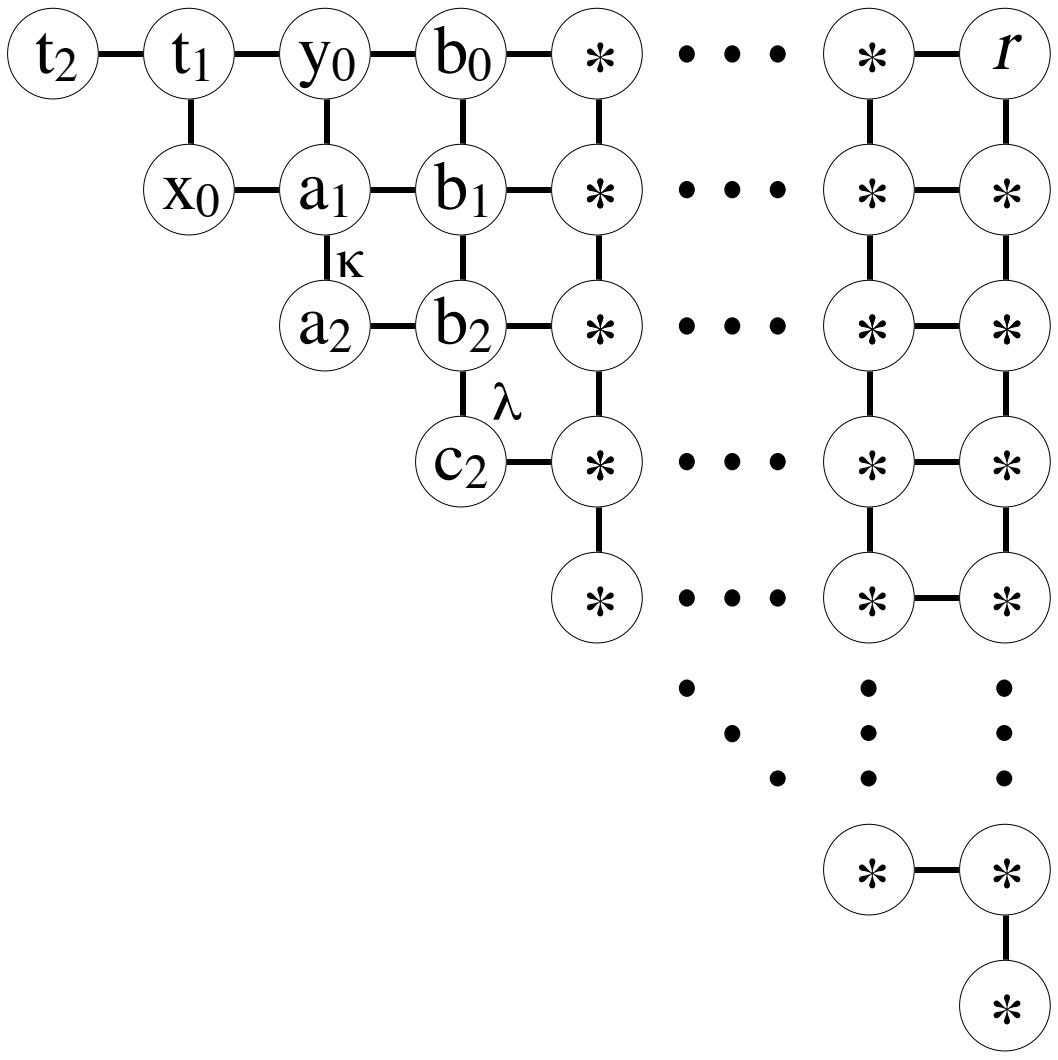

Figure 8.2

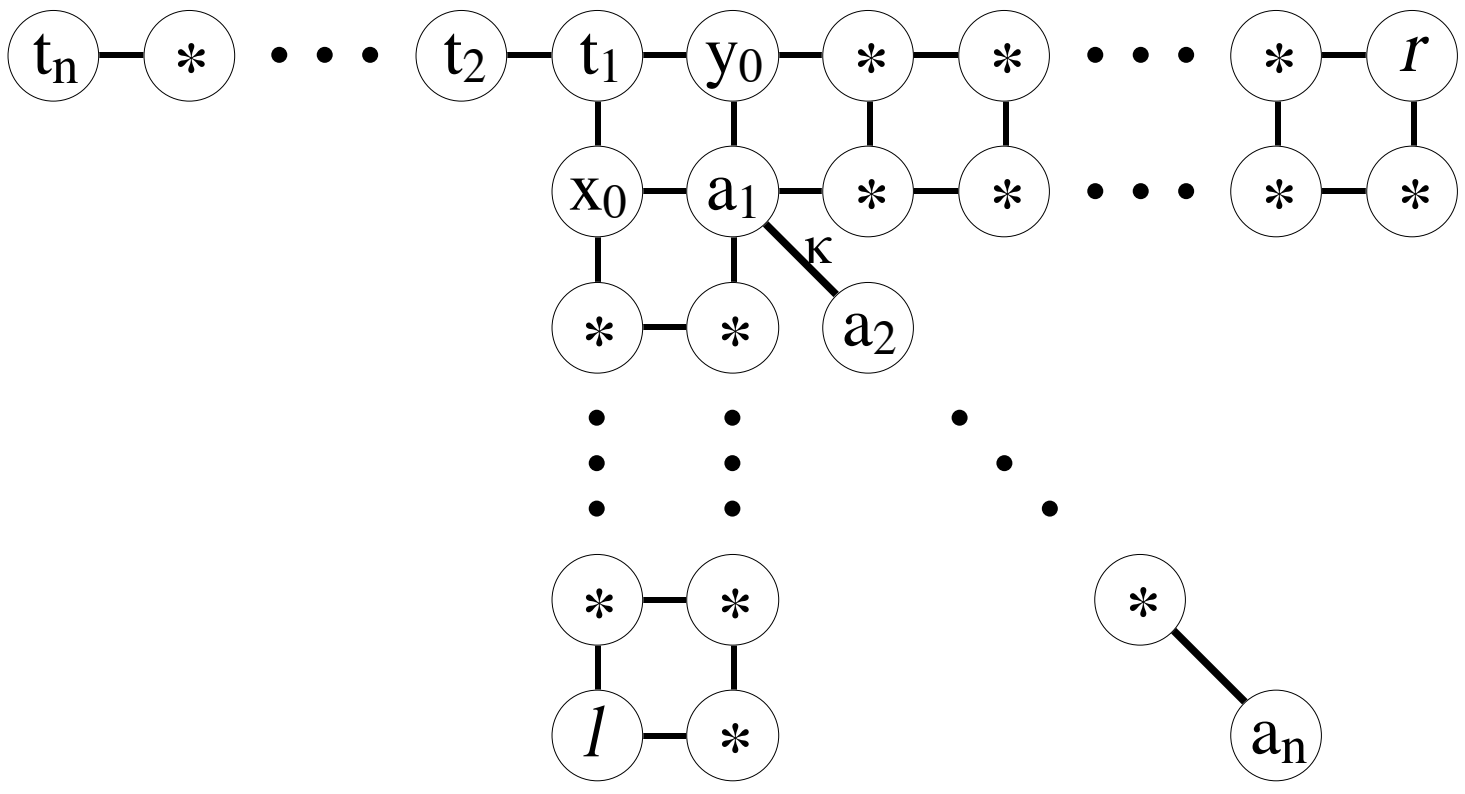

Figure 8.3

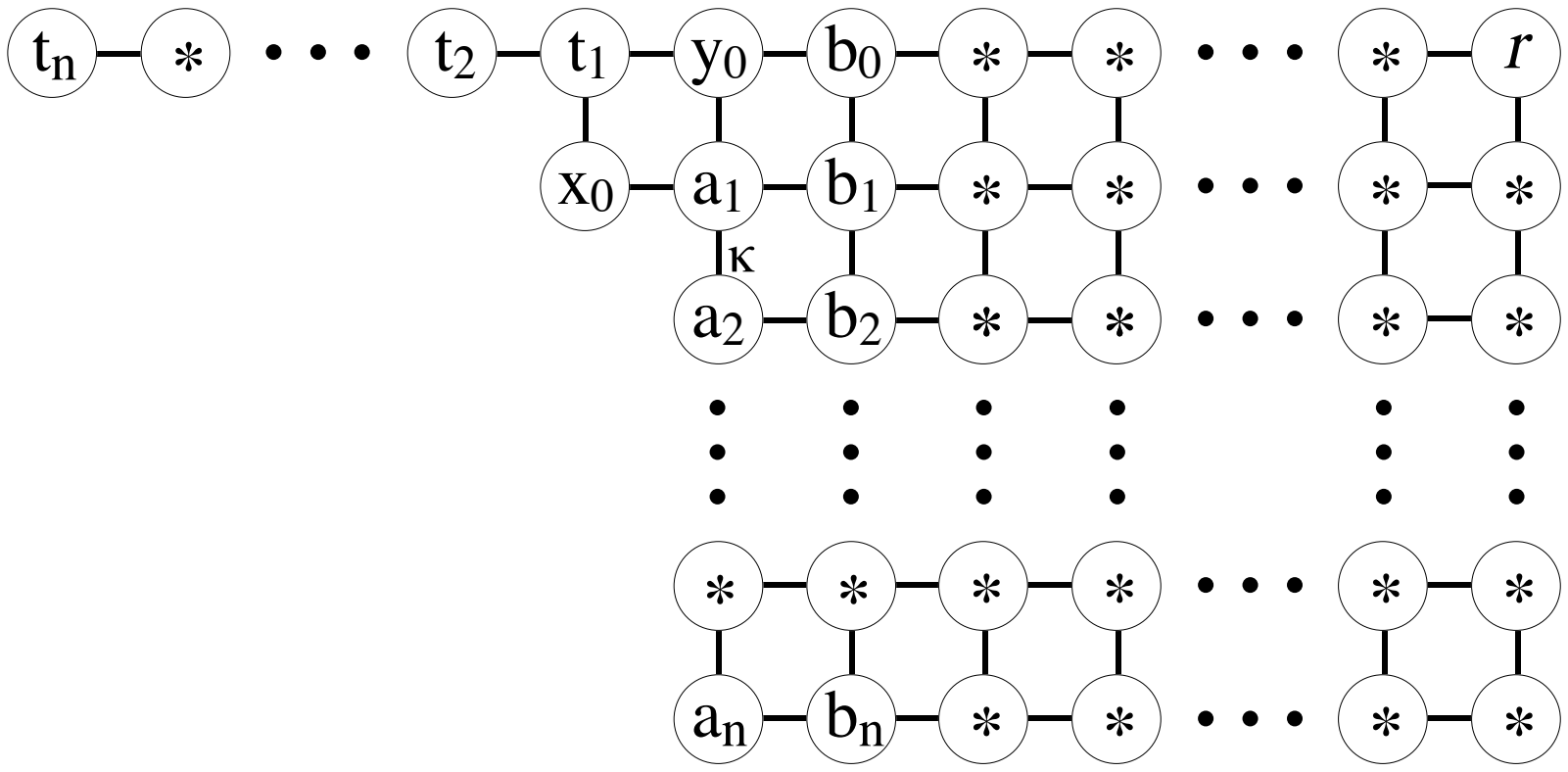

Figure 8.4

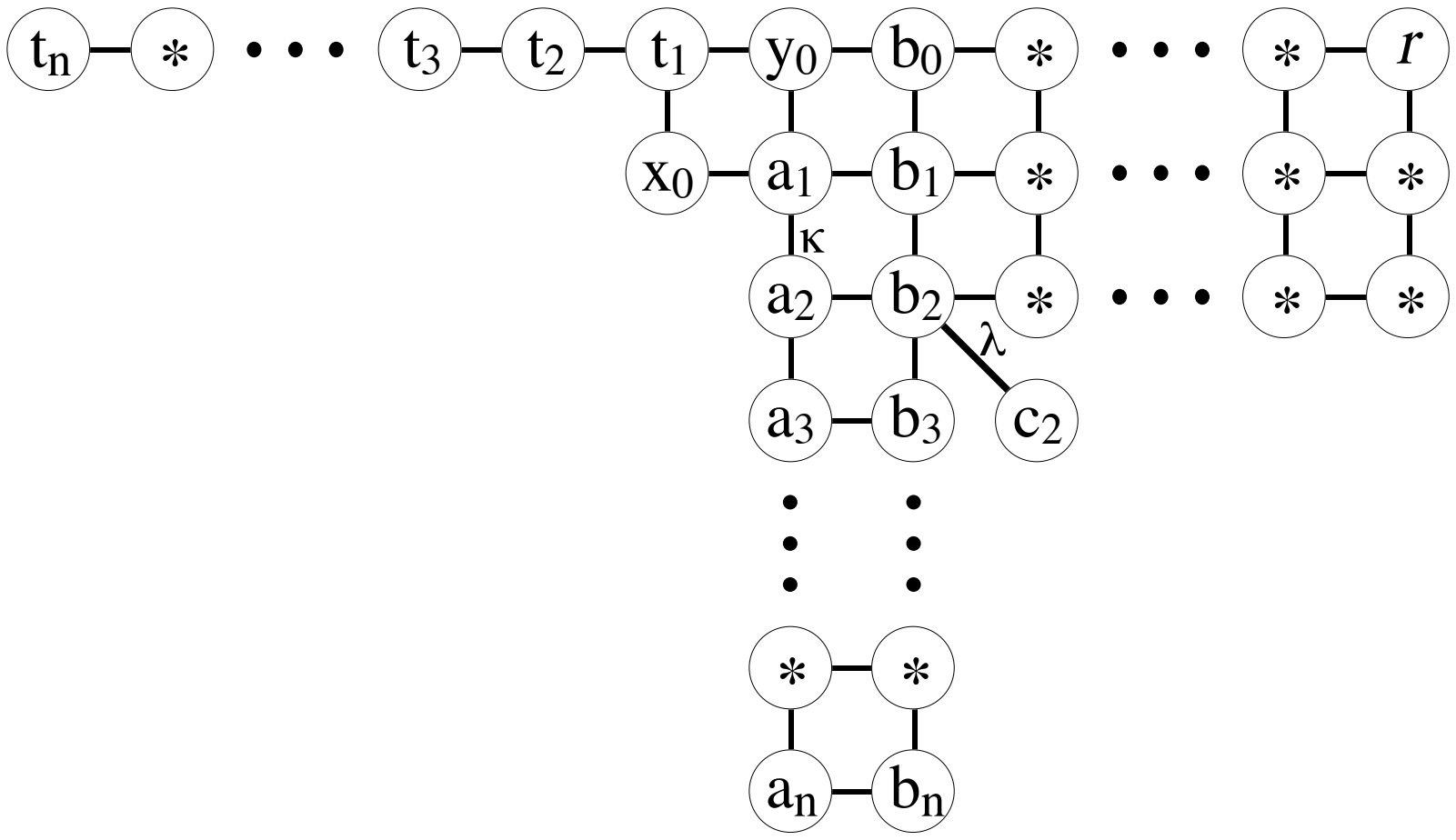

Figure 8.5

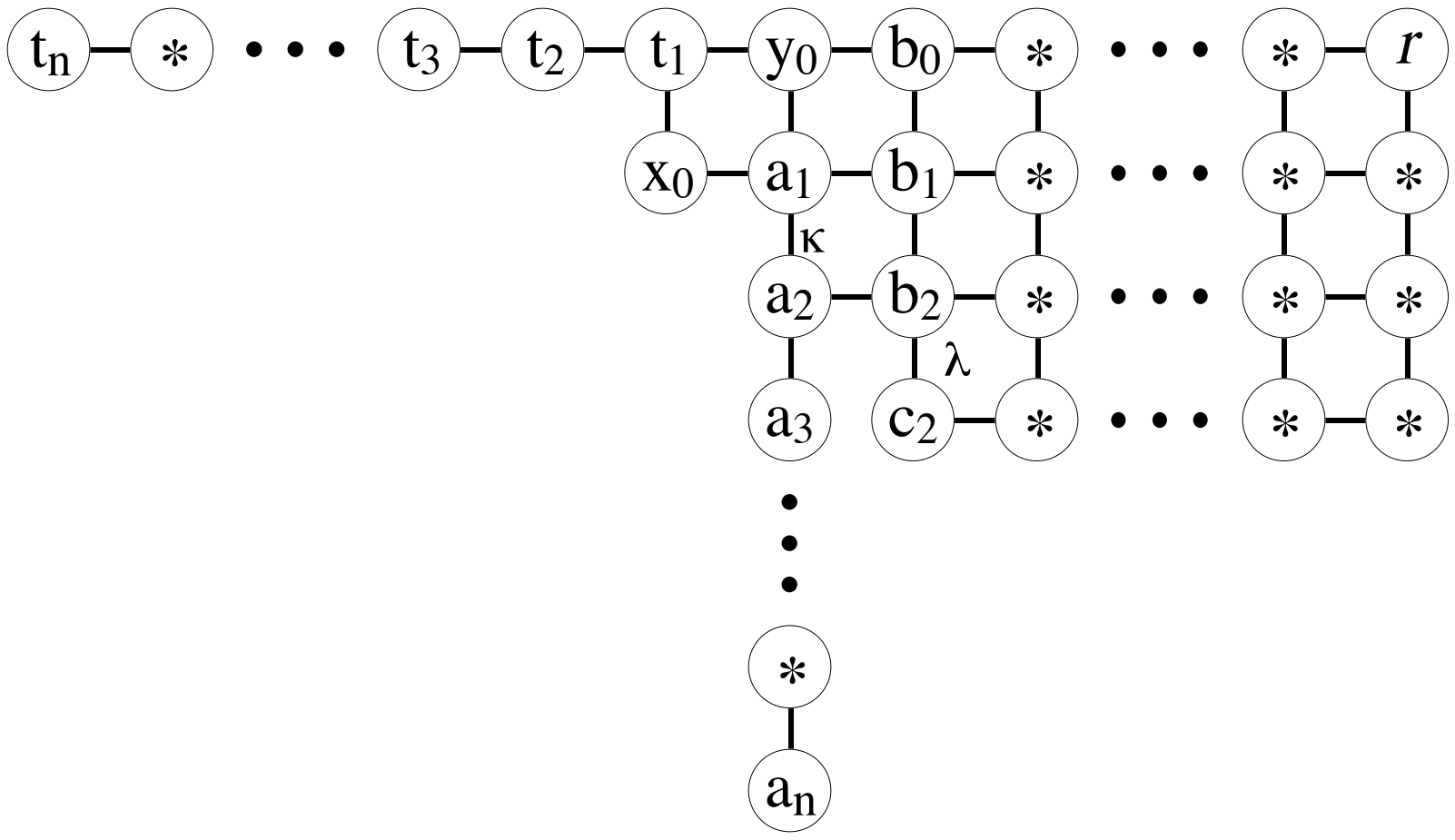

Figure 8.6

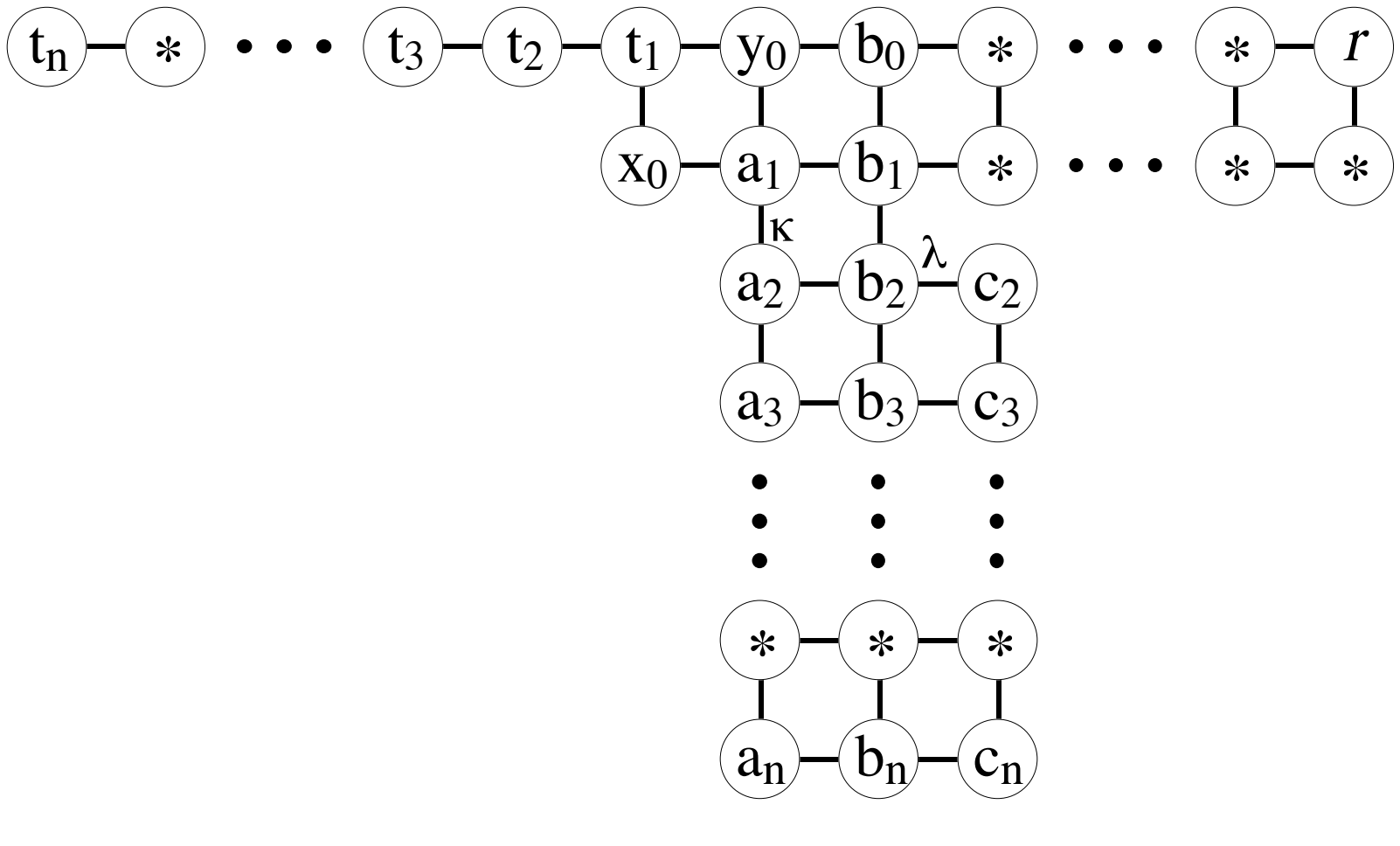

Figure 8.7

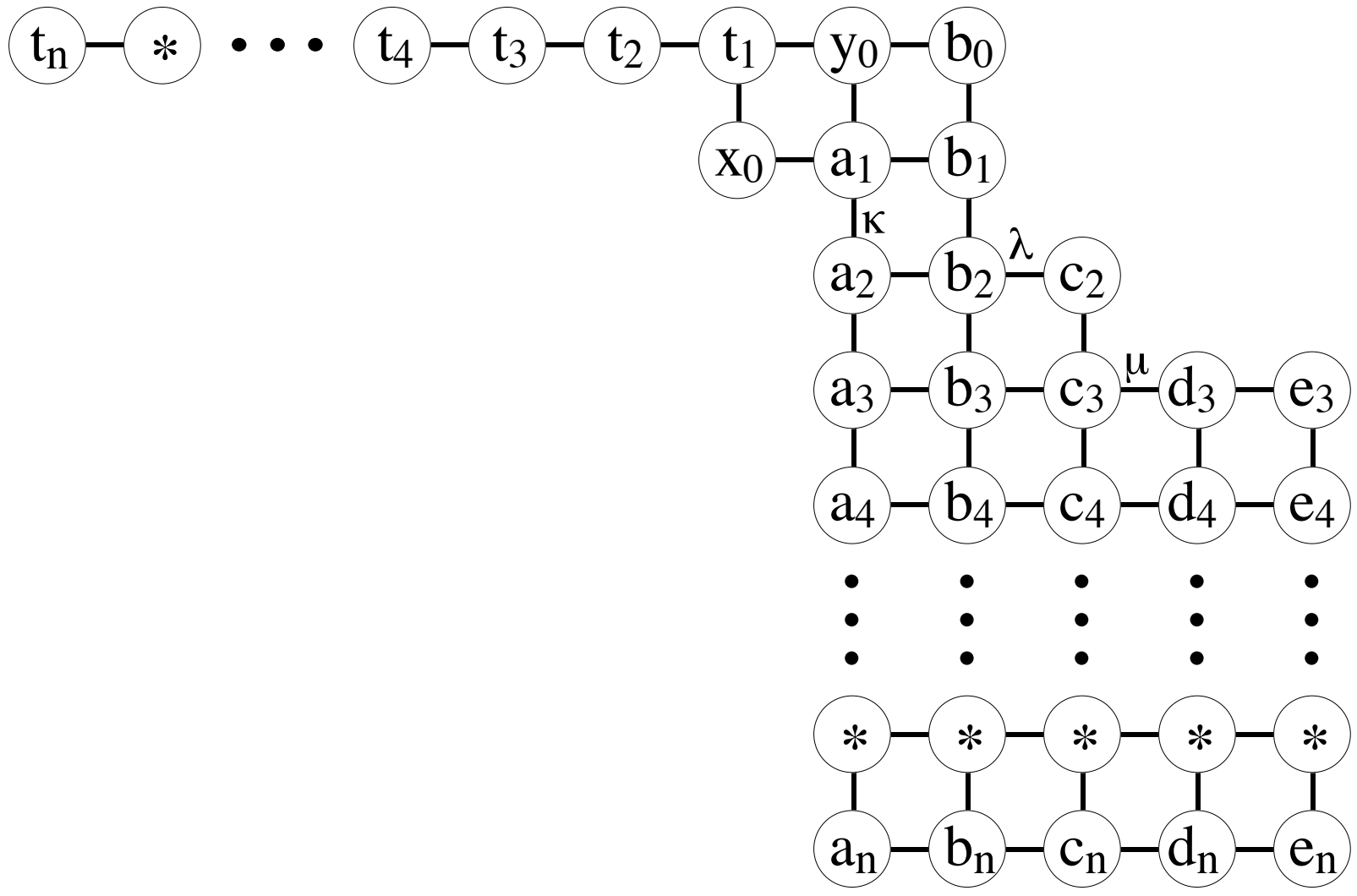

Figure 8.8

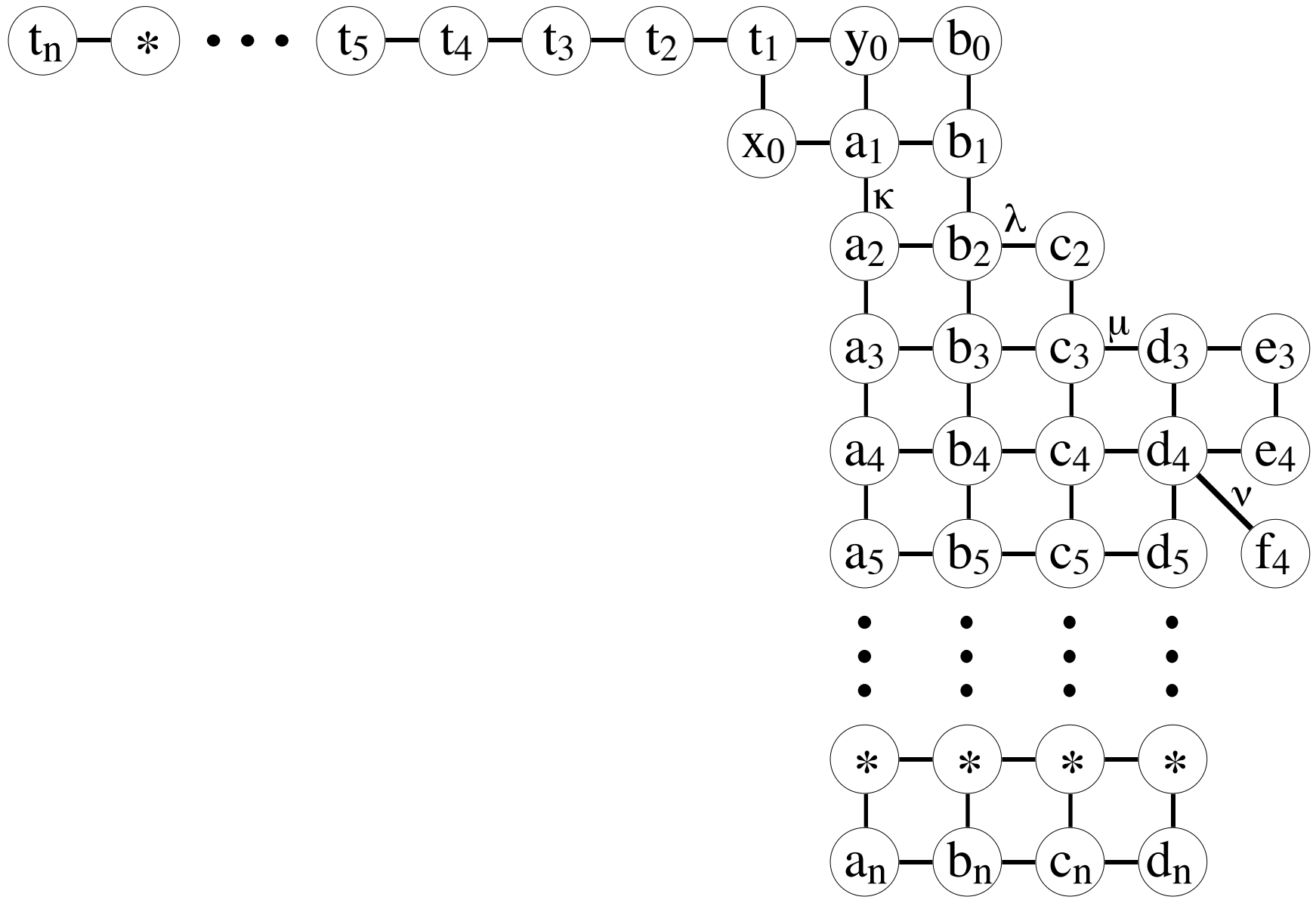

Figure 8.9

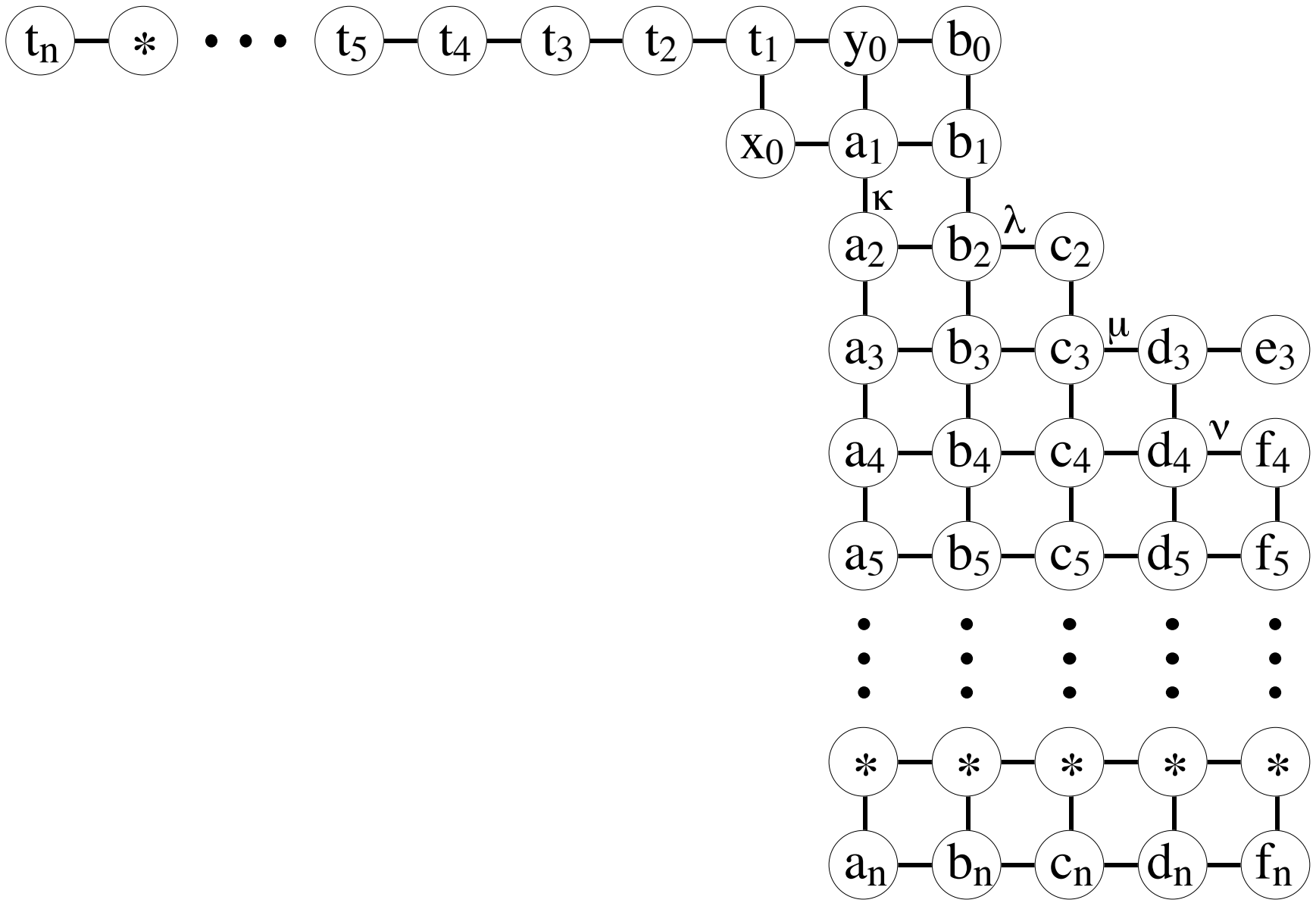

Figure 8.10

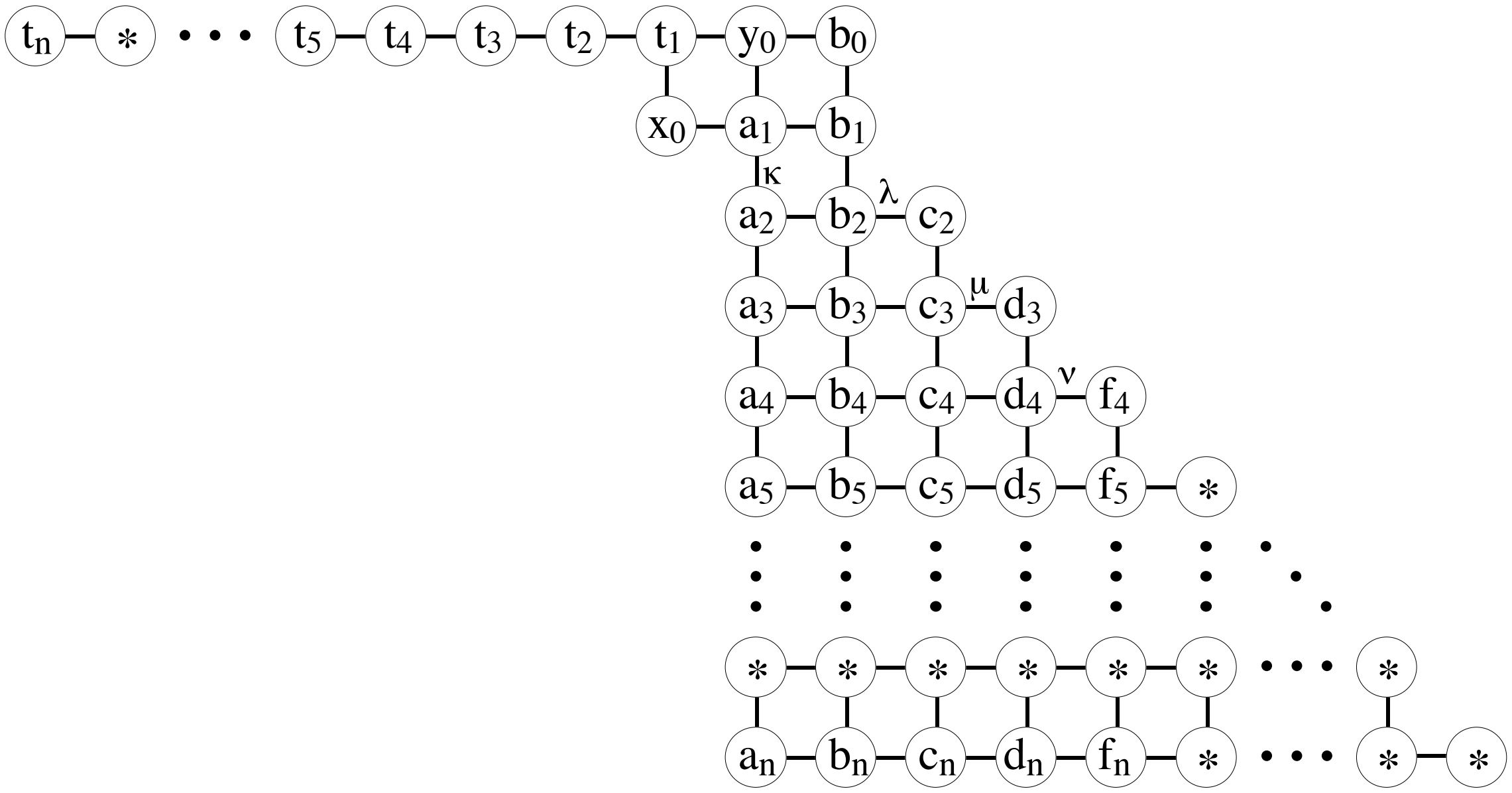

Figure 8.11

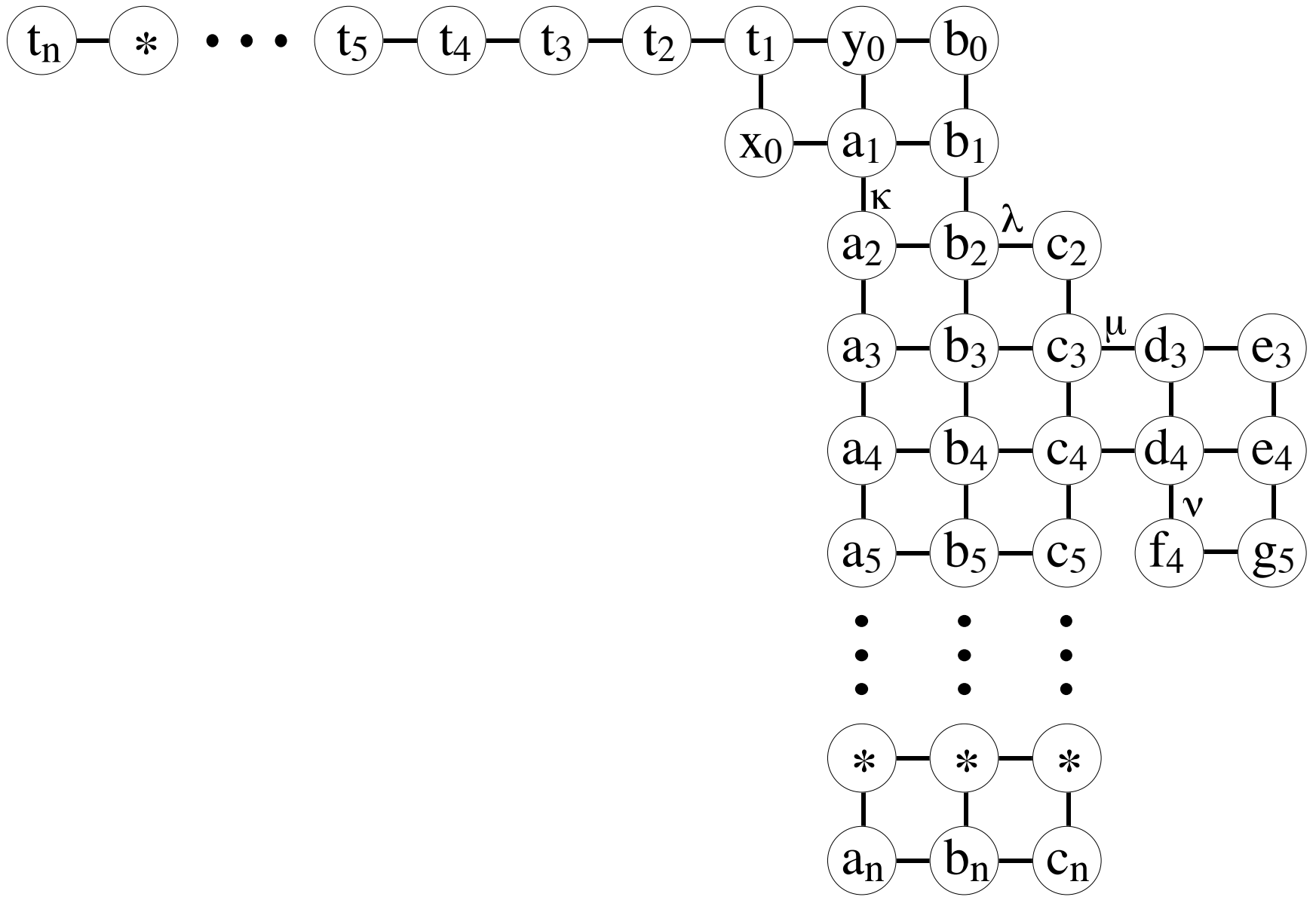

Figure 8.12

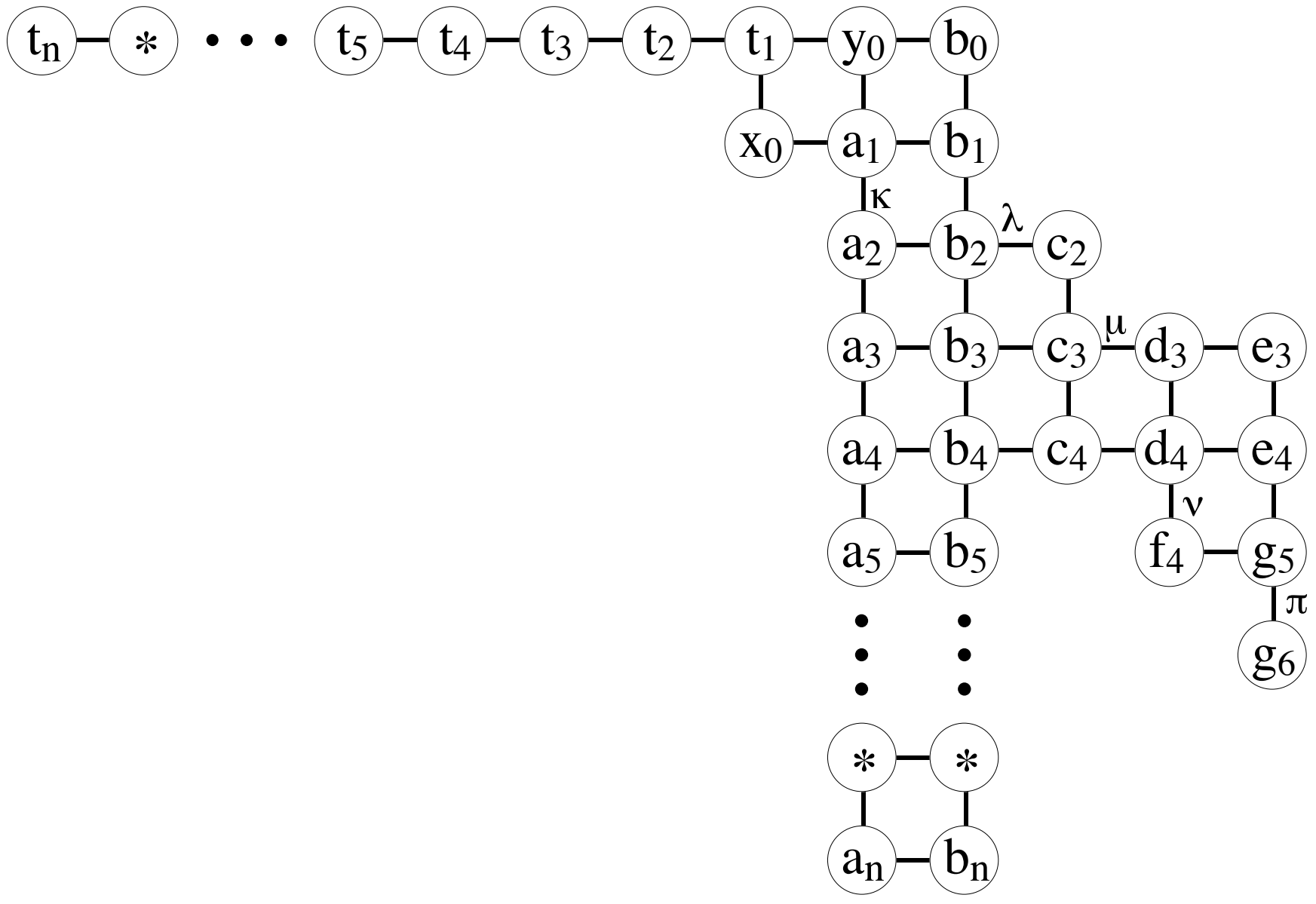

Figure 8.13

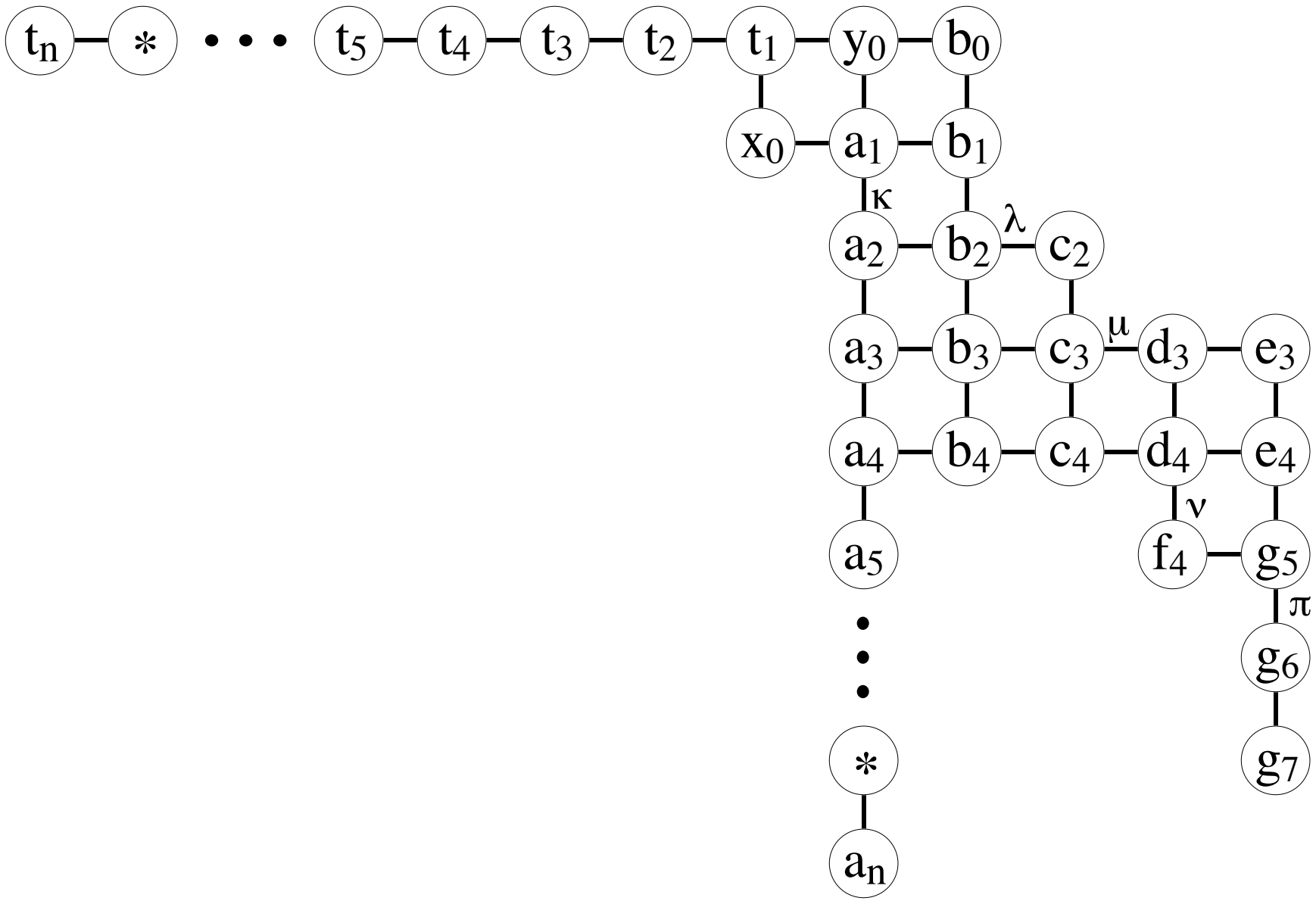

Figure 8.14

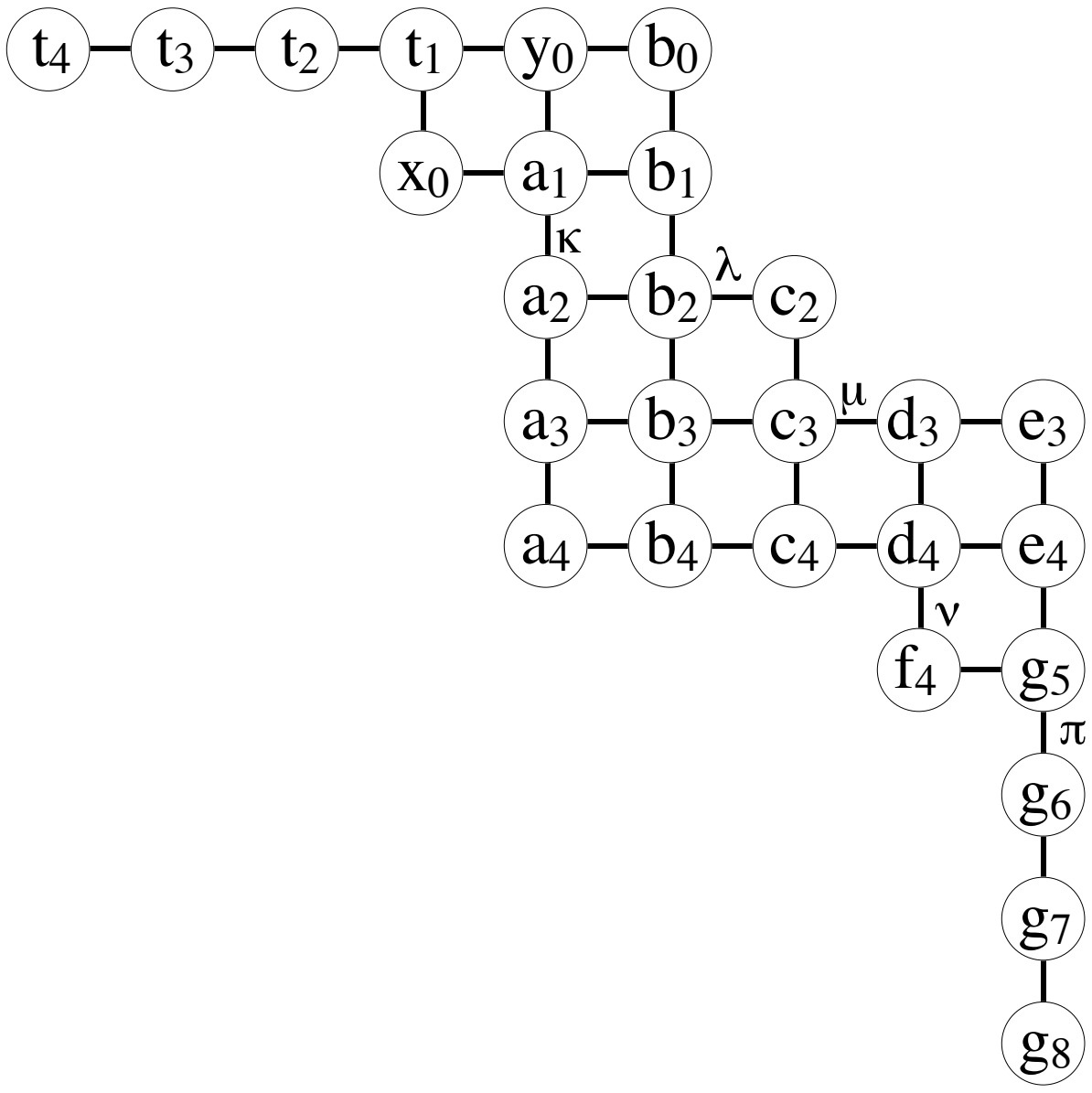

Figure 8.15